\documentclass[review, hidelinks, onefignum, onetabnum]{article}

\usepackage[top=20mm,bottom=20mm,left=20mm,right=20mm,headheight=24pt,headsep=15mm,footskip=14.5mm]{geometry}

\usepackage{amsmath, amssymb, amsthm}
\usepackage{xcolor}
\usepackage{xeCJK}
\usepackage{lipsum}
\usepackage{amsfonts}
\usepackage{graphicx}
\usepackage{algorithmic}
\usepackage{mathrsfs}
\usepackage[shortlabels]{enumitem} 

\usepackage{algorithm}
\usepackage{algorithmic}
\usepackage{hyperref} 

\newtheorem{definition}{Definition}[section]
\newtheorem{remark}{Remark}[section]

\definecolor{lightred}{RGB}{255,102,102} 
\definecolor{vermilion}{RGB}{197,39,45} 
\definecolor{eblue}{RGB}{1,126,218} 
\definecolor{teal}{RGB}{3,127,119} 
\definecolor{egreen}{RGB}{0,120,2} 
\definecolor{darkblue}{RGB}{26,0,217} 
\definecolor{yellowa}{RGB}{230,155,96} 
\definecolor{goldenorange}{RGB}{250,174,95} 

\definecolor{tealgreen}{RGB}{32, 178, 170}
\definecolor{brightorange}{RGB}{255, 140, 0}
\definecolor{crimsonred}{RGB}{220, 20, 60}
\definecolor{violetpurple}{RGB}{138, 43, 226}
\definecolor{dodgerblue}{RGB}{30, 144, 255}
\definecolor{darkgolden}{RGB}{229, 193, 0} 
\definecolor{deepgreen}{RGB}{0, 116, 0}
\definecolor{copperbrown}{RGB}{198, 122, 55}
\definecolor{emerald}{RGB}{58,191,153} 
\definecolor{cornflowerblue}{RGB}{24,129,244} 

\usepackage{version}
\includeversion{myversiona}\excludeversion{myversionb}

\usepackage{multicol,multirow}

\def\vx {{\boldsymbol u}}

\def\vw{{\boldsymbol w}}
\def\vx{{\boldsymbol x}}

\def\bA{\mathbf{A}}

\def\bD{\mathbf{D}}

\def\bI{\mathbf{I}}
\def\bX{\mathbf{X}}
\def\bY{\mathbf{Y}}

\def\bU{\mathbf{U}}
\def\bV{\mathbf{V}}
\def\bW{\mathbf{W}}
\def\bSig{\boldsymbol{\Sigma}}

\def\btheta{\boldsymbol{\theta}}

\def\sfE{\mathsf{E}}

\def\cA{\mathcal{A}}

\def\cD{\mathcal{D}}
\def\cE{\mathcal{E}}

\def\cG{\mathcal{G}}
\def\cJ{\mathcal{J}}
\def\cQ{\mathcal{Q}}

\def\cL{\mathcal{L}}
\def\cN{\mathcal{N}}
\def\cS{\mathcal{S}}
\def\cU{\mathcal{U}}
\def\cV{\mathcal{V}}

\def\cX{\mathcal{X}}

\def\ocA{\overline{\mathcal{A}}}

\def\ocL{\overline{\mathcal{L}}}

\def\ocU{\overline{\mathcal{U}}}

\def\ocX{\overline{\mathcal{X}}}

\def\E{\mathbb{E}}

\def\V{\mathbb{V}\mathrm{ar}}

\def\argmin{\mathop{\rm argmin}}
\def\argmax{\mathop{\rm argmax}}

\newcommand\norm[1]{\left\Vert#1\right\Vert}
\newcommand\abs[1]{\left|#1\right|}
\newcommand\trace[1]{\mathrm{Tr}\left( #1\right)}

\title{Variational Bayesian Inference for 
Tensor Robust Principal Component Analysis
\thanks{Corresponding author: Youwei Wen.
{This work was funded by the National Natural Science Foundation of China No. 12361089, 12571564, Guangdong Basic and Applied Research Foundation 2024A1515012347, HKRGC Grants No. LU13300125,
LU11309922, ITF Grant No. MHP/054/22, and LU BGR 105824. }
}}

\author{Chao Wang\thanks{Department of Statistics and Data Science,   Southern University of Science and Technology,  Shenzhen 518005,  Guangdong Province,  China
  ({wangc6@sustech.edu.cn}).}
  \and Huiwen Zheng\footnotemark[2]
\and Raymond Chan\thanks{Lingnan University,  Hong Kong SAR,  China 
  ({raymond.chan@ln.edu.hk}).}
\and Youwei Wen\thanks{Key Laboratory of Computing and Stochastic Mathematics (LCSM), School of Mathematics and Statistics, Hunan Normal University, Changsha, Hunan, China.  
  ({wenyouwei@gmail.com})} 
}

\usepackage{amsopn}
\DeclareMathOperator{\diag}{diag}

\begin{document}

\maketitle

\begin{abstract}
 Tensor Robust Principal Component Analysis (TRPCA)
holds a crucial position in machine learning and computer 
vision. It aims to recover underlying low-rank structures 
 and to characterize the sparse structures of noise. Current approaches often encounter difficulties in accurately capturing the low-rank properties of tensors and balancing the trade-off between low-rank and sparse components, especially in a mixed-noise scenario. To address these challenges,  we introduce a Bayesian framework for TRPCA,  which integrates a low-rank tensor nuclear norm prior and a generalized sparsity-inducing prior. 
By embedding the priors within the Bayesian framework, our method can automatically determine the optimal tensor nuclear norm and achieve a balance between the nuclear norm and sparse components. 
Furthermore,  our method can be efficiently extended to the weighted tensor nuclear norm model. 
Experiments conducted on synthetic and real-world datasets demonstrate the effectiveness and superiority of our method compared to state-of-the-art approaches.


\end{abstract}

%
%

\section{Introduction}
With data becoming ubiquitous from diverse fields and applications,  data structures are becoming increasingly complex with higher dimensions.
Tensor, a multidimensional array, is an efficient data structure with broad applications, including machine learning \cite{panagakis2021tensor} and computer vision \cite{  phan2020stable}. Meanwhile, high-dimensional data always lie near a low-dimensional manifold, which can be interpreted by their low rank. In matrix processing, the low-rank assumption allows two-dimensional data recovery from incomplete or corrupted data \cite{candes2012exact}. However, expanding the low-rank concept from matrices to tensors remains an unresolved challenge. A main challenge in tensor analysis is that the tensor rank is not well defined. Various definitions of tensor rank have been proposed. For example, the CANDECOMP/PARAFAC (CP) rank, as described in \cite{kolda2009tensor}, is based on the CP decomposition \cite{kiers2000towards} and identifies the smallest number of rank-one tensors needed to represent a tensor. The Tucker rank \cite{gandy2011tensor}, which stems from the Tucker decomposition \cite{tucker1966some}, consists of a vector where each component corresponds to the rank of a matrix obtained by unfolding the original tensor. Furthermore, developments in tensor singular value decomposition (t-SVD) \cite{kilmer2011factorization} have led to the tensor multi-rank \cite{gandy2011tensor} and tubal rank \cite{kilmer2013third}, both of which are analogous to the matrix singular value decomposition (SVD).

Among all these tensor applications, exploring low-rank features in sparse tensor decomposition has become essential, which is called 
Tensor Robust Principal Component Analysis (TRPCA) \cite{lu2019tensor}. It extends
the Robust Principal Component Analysis (RPCA) \cite{jolliffe2016principal} from matrices to tensors. Specifically, TRPCA seeks to extract the low tubal rank component, \( \cL \), and eliminate the noise component, \( \cS \), derived from noisy observations, $\cX$, expressed as \( \cX = \cL + \cS \). 
This is achieved through the optimization process \cite{lu2019tensor,YANG2022108311,mu2014square,goldfarb2014robust,wang2023guaranteed} described as 
\begin{equation}\label{equ:TRPCAfirst}
\min_{\cX=\cL+\cS} \norm{\cL}_*+\lambda\norm{\cS}_1.
\end{equation}
where  $\|\mathcal{L}\|_*$ is the tensor nuclear norm as the convex relaxation to a certain tensor rank. Note that minimizing the rank is an NP-hard problem. Various approximations have been proposed to approach different tensor ranks \cite{jiang2020multi,zheng2024scale,mu2020weighted}. 
Here,  $\|\mathcal{S}\|_1$ is the $\ell_1$ norm of sparsity,  and $\lambda> 0$ is the parameter used to balance low-rankedness and sparsity.

In the TRPCA model,  we can further reformulate the equality constraint by a penalty term and turn the optimization model \eqref{equ:TRPCAfirst} into
\begin{equation}\label{HSImodel}
\min_{\cS, \cL}\tfrac{\theta_1}{2}\|\mathcal{X}-\mathcal{S}-\mathcal{L}\|_F^2+\theta_2 \|\mathcal{S}\|_1+\theta_3 \|\mathcal{L}\|_*, 
\end{equation}
where $\theta_1, \theta_2$ and $\theta_3$ are tuning parameters. 
Note that \eqref{HSImodel} has broadened applications by 
assuming observation data is constructed not just by low-rank tensor and sparsity but also with certain bias or Gaussian noise, i.e., 
\begin{equation}\label{eq:3compon}
\mathcal{X}=\mathcal{L}+\mathcal{S}+\mathcal{E},  
\end{equation}
where  $\mathcal{E}$ is the corresponding bias and the Gaussian noise. This model is widely used in mixed noise removal \cite{zheng2019mixed,zhou2019bayesian} and hyperspectral denoising \cite{su2022hyperspectral}. 



The selection of the parameters in the model \eqref{equ:TRPCAfirst} and \eqref{HSImodel} is critical. 
Under the t-SVD framework, the optimal parameter for $\lambda$ in \eqref{equ:TRPCAfirst}  is suggested in \cite{lu2019tensor} for the tensor nuclear norm. 
Nevertheless, it cannot be extended to other forms of tensor low-rank regularization, such as the weighted tensor nuclear norm.
This issue becomes more serious when dealing with models involving multiple parameters in \eqref{HSImodel}. 
Traditional parameter selection methods, including the discrepancy principle \cite{Morozov}, $L$-curve \cite{Hansen1993}, GCV \cite{golub1979}, and RWP \cite{2013parameter,2023Poisson-whiteness}, are often customized to specific regularization formulations and need iterative minimizations, which makes it inadequate for our tensor recovery problem in \eqref{HSImodel}.

In this paper, we address the intricate task of simultaneously estimating tensors $\mathcal{L}$ and $\mathcal{S}$ and their regularization parameters $\theta_i$ for TRPCA. 
We introduce variational Bayesian inference (VBI) \cite{chantas2009variational} as a powerful tool to tackle this challenge, reformulating the optimization problem within a Bayesian framework. 
By treating regularization parameters $\theta_i$ as hyperparameters, we apply the inherent strengths of Bayesian methods, popular for their success in inverse problems \cite{uribe2023,cai2024a,jia2021,helin2022,glaubitz2024,zhou2024c,law2022} and established applications in matrix and tensor problems like matrix completion \cite{yang2018a}, tensor completion \cite{bazerque2013,tong2023}, and low-rank tensor approximation \cite{menzen2022}.

Despite these successes, the adoption of VBI in TRPCA remains limited.
To our best knowledge, only \cite{zhou2019bayesian} has explored VBI for TRPCA, employing a generalized sparsity-inducing prior.  
However, this method directly expresses the low-rank tensor as a t-product of two smaller factor tensors, presupposing the tubal rank, and models the sparse component $\mathcal{S}$ with independent Gaussian priors, which may not be optimal for sparse data.
In contrast, we propose an approach that employs a tensor nuclear norm prior to $\mathcal{L}$, eliminating the need for predefined tensor ranks. For the sparse component $\mathcal{S}$, we adopt a Laplace prior, which better captures sparse structures. 
This reformulation enhances model flexibility, offering a more principled and less restrictive approach to tensor recovery, thereby mitigating limitations posed by prior assumptions on tensor ranks or sparsity patterns.

In comparison, joint maximum a posteriori (MAP) estimation minimizes the negative log posterior to obtain point estimates for $\mathcal{S}$, $\mathcal{L}$, and $\btheta$, simultaneously recovering tensors and parameters. 
Our VBI framework, however, approximates the full posterior distribution, enabling uncertainty quantification alongside point estimates. 
For practical applications such as denoising and background subtraction, we use the expectation of the variational distribution as the point estimate for $\mathcal{S}$ and $\mathcal{L}$,
offering a robust and versatile approach to tensor recovery.

The primary contributions of this work are succinctly summarized as:
\begin{enumerate}[label=({\arabic*)},leftmargin=4ex, itemindent=0ex]
\item Innovative Variational Bayesian Tensor Recovery Model:
   This paper proposes a novel variational Bayesian inference model for tensor recovery. It characterizes low-rank tensors using the tensor nuclear norm and sparse tensors via the Laplacian distribution. This approach enables simultaneous inference of both low-rank and sparse components along with the hyperparameters (regularization parameters), eliminating the need for pre-specifying the tensor rank. 
\item Efficient Inference via Laplacian Approximation and MM Framework:
We introduce a Laplacian approximation methodology to tackle the computational intricacies associated with non-Gaussian posteriors arising from Laplace priors imposed on sparse tensor $\cS$ and low-rank tensor $\cL$. This approach directly tackles the $\ell_1$ norm minimization and tensor nuclear norm minimization problems in estimating the expectations of sparse tensor $\cS$ and low-rank tensor $\cL$. For covariance matrix computation, it integrates with the Majorization-Minimization (MM) framework, deriving a tight lower bound for the non-quadratic distributions encountered in the $\ell_1$ norm and tensor nuclear norm. This facilitates efficient variance computations, thereby significantly enhancing the efficiency and accuracy of inferring low-rank, sparse tensors as well as their hyperparameters.

\end{enumerate}


The rest of this paper is organized as follows. In Section \ref{sectionprepara},  we introduce the main preliminaries, including tensors and their decomposition. In Section \ref{Sectionbaypara},  we describe the hierarchical Bayesian model,  joint density, and hyperprior. In Section \ref{sectionvarpara},  we apply variational Bayesian inference to infer hyperparameters $\theta_i$ and solve the tensors $\cL, \cS$ at the same time.  In Section \ref{sectionexppara},  we provide the experimental results and show the superiority of our proposed methods. Finally, in Section \ref{sec:concl}, some conclusions are drawn.

\section{Preliminaries}\label{sectionprepara}
This section provides an overview of fundamental notations and definitions that will be utilized throughout the paper.

\subsection{Notations}

The set of natural numbers is denoted by $\mathbb{N}$, the set of real numbers by $\mathbb{R}$, and the set of complex numbers by $\mathbb{C}$. 
In the context of tensors, we adopt the convention of using boldface Euler script letters, exemplified by $\mathcal{A}$, to represent them. 
Matrices, on the other hand, are indicated with boldface uppercase letters, such as $\mathbf{A}$, with the identity matrix specifically denoted by $\mathbf{I}$. 
Vectors follow the convention of being written in boldface lowercase letters, like $\mathbf{a}$, whereas single values or scalars are represented by regular lowercase letters, for instance, $a$.
Regarding indexing, for a vector $\mathbf{a}$, the $i$-th element is denoted by $\mathbf{a}_i$. 
For a matrix $\mathbf{A}$, $\mathbf{A}_{i:}$ signifies the $i$-th row, $\mathbf{A}_{:j}$ denotes the $j$-th column, and the specific element located at the intersection of the $i$-th row and $j$-th column is represented by either $a_{ij}$ or, more commonly in matrix notation, $\mathbf{A}_{ij}$.
When dealing with a third-order tensor $\mathcal{A}$, each element positioned at the intersection of the $i$-th, $j$-th, and $k$-th dimensions is denoted by $a_{ijk}$ or, more conventionally for tensors, $\mathcal{A}_{ijk}$. This tensor can be dissected into distinct structural components: column fibers are designated as $\mathcal{A}_{:jk}$, row fibers as $\mathcal{A}_{i:k}$, and tube fibers as $\mathcal{A}_{ij:}$. Furthermore, the tensor can be analyzed through various slices: horizontal slices are noted as $\mathcal{A}_{i::}$, lateral slices as $\mathcal{A}_{:j:}$, and frontal slices as $\mathcal{A}_{::k}$.

We define the inner product of matrices $\mathbf{A}$ and $\mathbf{B}$ as $\langle \mathbf{A}, \mathbf{B}\rangle := \operatorname{Tr}(\mathbf{A}^* \mathbf{B})$, where $\mathbf{A}^*$ is the conjugate transpose of $\mathbf{A},$ and $\operatorname{Tr}(\cdot)$ represents the trace of a matrix. If $\mathbf{A}$ consists only of real numbers, $\mathbf{A}^{T}$ denotes its transpose.
 The $\ell_{2}$-norm of a vector $\mathbf{v}$ in the complex number space $\mathbb{C}^{n}$ is defined by $\|\mathbf{v}\|_{2} = (\sum_{i} |\mathbf{v}_i|^2)^{1/2}$, measuring the vector's magnitude in Euclidean space.

The inner product between two tensors ${\mathcal{A}}$ and ${\mathcal{B}}$ in ${\mathbb{C}^{n_{1} \times n_{2} \times n_{3}}}$ is defined as ${\langle\mathcal{A},  \mathcal{B}\rangle=\sum_{k=1}^{n_{3}}\left\langle \mathcal{A}_{::k},  \mathcal{B}_{::k}\right\rangle}$. The complex conjugate of ${\mathcal{A}}$,  which takes the complex conjugate of each entry of ${\mathcal{A}}$,  is denoted as conj$(\mathcal{A})$. 
The conjugate transpose of a tensor ${\mathcal{A}} \in\mathbb{C}^{n_1 \times n_2 \times n_{3}}$ is a tensor $\mathcal{A}^\ast$ obtained by conjugate transposing each of the frontal slices and then reversing the order of transposed frontal slices 2 through $n_3$.
The tensor ${\ell_{1}}$-norm of $\mathcal A$ is defined  as ${\|\mathcal{A}\|_{1}=\sum_{i j k}\left|a_{i j k}\right|}$,  
and the Frobenius norm as ${\|\mathcal{A}\|_{F}=}$ ${\sqrt{\sum_{i j l}\left|a_{i j k}\right|^{2}}}$.

\subsection{T-product and t-SVD}
Before introducing the definitions, we define three operators: 
\begin{equation}
 \operatorname{unfold}(\mathcal{A})=\left[\begin{array}{c}
\cA_{::1} \\
\cA_{::2} \\
\vdots \\
\cA_{::n_3}
\end{array}\right],  \text { fold }(\operatorname{unfold}(\mathcal{A}))=\mathcal{A},   
\end{equation}
and 
$$\operatorname{bcirc}(\mathcal{A}):=\left[\begin{array}{cccc}
\cA_{::1} & \cA_{::n_3} & \cdots & \cA_{::2} \\
\cA_{::2} & \cA_{::1} & \cdots & \cA_{::3} \\
\vdots & \vdots & \ddots & \vdots \\
\cA_{::n_3} & \cA_{::, n_3-1} & \cdots & \cA_{::1}
\end{array}\right] \in \mathbb{R}^{n_{1} n_{3} \times n_{2} n_{3}}. $$
Here unfold$(\cdot)$ maps ${\mathcal{A}}$ to a matrix of size ${n_{1} n_{3} \times n_{2}}$ and fold$(\cdot)$ is its inverse operator. 
We introduce the notation $\mathbf{A} := \operatorname{bdiag}(\mathcal{A})$ to concisely represent the block diagonal matrix derived from the tensor $\mathcal{A}$. Here, $\operatorname{bdiag}(\cdot)$ designates the block diagonalization operator, with the $i$-th block corresponding to $\mathcal{A}_{::i}$.

Now, we focus on applying the Discrete Fourier Transformation (DFT) to tensors. We represent the tensor $\mathcal{A}$ transformed by DFT along its third (tubal) dimension as $\overline{\mathcal{A}}$. This transformation is executed using the MATLAB command $\operatorname{fft}$, specifically performed as $\overline{\mathcal{A}} = \operatorname{fft}(\mathcal{A}, [], 3)$. Conversely, to revert the tensor to its original form from $\overline{\mathcal{A}}$, we use the inverse operation with $\mathcal{A} = \operatorname{ifft}(\overline{\mathcal{A}}, [], 3)$. 
We also introduce the notation $\overline{\mathbf{A}} := \operatorname{bdiag}(\overline{\mathcal{A}})$ to represent the block diagonal matrix constructed from the tensor $\overline{\mathcal{A}}$. 
Next, we introduce the definition of t-product.

\begin{definition}
\label{def:tprod}
{\rm (t-product\textup{ \cite{kilmer2011factorization}})}. Let ${\mathcal{A} \in \mathbb{R}^{n_{1} \times l \times n_{3}}}$ and ${\mathcal{B} \in \mathbb{R}^{l \times n_{2} \times n_{3}}}$,  then the t-product ${\mathcal{A} * \mathcal{B}}$ is defined by 
\begin{equation}
{ \mathcal{A} * \mathcal{B}=\operatorname{fold}(\operatorname{bcirc}(\mathcal{A}) \cdot \operatorname{unfold}(\mathcal{B})),  }
\end{equation}
resulting a tensor of size ${n_{1} \times n_{2} \times n_{3}}$. Note that $\mathcal{A} * \mathcal{B}= \mathcal{Z} $ if and only if $\overline{\mathbf{A}}\, \overline{\mathbf{B}}=\overline{\mathbf{Z}}.$ 
\end{definition}

Using the t-product framework, we define the identity tensor $\mathcal{I}\in\mathbb{R}^{n \times n \times n_3}$ as a tensor with its first frontal slice being the $n \times n$ identity matrix, while all subsequent frontal slices consist entirely of zeros. It is clear that ${\mathcal{A} * \mathcal{I}=\mathcal{A}}$ and ${\mathcal{I} * \mathcal{A}=\mathcal{A}}$ given the appropriate dimensions. 
In addition, a tensor ${\mathcal{H} \in}$ ${\mathbb{R}^{n \times n \times n_{3}}}$ is orthogonal if it satisfies 
$
\mathcal{H}^{*}*\mathcal{H}=\mathcal{H}*\mathcal{H}^{*}=\mathcal{I}.  
$
Moreover, we call a tensor $f$-diagonal if each of its frontal slices is a diagonal matrix.
Next, we define the tensor singular value decomposition as below:





\begin{definition} \label{skinnysvd}
 {\rm (tensor singular value decomposition: t-SVD \cite{kilmer2011factorization}).} 
 The t-SVD of $\mathcal{A}\in \mathbb{R}^{n_{1} \times n_{2} \times n_{3}} $ is given by 
\begin{equation}
{ \mathcal{A}=\mathcal{U} * \mathcal{D} * \mathcal{V}^{*} }, 
\end{equation} 
where ${\mathcal{U} \in \mathbb{R}^{n_{1} \times n_{1} \times n_{3}},  \mathcal{V} \in \mathbb{R}^{n_{2} \times n_{2} \times n_{3}}}$ are orthogonal tensors,  and ${\mathcal{D} \in \mathbb{R}^{n_{1} \times n_{2} \times n_{3}}}$ is an ${f}$-diagonal tensor. 
\end{definition}
It follows from Definition \ref{def:tprod} that $\mathcal{A}=\mathcal{U} * \mathcal{D} * \mathcal{V}^{*}$  if and only if $\overline{\mathbf{A}} =\overline{\mathbf{U}}\, \overline{\mathbf{D}}\, \overline{\mathbf{V}}^\ast.$
For tensor ${\mathcal{A} \in \mathbb{R}^{n_{1} \times n_{2} \times n_{3}}}$ with tubal rank $r$,  we also have skinny t-SVD similar as matrix. Let $r$ is the tubal rank of $\mathcal{A}$,  the skinny t-SVD of $\mathcal{A}$ is 
$\mathcal{A}=\mathcal{U} * \mathcal{D} * \mathcal{V}^{*}$,  where $\mathcal{U} \in \mathbb{R}^{n_{1} \times r \times n_{3}},  \mathcal{D} \in \mathbb{R}^{r \times r \times n_{3}},  \mathcal{V} \in \mathbb{R}^{n_{2} \times r \times n_{3}}$,  in which $\mathcal{U}^**\mathcal{U}=\mathcal{I}$ and $\mathcal{V}^**\mathcal{V}=\mathcal{I}$.

\begin{definition}
\label{def:tar}
    {\rm (tensor average rank and tubal rank \textup{\cite{lu2019tensor}})} 
    The  tensor average rank of ${\mathcal{A} \in \mathbb{R}^{n_{1} \times n_{2} \times n_{3}}}$, denoted as $\operatorname{rank}_{a}(\mathcal{A})$,  is defined as 
    \begin{equation*}
        \operatorname{rank}_{a}(\mathcal{A}) = \frac{1}{n_3} \operatorname{rank}(\operatorname{bcirc}(\mathcal{A})) = \frac{1}{n_3} \sum_{i=1}^{n_3}  \operatorname{rank} (\overline{\mathbf{A}}^{(i)}). 
    \end{equation*}
\end{definition}
The tensor tubal rank, denoted as $\operatorname{rank}_{t}(\mathcal{A})$, is defined as the number of nonzero singular tubes of $\mathcal{S}$, where $\mathcal{S}$ comes from the t-SVD of $\mathcal{A}$, i.e. $\mathcal{A}=\mathcal{U} * \mathcal{S} * \mathcal{V}^{*}$. In other words, one has
$$\operatorname{rank}_{t}(\mathcal{A})=\#\{i, \mathcal{S}(i,i,:)\neq0\}.$$

For tensor ${\mathcal{A} \in \mathbb{R}^{n_{1} \times n_{2} \times n_{3}}}$ with tubal rank $r$,  we also have skinny t-SVD similar as matrix. 
Minimizing the tubal rank is an NP-hard problem; we introduce a tensor nuclear norm as a convex relaxation. 

\begin{definition} \label{def:tnn}
{\rm (tensor nuclear norm \textup{\cite{lu2019tensor}}).} 
Let $\mathcal{A}=\mathcal{U}*\mathcal{D}*\mathcal{V}^{*}$ be the t-SVD of $\mathcal{A}\in \mathbb{R}^{n_{1} \times n_{2} \times n_{3}}$. Define $\sigma_{jk}(\mathcal{A})$ is the $j$-th  singular value of $\ocA_{::k}$,  or simply $\sigma_{jk}$ if the context is clear. The tensor nuclear norm (TNN)  of $\mathcal{A}$ is defined as
\[
\norm{\cA}_*=\frac{1}{n_{3}}\sum_{k=1}^{n_{3}}\norm{\ocA_{::k}}_*=\frac{1}{n_{3}}\sum_{k=1}^{n_{3}}\sum_{j=1}^{\min(n_{1},  n_{2}) } \sigma_{jk}. 
\]
\end{definition}

\subsection{Probability distribution}
Here, we define three kinds of probability distribution:
the uniform distribution, the Gamma distribution, and the multivariate Gaussian distribution.

The uniform distribution is a distribution that assigns equal probability mass to a region. 
For $a,b\in R$ and $a<b$, the uniform distribution for a random variable $x\in R$ is defined as
\[
p(x)=\begin{cases}
    \frac{1}{b-a}, & a\leq x\leq b,\\
    0, & \mathrm{otherwise}.
\end{cases} 
\]

The Gamma density function is  given as
\begin{equation}\label{gammafunpdf}
p(x)=\cG(x|a, b)\propto x^{a-1}\exp(-bx ), 
\end{equation}
where $a > 0$ and $b > 0$ represent shape and scale parameters respectively. We have its mean and variance of these Gamma distributions:
\begin{equation}\label{C6}
{\E(x)}=\frac{a }{b }, \quad {\V(x)}=\frac{a }{b^2}.
\end{equation}

The multivariate Gaussian distribution is fully characterized by a mean vector $\mu$ and a covariance matrix $\bSig$ and is defined as
\[
p(\vx)=\cN(\vx|\mu,\bSig)=(2\pi)^{-\tfrac{n}{2}}|\bSig|^{-\tfrac{1}{2}}\exp\left(-\tfrac{1}{2}\norm{\vx-\mu}_{\bSig^{-1}}^2 \right),
\]
where $\vx\in \mathbb{R}^n$ is a random variable.


\section{Bayesian model}\label{Sectionbaypara}

\subsection{The likelihood}

In \eqref{eq:3compon}, we assume the observed data $\cX$ can be decomposed into three parts: $\cE, \cS, \cL$. Note that 
$\cE$'s elements as independent and identically distributed (i.i.d.) from a zero-mean normal distribution with precision $\theta_1$. 
Then we obtain the likelihood function $p(\cX|\cS, \cL, \theta_1)$ characterizes the probability of observing $\cX$ conditioned on $\cS$, $\cL$, and $\theta_1$. 
By exploiting the properties of the normal distribution, the likelihood function is expressed as:
\begin{equation}\label{eq:likelihood}
p(\cX|\cS, \cL, \theta_1) \propto \theta_1^{\frac{n}{2}} \exp\left(-\frac{\theta_1}{2}\|\cX - \cS - \cL\|_F^2\right),
\end{equation}
where $\propto$ denotes ``proportional to''  and $n = n_1n_2n_3$ denotes the total flattened dimensionality of $\cX$. 
This formulation captures the probabilistic nature of the constraint violation, enhancing the robustness and applicability of the Bayesian inference process.

To facilitate further analysis and optimization,  we consider   the log-likelihood function 
\[
\log p(\cX|\cS, \cL,  \theta_1) = -\frac{\theta_1}{2}\|\cX-\cS-\cL\|_F^2 + \frac{n}{2} \log \theta_1 + C_1, 
\]
where $C_1$ is a constant term that does not depend on $\cS$,  $\cL$,  or $\theta_1$ and can be ignored in inference procedure.

\subsection{The prior distributions}

In Bayesian inference,  the selection of prior distributions is a fundamental step that shapes the posterior beliefs about the unknown parameters. 
These priors encode our prior knowledge or assumptions about the variables of interest. 
Here,  we choose appropriate prior distributions for $\mathcal{S}$ and $\mathcal{L}$,  which represent distinct latent variables with unique characteristics.
We remark that the choice of these priors is informed by regularization terms. 

\subsubsection{Prior distribution for $\mathcal{S}$}
For the sparse component $\mathcal{S}$, we employ a Laplace prior distribution that induces $\ell_1$-norm regularization. This choice is motivated by the well-established connection between Laplace priors and sparsity promotion in the Bayesian framework \cite{tibshirani1996regression}. Specifically, the prior density takes the form:
\begin{equation}\label{prior_S_full}
p(\mathcal{S}|\theta_2) \propto {\theta_2^{n}}  \exp\left(-\theta_2 \|\mathcal{S}\|_1\right), 
\end{equation}
where $\theta_2 > 0$ is a scale parameter. The $\ell_1$-norm arises naturally as the convex envelope of the $\ell_0$ pseudo-norm, making it the tightest convex relaxation for sparse recovery problems. From a probabilistic perspective, this corresponds to assuming independent exponentially distributed entries in $\mathcal{S}$, which favors exact zeros in the MAP estimate while maintaining computational tractability through convex optimization.

Taking the logarithm of the prior distribution,  we obtain:
\begin{equation}
\log p(\mathcal{S}|\theta_2) = -\theta_2 \|\mathcal{S}\|_1 + n \log \theta_2 + C_2, 
\end{equation}
where $C_2$ represents a constant term that does not depend on $\mathcal{S}$ or $\theta_2$.

\subsubsection{Prior distribution for $\mathcal{L}$}
For the variable $\mathcal{L}$, 
we employ a particular Gibbs prior~ \cite{lalush1992simulation} to promote a low-rank structure in $\mathcal{L}$. This prior takes the form of an exponential distribution with a tensor nuclear norm penalty, acting as a convex surrogate for the tensor average rank. 
It encourages $\mathcal{L}$ to have a low-rank representation,  which is often suitable for capturing the underlying low-dimensionality in the data. 
The prior distribution is given by:
\begin{equation}\label{prior_L_full}
p(\mathcal{L}|\theta_3) \propto {\theta_3^{n}} \exp\left(-\theta_3 \norm{\cL}_{*}\right).
\end{equation}
This characteristic encourages the low rank property in $\cL$ and is coherent with the regularization term $\norm{\cL}_*$ in \eqref{HSImodel}.
Taking the logarithm of the prior distribution,  we have:
\[
\log p(\mathcal{L}|\theta_3) = -\theta_3 \|\mathcal{L}\|_* + n \log \theta_3 + C_3, 
\]
where $C_3$  is a constant term that does not depend on $\mathcal{L}$ or $\theta_3$.

\subsection{The hyper-prior distribution}

In the field of statistical modeling,  the Gamma distribution has obtained significant attention as a versatile prior distribution for hyperparameters,  particularly in Bayesian frameworks \cite{ ABMK,  Gamma1, MCMC,  OBF,  Gamma-non}. 
The choice of a Gamma distribution as the prior for the hyperparameter $\theta_i$ is driven by two key reasons. First, it serves as a conjugate prior for precision parameters in exponential family distributions. For instance, when $\theta_i$ controls the precision of a Gaussian likelihood $p(x|\theta_i) \sim \mathcal{N}(0, \theta_i^{-1})$, the Gamma prior ensures the posterior distribution remains a Gamma distribution. This conjugacy simplifies posterior calculations in Bayesian inference, enabling efficient automatic updates of hyperparameters. Second, $\theta_i$ typically represents positive physical quantities like precision or rate. The Gamma distribution's support on $(0, +\infty)$ naturally aligns with this positivity constraint, eliminating the need for artificial non-negativity restrictions.

We assign independent Gamma priors to the hyperparameters $\theta_i$, which correspond to the mutually independent components $\mathcal{E}$, $\mathcal{S}$, and $\mathcal{L}$ in the model. This hierarchical structure preserves model consistency while enabling efficient computation. The independence assumption further facilitates automatic feature selection by factorizing the posterior distribution into marginal products over each $\theta_i$. Hence, we have

\[
p(\theta_i)=\cG(\theta_i|a_{\theta_i}, b_{\theta_i}),  i=1, 2, 3,
\]
where $a_{\theta_i}$ and $b_{\theta_i}$ are the shape and scale parameters for each hyperparameter $\theta_i$. 
However,  a key challenge in adopting the Gamma prior lies in the determination of optimal values for $a_{\theta_i}$ and $b_{\theta_i}$. 
In the absence of strong prior knowledge,  researchers often resort to weakly informative or non-informative priors,  where the influence of the prior is minimized \cite{MCMC,  ABMK,  OBF,  Gamma1,  Gamma-non}. 
This can be achieved by setting extremely small values for $a_{\theta_i}$ and $b_{\theta_i}$ (e.g.,  $a_{\theta_i} = b_{\theta_i} = 10^{-4}$),  thereby adopting an improper prior \cite{Gamma-non}.


\subsection{Joint distribution}

The estimation of the unknown tensors $\mathcal{L}$ and $\mathcal{S}$,  given the parameters $\theta_i( i=1, 2, 3)$,  can be tackled within the Maximum A Posteriori (MAP) estimation framework. 
This approach aims to maximize the posterior density $p(\mathcal{S}, \mathcal{L}|\mathcal{X}, \btheta)$ with respect to $\mathcal{L}$ and $\mathcal{S}$,  which is formulated as:
\[
( \mathcal{S}^\dagger, \mathcal{L}^\dagger) = \arg\max_{\mathcal{S}, \mathcal{L}} p(\mathcal{S}, \mathcal{L}|\mathcal{X}, \btheta). 
\]
Applying Bayes' theorem,  the maximization problem can be rewritten in terms of the likelihood function $p(\mathcal{X}|\mathcal{S}, \mathcal{L}, \btheta)$ and the prior densities $p(\mathcal{L}|\btheta)$ and $p(\mathcal{S}|\btheta)$:
\[
\arg\max_{\mathcal{S}, \mathcal{L}} p(\mathcal{X}|\mathcal{S}, \mathcal{L}, \btheta) p(\mathcal{L}|\btheta) p(\mathcal{S}|\btheta).
\]
We remark that,  in the MAP framework,  the hyperparameters $\btheta$ must be either pre-specified or estimated prior to the estimation of $\mathcal{L}$ and $\mathcal{S}$. For a more comprehensive estimation that includes the hyperparameters,  the joint maximum a posteriori (JMAP) estimation is employed: 
\begin{equation}\label{eq:slthetamax}
    (\mathcal{S}^\dagger,  \mathcal{L}^\dagger,  \btheta^\dagger) = \arg\max_{\mathcal{S}, \mathcal{L}, \btheta} p(\mathcal{S}, \mathcal{L}, \btheta|\mathcal{X})=\argmax_{\cS, \cL, \btheta}\frac{p(\mathcal{X},  \cS, \cL, \btheta)}{p(\cX)}.
\end{equation}

For simplicity,  we assume independence among the hyperparameters,  allowing us to express the joint density function of the variables $\mathcal{X}$,  $\mathcal{S}$,  $\mathcal{L}$,  and $\btheta$ as:
\[
p(\mathcal{X},  \cS, \cL, \btheta) = p(\mathcal{X}|\cS, \cL, \btheta) p(\mathcal{L}|\theta_3) p(\mathcal{S}|\theta_2) p(\theta_1) p(\theta_2) p(\theta_3).
\]
In the literature \cite{MCMC,ABMK,OBF,Gamma1,Gamma-non}, the Gamma distribution is commonly adopted as a prior for the hyperparameters $\theta_i$ ($i=1,2,3$) due to its conjugacy with certain likelihood functions, which facilitates analytical tractability in Bayesian inference. 
However, prior knowledge about the shape and scale parameters ($a_{\theta_i}$ and $b_{\theta_i}$) of the Gamma distribution is often lacking. 
To address this, a non-informative prior can be implemented by adopting an improper uniform prior distribution, defined as $p(x) \propto 1$ for $x \in \{\theta_i \mid i=1,2,3\}$ over the positive real line \cite{Gamma-non,weakly}. 
Hence we have
\begin{equation}\label{jointdensity}
p(\cX, \cS, \cL,  \btheta)\propto \theta_1^{n/2}  \theta_2^{n}  \theta_3^{n}\exp\left( -\frac{\theta_1}{2}\norm{\cX-\cS-\cL}_F^2 -\theta_2\norm{\cS}_1-\theta_3\norm{\cL}_{*} \right). 
\end{equation}


\section{Variational Bayesian inference}\label{sectionvarpara}

In Bayesian modeling,  inference involves conditioning on observed data $\mathcal{X}$ and estimating the posterior density $p(\cS, \cL,   \btheta|\cX)$. 
This task can be tackled via Markov Chain Monte Carlo (MCMC) sampling or optimization approaches. 
However, in this paper,  we adopt variational inference as the methodological framework to approximate the latent variables $\mathcal{L}$ and $\mathcal{S}$,  along with the parameter vector $\boldsymbol{\theta}$.

\subsection{Kullback-Leibler divergence and evidence lower bound}
The central goal of variational inference is to identify an optimal variational density $q( \cS, \cL,  \btheta)$  that closely approximates the posterior density $p(\mathcal{S}, \mathcal{L}, \boldsymbol{\theta}|\mathcal{X})$,  thereby facilitating efficient inference on the latent variables and parameters \cite{blei2017variational}.

Within this framework,  we define a family of densities $\mathcal{Q}$ over the latent variables and parameters. 
Each candidate $q(\mathcal{S},  \mathcal{L},  \boldsymbol{\theta}) \in \mathcal{Q}$ represents an approximation to the true posterior. 
The optimal candidate is chosen by minimizing the Kullback-Leibler (KL) divergence from the true posterior:
\begin{equation}
\begin{aligned}
\text{KL}(q(\mathcal{S},  \mathcal{L},  \boldsymbol{\theta}) \| p(\mathcal{S},  \mathcal{L},  \boldsymbol{\theta}|\mathcal{X})) &= \int_{\mathcal{S},  \mathcal{L},  \boldsymbol{\theta}} q(\mathcal{S},  \mathcal{L},  \boldsymbol{\theta}) \log \left( \frac{q(\mathcal{S},  \mathcal{L},  \boldsymbol{\theta})}{p(\mathcal{S},  \mathcal{L},  \boldsymbol{\theta}|\mathcal{X})} \right) \,  d\mathcal{L} \,  d\mathcal{S} \,  d\boldsymbol{\theta} \\
&= \mathbb{E}_{q(\mathcal{S},  \mathcal{L},  \boldsymbol{\theta})} \left[ \log \left( \frac{q(\mathcal{S},  \mathcal{L},  \boldsymbol{\theta})}{p(\mathcal{S},  \mathcal{L},  \boldsymbol{\theta}|\mathcal{X})} \right) \right].
\end{aligned}
\end{equation}
The  variational inference task simplifies to finding the variational density $q^\dagger(\mathcal{S},  \mathcal{L},  \btheta)$ that minimizes the Kullback-Leibler (KL) divergence from the variational density to the true posterior:
\[
q^\dagger(\mathcal{S},  \mathcal{L},  \btheta) = \argmin_{q(\mathcal{S},  \mathcal{L},  \btheta) \in \mathcal{Q}} \,  \text{KL}(q(\mathcal{S},  \mathcal{L},  \btheta) \,  \| \,  p(\mathcal{S},  \mathcal{L},  \btheta|\mathcal{X})).
\]
According to \eqref{eq:slthetamax}, the posterior density $p(\mathcal{S},  \mathcal{L},  \btheta|\mathcal{X})$ is the ratio between $p(\mathcal{X},  \mathcal{S},  \mathcal{L},  \btheta)$ and $p(\mathcal{X})$. 
The density $p(\mathcal{X})$ involves integrating out the latent variables from the joint density.
Unfortunately,  this integration is often intractable,  rendering direct computation of the posterior challenging. 
Expand the condition density,  we have 
\[
\text{KL}(q(\mathcal{S}, \mathcal{L}, \btheta) \,  \| \,  p(\mathcal{S}, \mathcal{L}, \btheta|\mathcal{X}))  = -\mathbb{E}_{q(\mathcal{S},  \mathcal{L},  \btheta)} \left[ \log \left( \frac{p(\mathcal{X}, \mathcal{S},  \mathcal{L},  \btheta)}{q( \mathcal{S},  \mathcal{L},  \btheta)} \right) \right] + \log p(\mathcal{X}).
\]
The second term is independent of latent variables and hyperparameters; therefore, it’s just a constant in the minimization problem, and we can ignore this term.
To circumvent the intractability,  we optimize an alternative objective that is equivalent to the KL divergence up to an additive constant. 
Specifically,  we minimize the first term on the right-hand side of the equation,  which constitutes the evidence lower bound (ELBO),  denoted  $\mathcal{J}(q(\mathcal{S}, \mathcal{L}, \btheta))$
\begin{equation}\label{ELBO_eq}
\mathcal{J}(q(\mathcal{S}, \mathcal{L}, \btheta)) \equiv \mathbb{E}_{q(\mathcal{S}, \mathcal{L}, \btheta)}\left[\log\left(\frac{p(\mathcal{X}, \mathcal{S}, \mathcal{L}, \btheta)}{q(\mathcal{S}, \mathcal{L}, \btheta)}\right)\right].
\end{equation}
This is 
\begin{equation}\label{eq:elbo_mini}
q^\dagger(\mathcal{S}, \mathcal{L}, \btheta)=\argmax_{{q(\mathcal{S}, \mathcal{L}, \btheta) \in \cQ}} \mathcal{J}(q(\mathcal{S}, \mathcal{L}, \btheta)). 
\end{equation}



\subsection{Mean-field variational family}
To fully specify the optimization problem,  we now consider the variational family. 
The complexity of this family directly impacts the difficulty of the optimization,  with more intricate families posing greater challenges.

In this paper,  we concentrate on the mean-field variational family,  which assumes mutual independence among the latent variables,  with each variable being governed by its individual variational factor \cite{blei2017variational}. 
This assumption simplifies the variational density into a factorized form:
\[
q(\mathcal{S}, \mathcal{L}, \btheta) = q(\mathcal{L})q(\mathcal{S})\prod_{i=1}^{3}q(\theta_i).
\]

The selection of variational densities $q(\mathcal{L})$,  $q(\mathcal{S})$,  and $q(\theta_i)$  is importance. 
For $q(\mathcal{L})$ and $q(\mathcal{S})$,  we adopt normal distributions due to their versatility and analytical convenience. 
The choice of the variational density $q(\theta_i) $ as a Gamma distribution is motivated by the conjugacy properties derived from the likelihood function (Eq. \ref{eq:likelihood}) and the prior distributions specified in Eqs. \eqref{prior_S_full} and \eqref{prior_L_full}. These indicate that the posterior distribution of $\theta_1$ and the conditional posteriors of $\theta_i$ for $i=2, 3$ follow Gamma distributions. Since the Gamma distribution is conjugate to itself, selecting $q(\theta_i) $ as a Gamma density ensures compatibility with the posterior, facilitating efficient variational inference.

Let \(\mathcal{Q}_{\mathcal{G}}\) denote the set of Gamma densities for the hyperparameters \(\theta_i\) (\(i=1,2,3\)), and \(\mathcal{Q}_{\mathcal{N}}\) denote the set of multivariate normal densities over the tensor space \(\mathbb{R}^{n_1 \times n_2 \times n_3}\).
The overall variational family $\mathcal{Q}$ can be expressed as the Cartesian product of these sets: $\mathcal{Q} = \mathcal{Q}_{\mathcal{N}} \times \mathcal{Q}_{\mathcal{N}} \times \mathcal{Q}_{\mathcal{G}}$.

\subsection{Laplacian approximation}

In \eqref{jointdensity},  the non-quadratic properties inherent in both the $\ell_1$ norm of $\cS$,  which represents the sum of the absolute values of all elements, and the tensor nuclear norm of $\cL$,  which is the weighted sum of its singular values,  pose significant obstacles for direct optimization within standard density families. 
These non-quadraticities complicate direct inference procedures,  rendering them computationally intractable.
To address this,  we utilize the Laplace approximation method,  involving mean calculation,  variance estimation,  and density function construction,  to approximate the density with a Gaussian distribution form.

Here, we consider a general density function $q(x)$ with a single random variable $x$ and simplify  \eqref{eq:elbo_mini}  as 
\[
q^\dagger(x) = \argmax_{q(x)\in \cQ_{\cN}} \int_{\Omega} q(x) \log \frac{p(x)}{q(x)} dx.
\]
where $\cQ_{\cN}$ is the set of all the density functions for the Gaussian distribution.  
According to Gibbs' inequality, for any two probability distributions $q(x)$ and $p(x)$ over a domain $\Omega$, the following holds:
\[
\int_{\Omega} q(x) \log \frac{p(x)}{q(x)} dx \leq 0,
\]
with equality achieved if and only if $q(x) = p(x)$, implying identical means and variances. 
However, the non-quadratic nature of the $\log p(x)$ term complicates the direct estimation of $q(x)$ in practice.
To address this, we employ the Laplacian approximation method to estimate $q(x)$. Since $q(x)$ is Gaussian, we have the following properties:
\begin{enumerate}[label=({\arabic*)},leftmargin=4ex, itemindent=0ex]
\item First-Order Condition for the Mean ($\sfE_x$): The gradient of $\log q(x)$ evaluated at $\sfE_x$ is zero, implying $\sfE_x$ is a maximum of $\log p(x)$.
\item Second-Order Condition for the Variance ($\sigma_x^2$): The negative of the Hessian (second-order derivative) of $\log q(x)$ evaluated at $\sfE_x$ equals the reciprocal of the variance. However, since we directly approximate $p(x)$, we use the Hessian of $\log p(x)$ evaluated at $\sfE_x$ to estimate $\sigma_x^2$:
\[
   \nabla^2 \log p(x)|_{x=\sfE_x} = -\frac{1}{\sigma_x^2}.
\]
\end{enumerate}
We now detail the estimation of $\sfE_x$ and $\sigma_x^2$ based on these conditions for some specific density function $p(x)$.

\subsubsection{Absolute value function}\label{LaplaceL1}
The $\ell_1$ norm of $\cS$,  as the sum of absolute element values,  necessitates approximating the distribution of absolute values to enable effective optimization within Gaussian density families. 
Given the log-probability density function $\log p(x) \propto -\frac{1}{2}(x-b)^2 - \beta |x|$, 
the first step of Laplace approximation involves computing the mean $\sfE_{x}$ of $q(x)$,
which corresponds to the maximum of $\log p(x)$:
\[
\sfE_{x} = \argmax_{x} \log p(x)=\argmin_{x} \frac{1}{2}(x-b)^2 + \beta |x|.
\]
In the paper, we utilize the sans serif font  $\sfE$ accompanied by a subscripted variable $x$ to denote the expectation of a random variable $x$. 
The  solution   is given by:
\[
\sfE_{x} = \begin{cases}
b - \beta \mathrm{sign}(b),  & \text{if } |b| > \beta \\
0,  & \text{otherwise}.
\end{cases}
\]
Proceeding to the second stage,  we estimate the variance,  $\sigma^2$. When $\sfE_{x} = 0$,  we directly set $\sigma_x^2 = 0$. 
For non-zero $\sfE_{x}$,  we leverage the inequality $\abs{x} \leq \frac{x^2}{2|y|} + \frac{|y|}{2}$ with equality at $x = y \neq 0$.
By setting $y = |\sfE_{x}|$, this facilitates a lower bound on $\log p(x)$:
\[
\log p(x) \geq -\frac{1}{2}(x-b)^2 - \frac{\beta}{2|\sfE_{x}|}x^2 + \text{const}.
\]
To simplify the variance estimation process,  we exclude the constant term and mean shift from consideration,  as they do not impact the variance calculation.
Approximating the second-order derivatives of this lower bound around $\sfE_{x}$,  we derive the variance estimate:
\[
\sigma_x^2 \approx \left(1 + \frac{\beta}{|\sfE_{x}|}\right)^{-1} = \frac{|\sfE_{x}|}{|\sfE_{x}| + \beta}.
\]
Ultimately,  utilizing the estimated mean $\sfE_{x}$ and variance $\sigma_x^2$,  we construct the optimal Gaussian density approximation:
\begin{equation}\label{densityabs}
q(x) = \mathcal{N}\left(\sfE_{x},  \frac{|\sfE_{x}|}{|\sfE_{x}| + \beta}\right).  
\end{equation}
Given the approximation $|x| \approx \frac{1}{2|\sfE_{x}|}x^2 + \frac{|\sfE_{x}|}{2}$ at $x = \sfE_{x}$,  we derive the expectation of the absolute value of $x$:
\begin{equation}\label{Eabsxsqrt}
\sfE|x| = |\sfE_{x}| + \frac{1}{2(|\sfE_{x}| + \beta)}. 
\end{equation}

\subsubsection{Nuclear norm}\label{LaplaceNN}

The Laplace approximation approach does not directly extend the Gaussian density approximation of absolute functions to the nuclear norm of matrices,  given its intrinsic nature as a sum of singular values rather than element-wise absolute values. 
When considering a density function incorporating the weighted nuclear norm of a matrix $\mathbf{X} \in \mathbb{R}^{n_1 \times n_2}$,  we assume $n_1 \leq n_2$ for generality. The targeted density is formulated as:
\[
p(\bX) \propto \exp\left(-\frac{\alpha}{2}\|\bX-\bA\|_F^2 - \beta \|\bX\|_{*}\right), 
\]
where $\bA$ is a given matrix,  $\beta$ is a regularization parameter,  and $\vw \in \mathbb{R}^{n_2}$ represents the vector of weights. 
We find the density   $q(\bX)$ that maximizes:
\[
q(\bX) = \arg\max_{q(\bX) \in \cQ_\mathcal{N}} \int q(\bX) \log \frac{p(\bX)}{q(\bX)} \,  d\bX .
\]
Applying the Laplacian approximation method,  the mean $\sfE_{\bX}$ of   $q(\bX)$   is obtained by solving:
\[
\sfE_{\bX} = \arg\min_{\bX} \left(\frac{\alpha}{2}\norm{\bX-\bA}_F^2 + \beta \norm{\bX}_{*}\right).
\]
Given the SVD  of $\bA$ as $\bA = \bU_{\bA}\bD_{\bA}\bV_{\bA}^T$, the minimizer  $\sfE_{\bX}$ for the aforementioned problem can be formulated as \cite{caijf2010}:
\[
\sfE_{\bX} = \bU_{\bA} \max\{\bD_{\bA} - \tfrac{\beta}{\alpha}\bI, 0\} \bV_{\bA}^T,
\]
where $\max\{\cdot, 0\}$ denotes an element-wise maximum operation applied to the diagonal matrix.
To compute the covariance of $\bX$,  we introduce an inequality derived from \cite{lefkimmiatis2023}:
\[
\norm{\bX}_{*} \leq \frac{1}{2}\text{Tr}(\omega(\bY)\bX\bX^T) +  \frac{1}{2}\norm{\bY}_{*}, 
\]
where $\omega(\bY) = (\bY\bY^T)^{-1/2}$, 
and equality holds when $\bX = \bY$. 
Hence we obtain
\begin{equation}\label{Xnuclearnormapprox}
\norm{\bX}_{*}\approx \frac{1}{2}\text{Tr}(\omega(\sfE_{\bX})\bX\bX^T)
+ \frac{1}{2}\|\sfE_{\bX}\|_{*}.  
\end{equation}
Considering the $j$-th columns of $\bX$ and $\bA$ denoted by $\bX_{:j}$ and $\bA_{:j}$ respectively,  we bound the log-likelihood $\log p(\bX)$ as follows: 
\[
\log p(\bX) \geq -\sum_j \left(\frac{\alpha}{2}\|\bX_{:j}-\bA_{:j}\|_2^2 + \frac{\beta}{2} \bX_{:j}^T\omega({\sfE_{\bX}})\bX_{:j}\right) + \text{const}.
\]
Evaluating the second-order derivatives of the lower bound with respect to $\bX_{:j}$ yields the inverse covariance matrix $\bSig_{\bX_{:j}}^{-1} = \alpha\bI + \beta \omega({\sfE_{\bX}}).$
Let $\sfE_{\bX} =\bU\bD\bV^T$ be the skinny SVD of $\sfE_{\bX} $,  we have 
$\omega({\sfE_{\bX}}) = \bU  
\bD^{-1} {\bU}^T. $
Hence  we have 
\begin{equation}\label{eq:sigma_x}
    \bSig_{\bX_{:j}}={\bU} {\bD}  \left(\alpha{\bD}  + \beta \bI\right)^{-1}{\bU}^T.
\end{equation}
Finally,   the optimal density approximation $q(\bX)$  is expressed as:
\begin{equation}\label{XGauusianapproxdensity}
q(\bX) = \prod_j \mathcal{N}(\bX_{:j}|\sfE_{\bX_{:j}}, {\bU} {\bD}  \left(\alpha{\bD}  + \beta \bI\right)^{-1}{\bU}^T).
\end{equation}
We have 
\begin{equation}\label{EmatrixX2}
\E\norm{\bX}_F^2 = \norm{\sfE_{\bX}}_F^2 +\sum_j\trace{\bSig_{\bX_{:j}}}= \norm{\sfE_{\bX}}_F^2+ n_2\sum_i \frac{\bD_i}{\alpha\bD_i+\beta}.  
\end{equation}
According to the approximation of the nuclear norm in \eqref{Xnuclearnormapprox},  
we have
\begin{equation}\label{Ematrixsq}
\E\norm{\bX}_{*}= \frac{1}{2} \sum_j\E \left[\bX_{:j}^T\omega({\sfE_{\bX}})\bX_{:j}\right]
+ \frac{1}{2}\|\sfE_{\bX}\|_{*}. 
\end{equation}
Let ${\bD_{i}}$ be the $i$-th singular values of $\sfE_{\bX}$. 
Then we have
\begin{equation}\label{Ematrixwnn}
\E\norm{\bX}_{*}
= \norm{\sfE_{\bX}}_{*} +\frac{n_2}{2}\trace{ \left(\alpha{\bD}  + \beta \bI\right)^{-1}}
=\norm{\sfE_{\bX}}_{*} +\frac{n_2}{2}\sum_i \frac{1}{\alpha\bD_i+\beta}.
\end{equation}

\subsection{ Coordinate ascent variational inference} 
In order to maximize the ELBO $\cJ (q(\mathcal{S}, \mathcal{L}, \btheta))$,  we apply coordinate ascent variational inference (CAVI) \cite{blei2017variational, wang2013variational}. Starting from an initial density
$(q_0(\btheta), q_0(\mathcal{L}), q_0(\mathcal{S}))$,  the densities of $\mathcal{S}$,  $\mathcal{L}$ and $\btheta$ are updated as follows:
\begin{eqnarray}
q_\ell(\mathcal{S}) & = &\underset{q(\mathcal{S}) \in \mathcal{Q}_{\mathcal{N}}}{\operatorname{argmax}} \cJ \left(q(\mathcal{S}) q_{\ell-1}(\mathcal{L}) q_{\ell-1}(\btheta) \right),  \label{def:qks} \\
 q_\ell(\mathcal{L})  & = &\underset{q(\mathcal{L}) \in \mathcal{Q}_{\mathcal{N}}}{\operatorname{argmax}} \cJ \left(q_{\ell}(\mathcal{S}) q(\mathcal{L}) q_{\ell-1}(\btheta)\right),  \label{def:qkx}\\
q_\ell(\btheta) & = &\underset{q(\btheta) \in \mathcal{Q}_{\mathcal{G}}}{\operatorname{argmax}} \cJ \left(q_{\ell}(\mathcal{S})q_{\ell}(\mathcal{L}) q(\btheta)\right), \label{def:qktheta}   
\end{eqnarray}
where $q_\ell(\btheta),  q_\ell(\mathcal{L}),  q_\ell(\mathcal{S})$ refer to the variational densities obtained in the $\ell$-th iteration.

\subsubsection{The density  $q_\ell(\mathcal{S})$}

In accordance with  \eqref{ELBO_eq},  we formulate the optimization problem as maximizing the evidence lower bound (ELBO) with respect to the variational distribution $q(\mathcal{S})$:
\begin{equation}
\label{eq:max_q_S}
    \argmax_{q(\cS)\in \cQ_{\cN} }\mathcal{J}(q(\cS), q_{\ell-1}(\cL, \btheta)) = \argmax_{q(\cS)\in \cQ_{\cN}}\int q(\cS){\mathbb{E}_{q_{\ell-1}(\cL, \btheta)}\log \frac{p(\cX, \cS, \cL, \btheta)}{ q(\cS)} }
d\cS.
\end{equation}
Given the joint density distribution as defined in \eqref{jointdensity},  we can express the expectation term within the ELBO as:
\begin{eqnarray*}
& &\E_{q_{\ell-1}(\cL, \btheta)}\left[\log p(\cX, \cS, \cL, \btheta)\right] \\
&=& -\sum_{ijk}\left(\frac{\sfE_{\theta_1}^{\ell-1}}{2}\left(\cX_{ijk} - \sfE_{\cL_{ijk}}^{\ell-1} - \cS_{ijk}\right)^2+ \sfE_{\theta_2}^{\ell-1}\abs{\cS_{ijk}}\right) + \text{const},     
\end{eqnarray*}
where $\text{const}$ is a constant independent of $\cS$.
According to the discussion in Section \ref{LaplaceL1},  the mean of $\mathcal{S}_{ijk}$ is given by 
\[
\sfE_{\cS}^{\ell} 
=\argmin_{\cS}\frac{\sfE_{\theta_1}^{\ell-1}}{2}\norm{\mathcal{X} - \sfE_{\cL}^{\ell-1} - \mathcal{S}}_F^2 + {\sfE_{\theta_2}^{\ell-1}}\norm{\mathcal{S}}_1. 
\]
It is known that the minimizer is the well-known soft threshold: 
\begin{equation}\label{expectationSijk}
\sfE_{\cS_{ijk}}^{\ell} =\left\{\begin{array}{ll}
\mathcal{X}_{ijk} - \sfE_{\cL_{ijk}}^{\ell-1}- \frac{\sfE_{\theta_2}^{\ell-1}}{\sfE_{\theta_1}^{\ell-1}}, & \mbox{if   }\mathcal{X}_{ijk} - \sfE_{\cL_{ijk}}^{\ell-1}\geq \frac{\sfE_{\theta_2}^{\ell-1}}{\sfE_{\theta_1}^{\ell-1}},\\
\mathcal{X}_{ijk} - \sfE_{\cL_{ijk}}^{\ell-1}+ \frac{\sfE_{\theta_2}^{\ell-1}}{\sfE_{\theta_1}^{\ell-1}}, & \mbox{if   }\mathcal{X}_{ijk} - \sfE_{\cL_{ijk}}^{\ell-1}\leq -\frac{\sfE_{\theta_2}^{\ell-1}}{\sfE_{\theta_1}^{\ell-1}},\\
0, & \mbox{others}.\\
\end{array} \right. 
\end{equation}
Applying \eqref{densityabs},  the variance of  $\mathcal{S}_{ijk}$ is given by 
\[
\bSig_{\mathcal{S}_{ijk}}^{\ell}  = \frac{\sfE_{\theta_1}^{\ell-1}|\sfE_{\cS_{ijk}}^{\ell}|}{ \sfE_{\theta_1}^{\ell-1}|\sfE_{\cS_{ijk}}^{\ell}|+\sfE_{\theta_2}^{\ell-1}}.
\]
Then  the density function of  $q(\mathcal{S})$ is given:
\begin{equation}\label{densityqS}
q_{\ell}(\mathcal{S}_{ijk}) = \mathcal{N}(\mathcal{S}|\sfE_{\cS_{ijk}}^\ell,  \bSig_{\mathcal{S}_{ijk}}^\ell).  
\end{equation}

\subsubsection{The density $q_\ell(\mathcal{L})$}

In accordance with  \eqref{ELBO_eq},  we formulate the optimization problem as maximizing the evidence lower bound (ELBO) with respect to the variational distribution $q(\mathcal{L})$:
\begin{equation}
\label{eq:max_q_L}
\argmax_{q(\cL)\in \cQ_{\cN} }\cJ(q_{\ell}(\cS), q(\cL), q_{\ell-1}(\btheta)) = \argmax_{q(\cL)\in \cQ_{\cN}}\int q(\cL){\mathbb{E}_{q_{\ell}(\cS)q_{\ell-1}(\btheta)}\log\tfrac{ p(\cX, \cS, \cL, \btheta)}{q(\cL)}}
d\cL.
\end{equation}
Given the joint density distribution as defined in \eqref{jointdensity},  we can express the expectation term within the ELBO as:
\begin{equation*}
\begin{split}
    \E_{q_{\ell}(\cS)q_{\ell-1}(\btheta)}\left[\log p(\cX, \cS, \cL, \btheta)\right]&=- \left(\frac{\sfE_{\theta_1}^{\ell-1}}{2}\norm{\cL-(\cX-\sfE_{\cS}^{\ell})}_F^2+
\sfE_{\theta_3}^{\ell-1}\norm{\cL}_{*} \right)+ \text{const} \\ 
 &=- \left(\frac{\sfE_{\theta_1}^{\ell-1}}{2n_3}\norm{\ocL-(\ocX-\overline{\sfE_{\cS}^{\ell}})}_F^2+
\frac{\sfE_{\theta_3}^{\ell-1}}{n_3}\norm{\ocL}_{*} \right)+ \text{const} 
\end{split}
\end{equation*}
According to the discussion in Section \ref{LaplaceNN},  the mean of $\mathcal{L}$ is given by 
\[
\sfE_{\cL}^{\ell}=\argmin_{\cL}\tfrac{\sfE_{\theta_1}^{\ell-1}}{2}\norm{\cX - \cL- \sfE_{\cS}^{\ell} }_F^2 +{\sfE_{\theta_3}^{\ell-1}} \norm{\cL}_{*}.
\]
This subproblem is to solve a proximal operator of the tensor nuclear norm, which has a closed-form solution as tensor singular value thresholding (t-SVT) \cite{lu2019tensor}. 
Let the SVD of $\cX - \sfE_{\cS}^{\ell}$  is given by  
$\cX - \sfE_{\cS}^{\ell} = \cU^{\ell}*\cD^{\ell}*{\cV^{\ell}}^T. $ The update of $\sfE_{\cL}^\ell$ is 
\begin{equation}\label{eq:sol_tst}
    \sfE_{\cL}^\ell = \cU^{\ell}*\cD^{\ell}_\tau*{\cV^{\ell}}^T,
\end{equation}
where ${\cD^\ell_{\tau}}$ is an ${n_{1} \times n_{2} \times n_{3}}$ tensor that satisfies $\overline{\mathcal{D}}^\ell_{\tau}=\max \{\overline{\mathcal{D}^\ell}-\tau,  0\}$ with $\tau = \sfE_{\theta_1}^{\ell-1}/\sfE_{\theta_3}^{\ell-1}$. 
Recall that
we adopt the notation of an overline $\overline{\mathcal{A}}$  to signify the application of the  DFT to the tensor $\mathcal{A}$ specifically along its third dimension.


We apply \eqref{eq:sigma_x} and then obtain the covariance matrix of the vector $\ocL_{:jk}$
\[
\bSig_{\ocL_{:jk}}^\ell= n_3\ocU_{::k}^\ell{\overline{\cD}_\tau^{\ell}}_{::k}   \left(\sfE_{\theta_1}^{\ell-1}{\overline{\cD}_\tau^{\ell}}_{::k}  + {\sfE_{\theta_3}^{\ell-1}} \bI\right)^{-1}{\ocU_{::k}^\ell}^T.
\]
Thus,  
we construct the density function and parameterize a normal density $q(\ocL)$ as:
\begin{equation}\label{densityqL}
q_{\ell}(\ocL)=\prod_{jk}\cN\left(\ocL_{:jk}|\sfE_{\ocL_{:jk} }^\ell, \bSig_{\ocL_{:jk}}^\ell\right).
\end{equation}

\subsubsection{The density $q_\ell(\btheta)$}
In accordance with  \eqref{ELBO_eq},  we frame the optimization problem as maximizing the evidence lower bound (ELBO) with respect to the variational distribution $q(\btheta)$, 
\begin{equation}
\label{eq:max_q_theta}
\argmax_{q(\btheta)\in \cQ_{\cG} }\cJ(q_{\ell}(\cS), q_{\ell}(\cL), q(\btheta)) = \argmax_{q(\btheta)\in \cQ_{\cG}}\int q(\btheta)\E_{q_\ell(\cS, \cL)}\log\tfrac{ p(\cX, \cS, \cL, \btheta)}{ q(\btheta)}
d\btheta,
\end{equation}
where $\cQ_{\cG}$ is the set of all the density functions for the Gamma distribution. 
By taking the partial derivative of the objective function in \eqref{eq:max_q_theta} with respect to $q(\btheta)$, and letting it be equal to 0, we obtain that the optimal $q(\btheta)$ is proportional to 
\[
q(\btheta)\propto \exp{\E_{q_\ell(\cS, \cL)}\log p(\mathcal{X},  \cS, \cL,  \btheta)}
\]
with
\begin{eqnarray*}
\E_{q_\ell(\cS, \cL)}\log p(\mathcal{X},  \cS, \cL,  \btheta)
&=& -\tfrac{\theta_1}{2} \E_{q_\ell(\cS, \cL)}\norm{\cX-\cS-\cL}_F^2 -\theta_2\E_{q_\ell(\cS)}\norm{\cS}_1\\
&&   -\theta_3\E_{q_\ell(\cL)}\norm{\cL}_{*}  +\tfrac{n}{2} \log  \theta_1 +n\log  \theta_2 
 +n\log \theta_3+ \text{const}. 
\end{eqnarray*}
Since \(q(\btheta) = \prod_{i=1}^3 q(\theta_i)\) is assumed to factorize due to the independence of the model components, we derive the form of each \(q(\theta_i)\) by comparing the coefficients of \(\theta_i\) and \(\log \theta_i\)  with the log-density of a Gamma distribution
\[
\log q(\btheta)=\sum_{i=1}^3\left((a_{\theta_i}^\ell-1)\log \theta_i - b_{\theta_i}^\ell \theta_i\right) + \text{const},
\]
where \(a_{\theta_i}^\ell\) and \(b_{\theta_i}^\ell\) are the shape and rate parameters, respectively.
By comparing the coefficients in $\E_{q_\ell(\cS, \cL)}\log p(\mathcal{X},  \cS, \cL,  \btheta)$ with those of a Gamma density ($\cG(x|a, b)\propto x^{a-1}\exp(-bx )$),  we can infer the shape $a_{\theta_i}^\ell$ and scale $b_{\theta_i}^\ell$ parameters for each $\theta_i$. Consequently,  the shape parameters are given by:

\[
a_{\theta_1}^\ell = \tfrac{n}{2}+1,  \quad a_{\theta_2}^\ell = n+1,  \quad a_{\theta_3}^\ell = n+1.
\]
While the scale parameters are expressed as expectations over the variational distributions $q_\ell(\cS)$ and $q_\ell(\cL)$,  as defined in the following system of equations:
\begin{equation}\label{scale_parameters}
\left\{
\begin{array}{l}
b_{\theta_1}^\ell = \frac{1}{2}\mathbb{E}_{q_{\ell}(\mathcal{L})q_{\ell}(\mathcal{S})}\left[\norm{\mathcal{X}-\mathcal{S}-\mathcal{L}}_F^2\right],  \\
b_{\theta_2}^\ell = \mathbb{E}_{q_{\ell}(\mathcal{S})} \norm{\mathcal{S}}_1 ,  \\
b_{\theta_3}^\ell = \mathbb{E}_{q_{\ell}(\mathcal{L})} \norm{\mathcal{L}}_{*} .
\end{array}
\right.
\end{equation}

The computation for $b_{\theta_1}^\ell $ involves the expectations of both $\norm{\cS}_F^2$ and $\norm{\cL}_F^2$. 
It is obvious that 
\[
\E_{q_{\ell}(\mathcal{S})}\norm{\cS}_F^2=\sum_{ijk}\E_{q_{\ell}}\abs{\cS_{ijk}}^2=\sum_{ijk}\left(\abs{\sfE_{\cS_{ijk}}^\ell}^2+\bSig_{\cS_{ijk}}^{\ell}\right). 
\]
According to \eqref{EmatrixX2}, we have
\[
\E_{q_{\ell}(\cL)} \norm{\cL}_F^2=\frac{1}{n_3}\sum_{j,k}\E_{q_{\ell}(\ocL)}\norm{\ocL_{:jk}}_F^2=\frac{1}{n_3}\sum_{j,k}\left(\norm{\sfE_{\ocL_{:jk}}^\ell}_F^2+\trace{\bSig_{\ocL_{:jk}}^\ell} \right).
\]
Hence we have 
\[
\mathbb{E}_{q_{\ell}(\mathcal{L})q_{\ell}(\mathcal{S})}(\|\mathcal{X}-\mathcal{S}-\mathcal{L}\|_F^2) = \norm{\mathcal{X} -\sfE_{\mathcal{L}}^{\ell}-\sfE_{\mathcal{S}}^{\ell}}_F^2 + \frac{1}{n_3}\sum_{j,k}\trace{\bSig_{\ocL_{:jk}}} + \sum_{ijk} \bSig_{\cS_{ijk}}^{\ell}.
\]

For the expectation of $\|\mathcal{S}\|_1$, according to  \eqref{Eabsxsqrt}, we have
\[
\mathbb{E}_{q_\ell(\mathcal{S})}\|\mathcal{S}\|_1 = \sum_{ijk} \mathbb{E}_{q_\ell(\mathcal{S})} |\mathcal{S}_{ijk}| = \|\sfE_{\mathcal{S}}^\ell\|_1 + \frac{1}{2}\sum_{ijk}\left( \sfE_{\theta_1}^{\ell} |\sfE_{\cS_{ijk}}^{\ell}| + \sfE_{\theta_2}^{\ell} \right)^{-1}.
\]
For the expectation of the nuclear norm $\norm{\cL}_{*} $, it is the arithmetic mean of each slice $\ocL_{::k}$ of the tensor $\mathcal{L}$.
Hence we need to evaluate $\mathbb{E}_{q_{\ell}(\ocL_{::k})}\|\ocL_{::k}\|_{*}$. 
According to \eqref{Ematrixwnn}, we have
\[
\E_{q_{\ell}(\ocL_{::k})} \|\ocL_{::k}\|_{*}=\norm{\sfE_{\ocL_{::k}}^\ell}_{*}+\frac{n_2n_3}{2}\trace{ \left(\sfE_{\theta_1}^{\ell-1}{\overline{\cD}_\tau^{\ell}}_{::k}  + {\sfE_{\theta_3}^{\ell-1}} \bI\right)^{-1}}.
\]
Hence
\[
\E_{q_\ell(\cL)}\norm{\cL}_{*}=\norm{\sfE_{\cL}^\ell}_{*}+\frac{n_2}{2}\sum_{k}\trace{  \left(\sfE_{\theta_1}^{\ell-1}{\overline{\cD}_\tau^{\ell}}_{::k} + {\sfE_{\theta_3}^{\ell-1}} \bI\right)^{-1}}.
\]
Now,  we focus on the expectation of the nuclear norm $\norm{\cL}_{*} $,  which requires evaluating $\mathbb{E}_{q_{\ell}(\ocL_{::k})} \|\ocL_{::k}\|_{*} $ for each slice $\ocL_{::k}$ of the tensor $\mathcal{L}$.

We summarize the proposed adaptive method in Algorithm \ref{alg:paraalgmain}. For simplicity, we refer to our proposed algorithm for solving the tensor nuclear norm model as VBI$_{\rm TNN}$.

\begin{algorithm}[th ] 
\caption{VBI$_{\rm TNN}: $Variational Bayesian inference for the TNN-based TRPCA. } 
\label{alg:paraalgmain} 
\begin{algorithmic}[1] 
\STATE{{\bf Initialization: } $\sfE_{\theta_1}, \sfE_{\theta_2}, \sfE_{\theta_3}$,  $\sfE_{\mathcal{L}}^0,  \sfE_{\mathcal{S}}^0$,  $\boldsymbol{\Sigma}_{\mathcal{L}}^0, \boldsymbol{\Sigma}_{\mathcal{X}}^0$}
\STATE Let $a_{\theta_1} = \tfrac{n}{2}+1,  \quad a_{\theta_2} = n+1,  \quad a_{\theta_3} = n+1.$
\WHILE {$\ell \leq \ell_{\rm{Max}}$ or not converged}
\STATE $\sfE_{\cS_{ijk}}^{\ell} =\left\{\begin{array}{ll}
\mathcal{X}_{ijk} - \sfE_{\cL_{ijk}}^{\ell-1}- \frac{\sfE_{\theta_2}^{\ell-1}}{\sfE_{\theta_1}^{\ell-1}}, & \mbox{if   }\mathcal{X}_{ijk} - \sfE_{\cL_{ijk}}^{\ell-1}\geq \frac{\sfE_{\theta_2}^{\ell-1}}{\sfE_{\theta_1}^{\ell-1}}\\
\mathcal{X}_{ijk} - \sfE_{\cL_{ijk}}^{\ell-1}+ \frac{\sfE_{\theta_2}^{\ell-1}}{\sfE_{\theta_1}^{\ell-1}}, & \mbox{if   }\mathcal{X}_{ijk} - \sfE_{\cL_{ijk}}^{\ell-1}\leq -\frac{\sfE_{\theta_2}^{\ell-1}}{\sfE_{\theta_1}^{\ell-1}}\\
0, & \mbox{others }\\
\end{array} \right.$
\STATE Take the SVD of $\cX - \sfE_{\cS}^{\ell}$  as  
$\cX - \sfE_{\cS}^{\ell} = \cU^{\ell}*\cD^{\ell}*{\cV^{\ell}}^T $
\STATE $\sfE_{\cL}^\ell = \cU^{\ell}*\cD^{\ell}_\tau*{\cV^{\ell}}^T$
\STATE $\bSig_{\mathcal{S}_{ijk}}^{\ell}  = \frac{\sfE_{\theta_1}^{\ell-1}|\sfE_{\cS_{ijk}}^{\ell}|}{ \sfE_{\theta_1}^{\ell-1}|\sfE_{\cS_{ijk}}^{\ell}|+\sfE_{\theta_2}^{\ell-1}}, $
and  $\bSig_{\ocL_{:jk}}^\ell= n_3\ocU_{::k}^\ell{\overline{\cD}_\tau^{\ell}}_{::k}   \left(\sfE_{\theta_1}^{\ell-1}{\overline{\cD}_\tau^{\ell}}_{::k}  + {\sfE_{\theta_3}^{\ell-1}} \bI\right)^{-1}{\ocU_{::k}^\ell}^T$
\STATE $q(\mathcal{S}_{ijk}) = \mathcal{N}(\mathcal{S}|\sfE_{\cS_{ijk}}^\ell, \bSig_{\mathcal{S}_{ijk}}^\ell)$ and $q(\ocL)=\prod_{jk}\cN\left(\ocL_{:jk}|\sfE_{\ocL_{:jk} }^\ell, \bSig_{\ocL_{:jk}}^\ell\right).$
\STATE $b_{\theta_1}^\ell =  \norm{\mathcal{X} -\sfE_{\mathcal{L}}^{\ell}-\sfE_{\mathcal{S}}^{\ell}}_F^2/2 + \frac{1}{2n_3}\sum_{j,k}\trace{\bSig_{\ocL_{:jk}}} + \sum_{ijk} \bSig_{\cS_{ijk}}^{\ell} /2$
\STATE $b_{\theta_2}^\ell = \|\sfE_{\mathcal{S}}^\ell\|_1 + \frac{1}{2}\sum_{ijk}\left( \sfE_{\theta_1}^{\ell} |\sfE_{\cS_{ijk}}^{\ell}| + \sfE_{\theta_2}^{\ell} \right)^{-1}$
\STATE $b_{\theta_3}^\ell=\norm{\sfE_{\cL}^\ell}_{*}+\frac{n_2}{2}\sum_{k}\trace{  \left(\sfE_{\theta_1}^{\ell-1}{\overline{\cD}_\tau^{\ell}}_{::k}  + {\sfE_{\theta_3}^{\ell-1}} \bI\right)^{-1}}$
\STATE $q(\theta_i)=\cG(\theta_i|a_{\theta_i}, b_{\theta_i}^\ell),  \text{ and } \sfE_{\theta_i}^\ell = a_{\theta_i}/b_{\theta_i}^\ell  i=1, 2, 3  $
\ENDWHILE
\RETURN $\mathcal{L} =\sfE_\mathcal{L}^{\ell}$,  $\mathcal{S} =\sfE_\mathcal{S}^{\ell}$  
\end{algorithmic}
\end{algorithm}

\begin{remark}
According to \cite{bhattacharya2025}, a general theoretical treatment of analyzing the convergence of CAVI is missing in the literature. This is due to the lack of tractability of the updating formula involving unwieldy normalization constants and the technical challenge of dealing with optimization over infinite-dimensional distributions. Here, we will empirically show the convergence in Section~\ref{sectionexppara}.
\end{remark}

\subsection{Variational Bayesian inference for weighted tensor nuclear norm}

In this subsection, we consider a variant of the tensor nuclear norm by reweighting the singular values \cite{jiang2020multi,chang2020weighted}. 
Note that the standard tensor nuclear norm can be regarded as a special version of the weighted tensor nuclear norm, where the weighting matrix consists of elements that are all equal to one.
Formally, for a non-negative matrix $\bW \in \mathbb{R}^{\min(n_1, n_2) \times n_3}$ with column vectors $\bW_{:k}$, the weighted tensor nuclear norm $\norm{\cA}_{\bW *}$ is defined as:
\[
\norm{\cA}_{\bW *} = \frac{1}{n_3} \sum_{k=1}^{n_3} \sum_{j=1}^{\min(n_1, n_2)} \bW_{jk} \sigma_{jk},
\]
where $\sigma_{jk}$ denotes the $j$-th singular value of the $k$-th frontal slice $\cA_{::k}$ of tensor $\cA$.
To incorporate this weighted norm, we modify the robust principal component model \eqref{HSImodel} as follows:

\begin{equation}\label{WHSImodel}
\min_{\cS, \cL} \left\{ \frac{\theta_1}{2} \|\mathcal{X} - \mathcal{L} - \mathcal{S}\|_F^2 + \theta_2 \|\mathcal{S}\|_1 + \theta_3 \|\mathcal{L}\|_{\bW *} \right\}.
\end{equation}

During the inference of $\cL$, we update the expectation of $\ocL_{::k}$ in \eqref{densityqL} to:
\begin{equation}\label{eq:sol_wtst}
    \sfE_{\cL}^\ell = \cU^{\ell}*\cD^{\ell}_{\bW}*{\cV^{\ell}}^T,
\end{equation}
where ${\cD^\ell_{\bW}}$ is an ${n_{1} \times n_{2} \times n_{3}}$ tensor that satisfies $${\overline{\mathcal{D}}^\ell_{\bW}}_{::k}=\max \left\{{\overline{\mathcal{D}^\ell}}_{::k}-\tfrac{\sfE_{\theta_1}^{\ell-1}}{\sfE_{\theta_3}^{\ell-1}}\diag (\bW_{:k}),  0\right\}.$$ 
Concurrently, the covariance matrix of $\ocL_{::k}$ is adjusted to:
\[
\bSig_{\ocL_{:jk}}^\ell = n_3\ocU_{::k}^\ell{\overline{\cD}^{\ell}}_{::k}(\tau^\ell) \left( \sfE_{\theta_1}^{\ell-1}{\overline{\cD}_\tau^{\ell}}_{::k} + \sfE_{\theta_3}^{\ell-1} \diag(\bW_{:k}) \right)^{-1} {\ocU_{::k}^\ell}^T.
\]
Given these updates, the computation of $b_{\theta_3}^\ell = \mathbb{E}_{q_{\ell}(\mathcal{L})} \norm{\mathcal{L}}_{\bW*}$ necessitates a corresponding adjustment:

\[
\E_{q_\ell(\cL)}\norm{\cL}_{\bW*} = \norm{\sfE_{\cL}^\ell}_{\bW*} + \frac{n_2}{2} \sum_{k=1}^{n_3} \trace{\left( \sfE_{\theta_1}^{\ell-1} {\overline{\cD}_\tau^{\ell}}_{::k} + \sfE_{\theta_3}^{\ell-1} \diag(\bW_{:k})^{-1} \right)^{-1}}.
\]
Note the subtle yet crucial change in the trace term, ensuring consistency with the weighted norm definition.

\section{Experiments}\label{sectionexppara}
In this section,  we give experimental results to illustrate the performance of the proposed method. All the experiments are implemented using MATLAB (R2022b) on the Windows 10 platform with Intel Core
i5-1135G7 2.40 GHz and 16 GB of RAM.

\subsection{Validation on synthetic data}
Here, we generate each observation $\mathcal{X}$ in $\mathbb{R}^{n_1 \times n_2 \times n_{3}}$ by combining a low-rank tensor $\mathcal{L}_{0}$ and  a sparse tensor $\mathcal{S}_{0}$ with a Gaussian noise $\mathcal{E}_{0}$ in the  the same dimensions. The low-rank tensor $\mathcal{L}_{0}$ is derived from the t-product of two smaller tensors, namely $\mathcal{P}$ in $\mathbb{R}^{n_1 \times r \times n_{3}}$ and $\mathcal{H}$ in $\mathbb{R}^{r \times n_2 \times n_{3}}$, where $r$ is significantly smaller than $n_2$. The tubal rank of $\mathcal{L}_{0}$ does not exceed $r$. The entries of $\mathcal{P}$ are independently and identically distributed according to a Gaussian distribution $\mathcal{N}(0, 1/n_1)$, and those of $\mathcal{H}$ follow $\mathcal{N}(0, 1/n_2)$. The sparse tensor $\mathcal{S}_{0}$ has entries determined by a Bernoulli process, where each element is either $+1$ or $-1$ with a probability $\rho$, and $0$ with a probability $1-2\rho$. The entries in Gaussian noise $\mathcal{S}_{0}$ follow $\mathcal{N}(0, \sigma^2)$.

\begin{figure*}[h]
		\begin{center}
			\begin{tabular}{cccccccccc}
	            \includegraphics[width=0.3\textwidth]{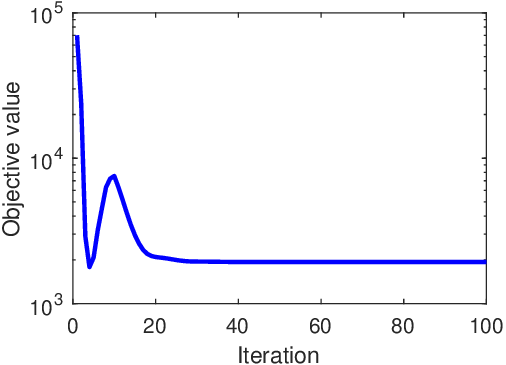}  &	
                    \includegraphics[width=0.3\textwidth]{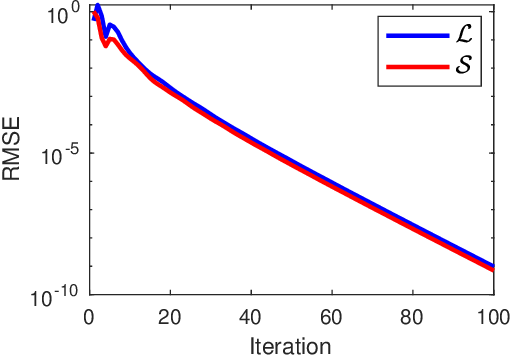}
                     &	
                    \includegraphics[width=0.3\textwidth]{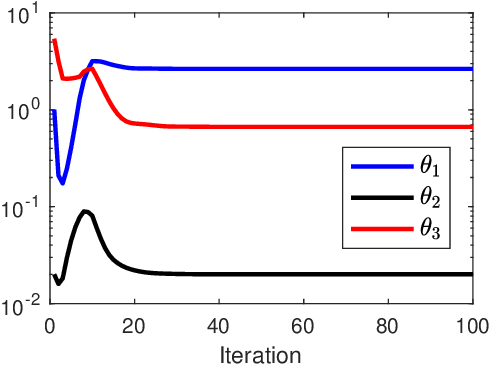}

			\end{tabular}
		\end{center}
		\caption{Empirical evidence on convergence. Left: objective function, middle: RMSE, right: parameters: $\theta_1, \theta_2$, and $\theta_3$, generated by \eqref{alg:paraalgmain} across iterations. 
		}\label{fig:conv}
	\end{figure*}

We initiate our analysis by examining the convergence properties using a third-order tensor with dimensions $40 \times 40 \times 30$. The rank parameter $r$ is set to 3, with the  parameter $\rho$ at 0.1 and the noise level $\sigma$ at $10^{-2}$. The algorithm is allowed a maximum of 100 iterations, starting with initial guesses for $\mathcal{L}$ and $\mathcal{S}$ as $\mathcal{X}$ and $\mathcal{O}$, respectively.
The convergence of the algorithm is monitored using the relative mean square error (RMSE) for $\mathcal{L}$ and $\mathcal{S}$, defined as $\frac{\|\sfE_{\cL}^{\ell}-\sfE_{\cL}^{\ell-1}\|_F}{\|\sfE_{\cL}^{\ell}\|_F}$ and $\frac{\|\sfE_{\cS}^{\ell}-\sfE_{\cS}^{\ell-1}\|_F}{\|\sfE_{\cS}^{\ell}\|_F}$, respectively. The progression of the objective values, RMSE, and parameters ($\theta_1, \theta_2, \theta_3$) is plotted across iterations in \eqref{fig:conv}. 
Due to the nonlinear and nonconvex nature of simultaneously optimizing three tensors and their associated parameters, initial fluctuations in the objective values are observed. However, after approximately ten iterations, the objective values begin to decrease steadily and achieve convergence by the 30th iteration. The parameter values similarly stabilize within these iterations. Both RMSE metrics show a sharp decline, reaching as low as $10^{-4}$ by the 30th iteration.
Given these observations, we establish a stopping criterion where the algorithm terminates when RMSE falls below $10^{-4}$ or when 50 iterations are reached, whichever occurs first. This criterion ensures efficient and effective convergence to an optimal solution within a reasonable number of iterations.

\begin{figure*}[h]
		\begin{center}
			\begin{tabular}{cccccccccc}
	            \includegraphics[width=0.3\textwidth]{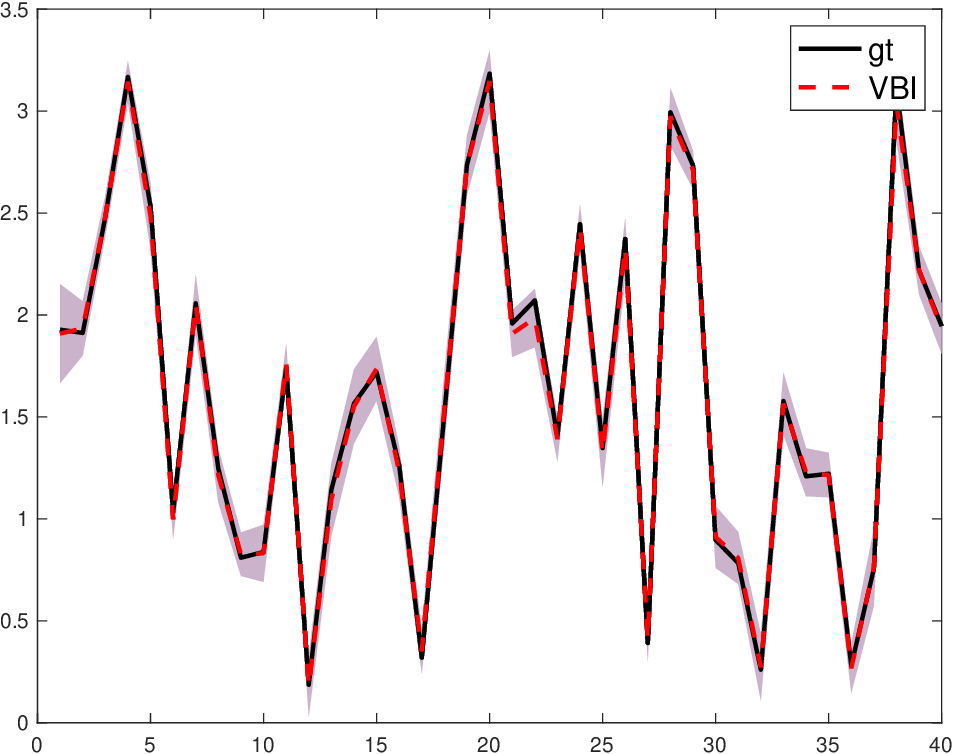}  &	
                    \includegraphics[width=0.3\textwidth]{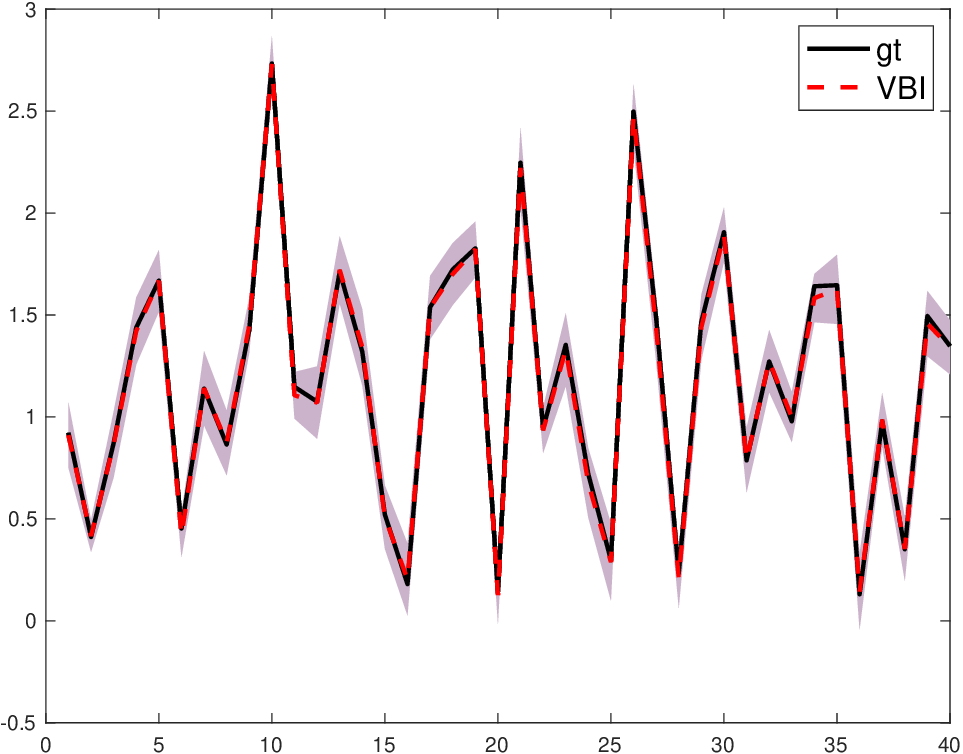}
                     &	
                    \includegraphics[width=0.3\textwidth]{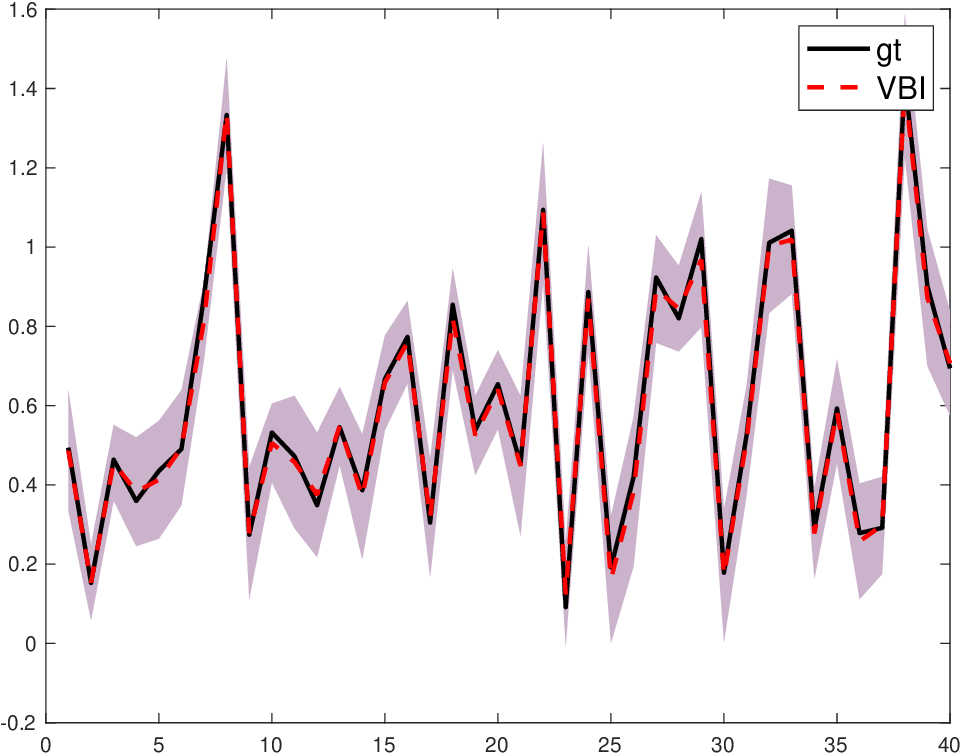}

			\end{tabular}
		\end{center}
		\caption{Uncertainty quantification: recovery of $\ocL_{:ij}$ with 99.73\% credible interval (shaded area) where $i = 20, j = 5, 10, 25$. }\label{fig:uncertain}
	\end{figure*}

Here, we further evaluate the uncertainty quantification performance of our Variational Bayesian Inference (VBI) algorithm using the same simulated tensor as previously described. \eqref{fig:uncertain}.  presents the mean estimates and 99.73\% credible intervals for the recovery of tensor filter $\ocL_{:ij}$ with $i = 20, j = 5, 15, 20$. 
The mean values consistently align with the ground truth across all fibers, while remarkably narrow credible intervals (indicated by minimal shading) demonstrate the high precision of our method. This precision is further corroborated by the low parameter standard deviations.

As part of a proof-of-concept study, we employ a partial sum of the tubal nuclear norm \cite{jiang2020multi} as a representative example for a weighted TNN in our numerical experiments. We aim to compare our proposed algorithms, VBI$_{\rm TNN}$ and VBI$_{\rm PSTNN}$, against two established methods in tensor rank approximation: TNN \cite{lu2019tensor} and PSTNN \cite{jiang2020multi}. 
For this comparative analysis, we set the noise levels $\sigma$ at $10^{-3}$, $10^{-2}$, and $10^{-1}$, the rank $r$ at 3 and 5, and the parameter $\rho$ at 0.01 and 0.1. We assess the performance of these methods by calculating the relative square error between the recovered tensors, $\hat{\mathcal{L}}$ and $\hat{\cS}$, and the ground-truth tensors, $\mathcal{L}_{\rm GT}$ and $\cS_{\rm GT}$. These errors are quantified as follows: ${\rm error}_\cL = \frac{\|\hat{\mathcal{L}} - \mathcal{L}_{\rm GT}\|_F}{\|\mathcal{L}_{\rm GT}\|_F}$ for the low-rank component and ${\rm error}_\cS = \frac{\|\hat{\mathcal{S}} - \mathcal{S}_{\rm GT}\|_F}{\|\mathcal{S}_{\rm GT}\|_F}$ for the sparse component.

As shown in \eqref{tab:val},  VBI$_{\rm TNN}$ generally outperforms TNN across most tested scenarios, while VBI$_{\rm PSTNN}$ is better than PSTNN.  Moreover, VBI$_{\rm PSTNN}$ consistently delivers the best performance, indicating its superior ability to recover both the low-rank and sparse components of tensors under various noise and rank conditions. This comparative analysis underscores the effectiveness of our proposed methods, particularly VBI$_{\rm PSTNN}$, in handling complex tensor decomposition with higher accuracy and robustness against noise.

\begin{table}[htbp]
\footnotesize\tabcolsep=1.5mm
\caption{Recovery results on the synthetic datasets with different settings. }\label{tab:val}
\begin{center}
  \begin{tabular}{|c|c|c|c|c|c|c|c|c|c|c|c|c|} \hline
   \multicolumn{3}{|c|}{Method} & \multicolumn{2}{|c|}{TNN} &   \multicolumn{2}{|c|}{VBI$_{\rm TNN}$} & \multicolumn{2}{|c|}{PSTNN} &   \multicolumn{2}{|c|}{VBI$_{\rm PSTNN}$}\\ \hline
   $\sigma$ & $r$ & $\rho$ & error$_\mathcal{L}$ & error$_\mathcal{S}$ & error$_\mathcal{L}$ & error$_\mathcal{S}$  & error$_\mathcal{L}$ & error$_\mathcal{S}$ & error$_\mathcal{L}$ & error$_\mathcal{S}$ 
   \\ \hline
    \multirow{4}{*}{$10^{-3}$} & \multirow{2}{*}{3} & 0.01 & 0.0029 & 0.0075 & 0.0025 & 0.0056 & 0.0028 & 0.0064 & {\bf 0.0023} &  {\bf 0.0052} \\
    & & 0.1 & 0.0034 & 0.0027 & 0.0032 & 0.0025 & 0.0033 & 0.0024 & {\bf  0.0029}  & {\bf  0.0023} \\ \cline{2-11}
    & \multirow{2}{*}{5} & 0.01 & 0.0026 & 0.0083 & 0.0025 & 0.0063  & 0.0024 & 0.0070 & {\bf  0.0022} & {\bf  0.0058} \\
    & & 0.1 & 0.0033 & 0.0033 & 0.0036 & 0.0032 & {\bf 0.0030} & {\bf  0.0028} & 0.0031 & 0.0029  \\ \hline
       \multirow{4}{*}{$10^{-2}$} & \multirow{2}{*}{3} & 0.01 & 0.0286 & 0.0738 & 0.0248 & 0.0556  & 0.0276 & 0.0638 &  {\bf 0.0230}  & {\bf 0.0523} \\
    & & 0.1 & 0.0344 & 0.0274 & 0.0302 & 0.0238  & 0.0325 & 0.0240 & {\bf 0.0275} & {\bf 0.0223} \\ \cline{2-11}
    & \multirow{2}{*}{5} & 0.01 & 0.0257 & 0.0820 & 0.0242 & 0.0620 & 0.0240 & 0.0700 & {\bf 0.0219} & {\bf 0.0576}  \\
    & & 0.1 & 0.0331 & 0.0329 & 0.0322 & 0.0294  & 0.0298 & 0.0281 & {\bf 0.0281} & {\bf 0.0267} \\ \hline
       \multirow{4}{*}{$10^{-1}$} & \multirow{2}{*}{3} & 0.01 & 0.2744 & 0.7227 & 0.2317 & 0.5435  & 0.2769 & 0.6398 & {\bf 0.2255} & {\bf 0.5195} \\
    & & 0.1 & 0.3222 & 0.2623 & 0.2730 & 0.2262 & 0.3264 & 0.2410 & {\bf 0.2661} & {\bf 0.2187} \\ \cline{2-11}
    & \multirow{2}{*}{5} & 0.01 & 0.2392 & 0.7841 & 0.2201 & 0.5921  & 0.2346 & 0.6896 & {\bf 0.2077} & {\bf 0.5620} \\
    & & 0.1 & 0.2903 & 0.2961 & 0.2692 & 0.2589 & 0.2864 & 0.2705 & {\bf 0.2543} & {\bf 0.2484}  \\ \hline
  \end{tabular}
\end{center}
\end{table}

\subsection{Image denoising}
In this section,  we evaluate the performance of the proposed method on image denoising. The peak signal-to-noise ratio (PSNR) \cite{lu2019tensor} and the structural similarity index (SSIM) \cite{wang2004image} are used to evaluate the recovery performance quantitatively. 

\subsubsection{Image with sparse noise}
We conduct experiments on four images: ``house", ``moto", ``face", and ``hat". In this study, we model the clean images as the low-rank component and random corruptions as sparse outliers. Each image is corrupted by setting 10 percent of the pixels to random values ranging from 0 to 255, with the locations of these distortions unspecified. We compare our proposed method with several existing techniques, including LRTV \cite{he2015total}, $S_{wp}(0.9)$ \cite{YANG2022108311}, BTRTF \cite{zhou2019bayesian}, TNN \cite{lu2019tensor}, and PSTNN \cite{jiang2020multi}, using the original implementations provided by the respective authors.
Given the absence of Gaussian noise in this task, the parameter $\theta_1$ is set to a high value of 100 to accommodate this condition, while $\theta_2$ and $\theta_3$ are set to 1. The truncation parameter $K$ for VBI$_{\rm PSTNN}$ is consistently set at 50 across all cases.

\begin{figure*}[h]
\scriptsize
\centering
\setlength{\tabcolsep}{0.2em}
\begin{tabular}{cccccccccc}
{   Clean} & { LRTV}  & { $S_{wp}$(0.9)} & { BTRTF }  & { TNN}  & { PSTNN} & { VBI$_{\rm TNN} $ }  & { VBI$_{\rm PSTNN}$ }\\
\includegraphics[width=0.115\textwidth]{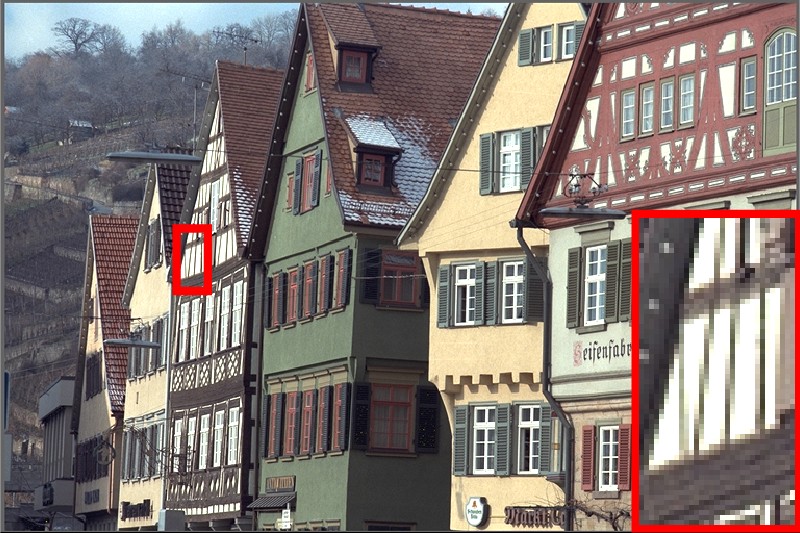}   &
\includegraphics[width=0.115\textwidth]{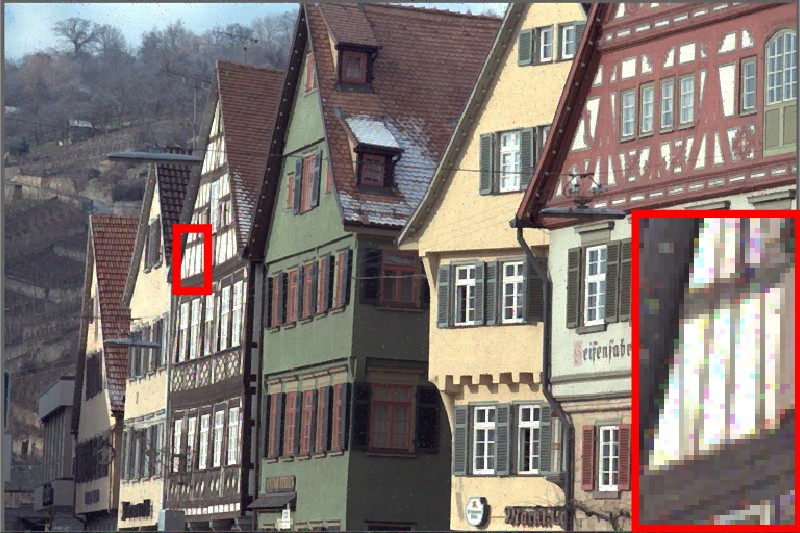}   &
\includegraphics[width=0.115\textwidth]{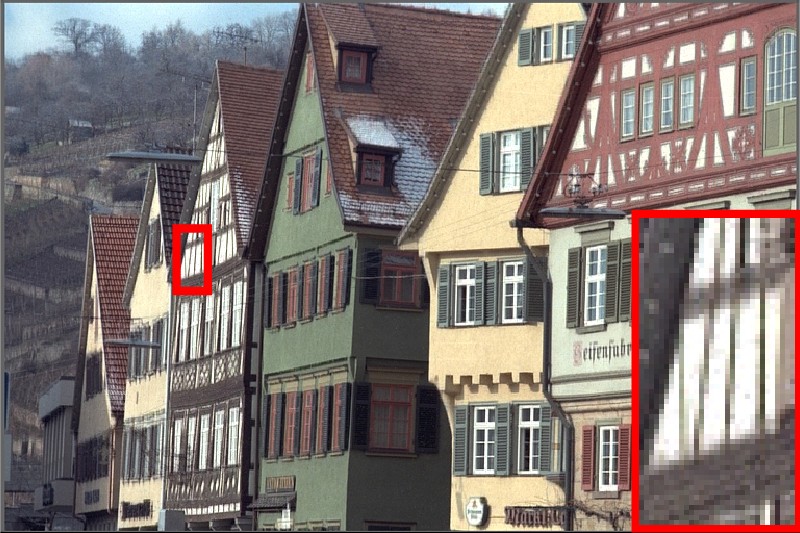}  &
\includegraphics[width=0.115\textwidth]{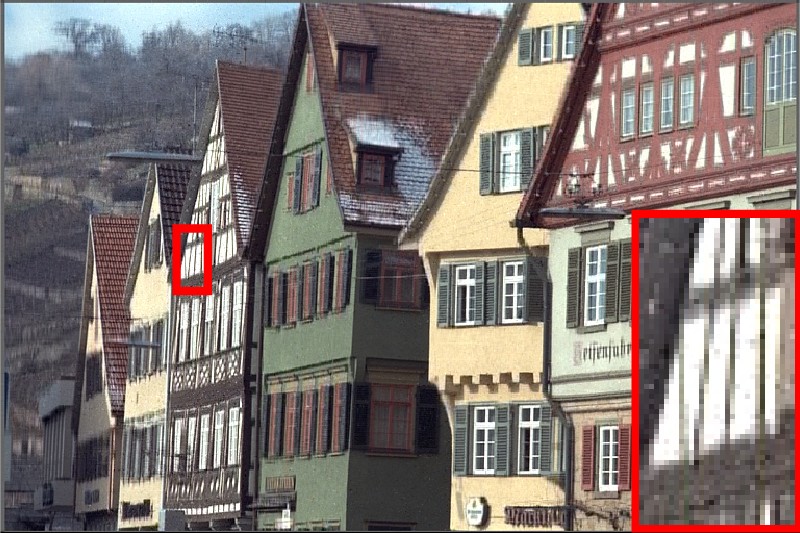}  & 
\includegraphics[width=0.115\textwidth]{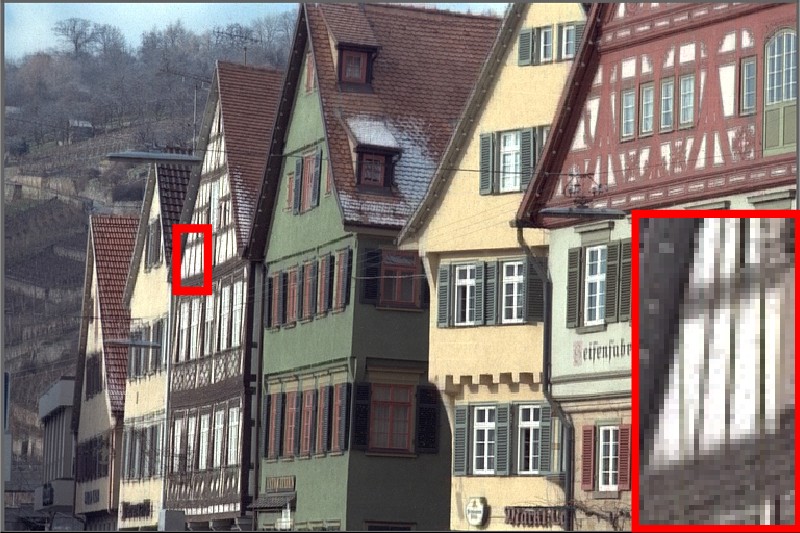}   &
\includegraphics[width=0.115\textwidth]{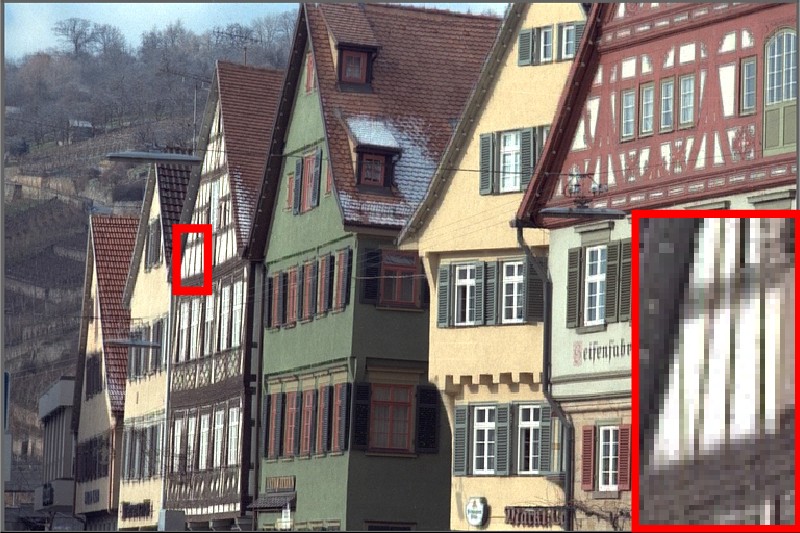}   &
\includegraphics[width=0.115\textwidth]{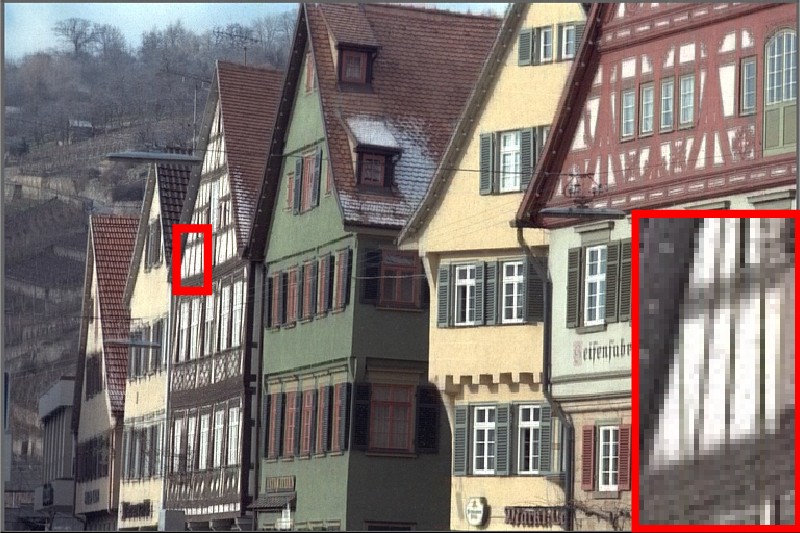}  & 
\includegraphics[width=0.115\textwidth]{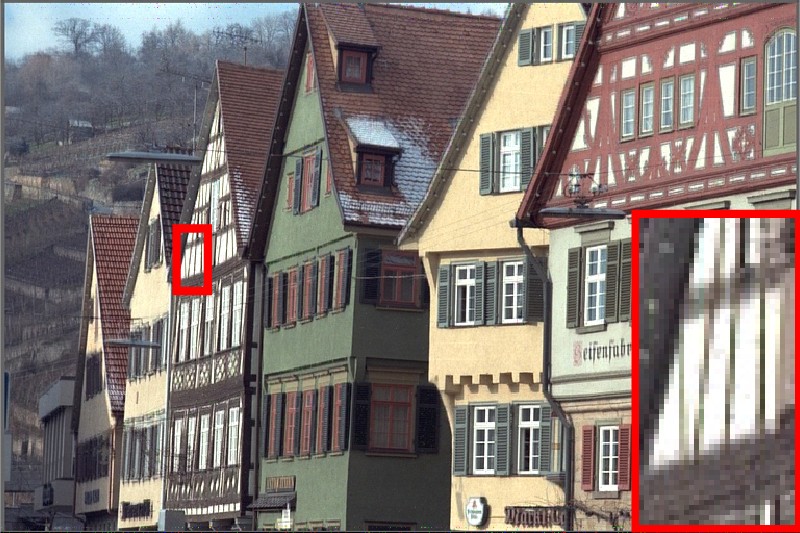}\\

\includegraphics[width=0.115\textwidth]{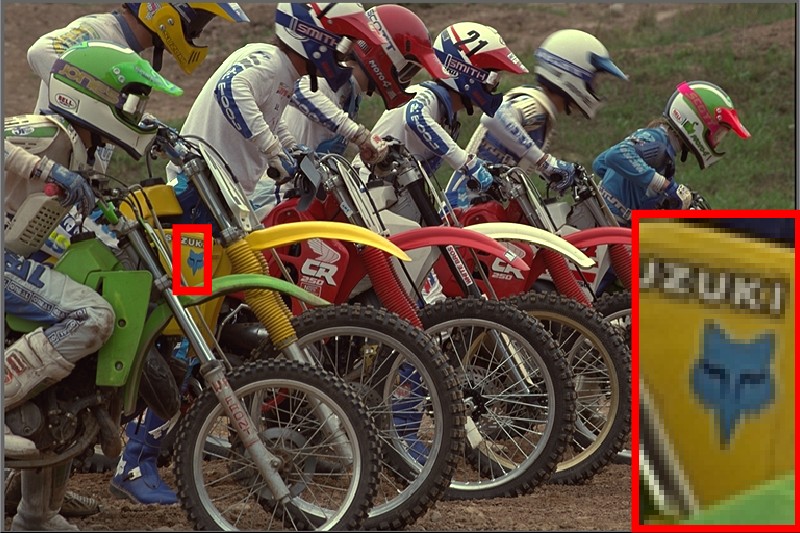}   &
\includegraphics[width=0.115\textwidth]{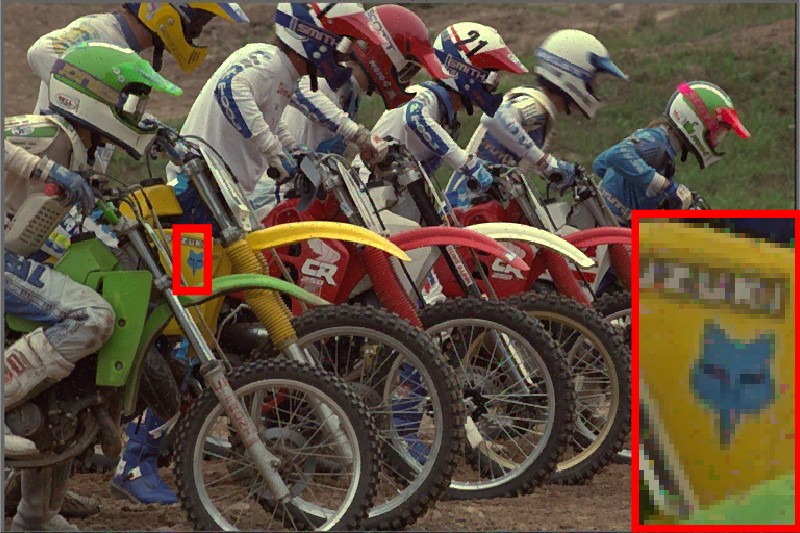}   &
\includegraphics[width=0.115\textwidth]{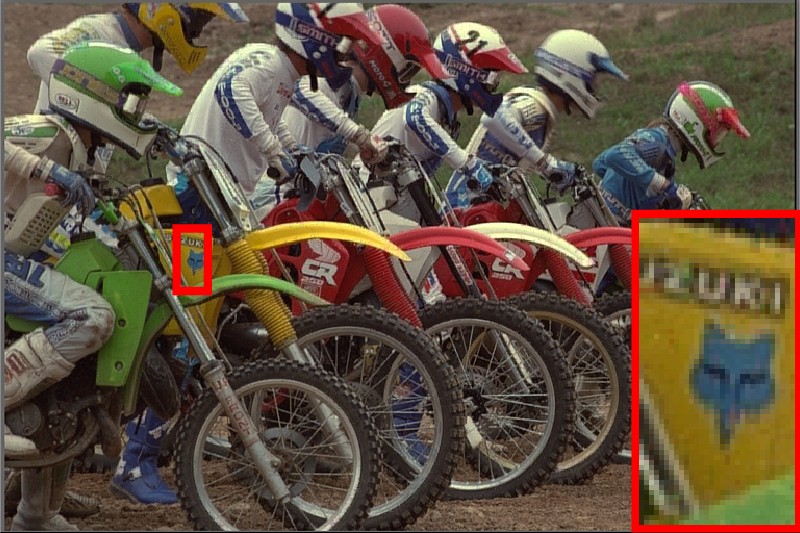}  &
\includegraphics[width=0.115\textwidth]{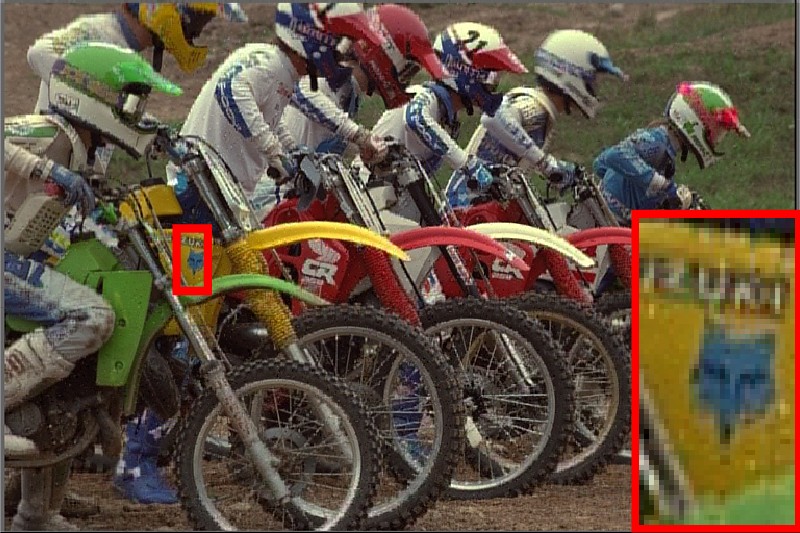}  & 
\includegraphics[width=0.115\textwidth]{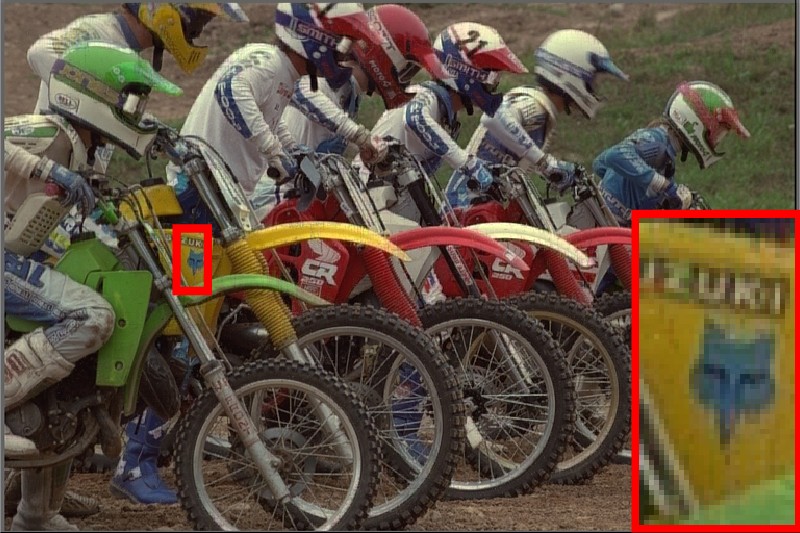}   &
\includegraphics[width=0.115\textwidth]{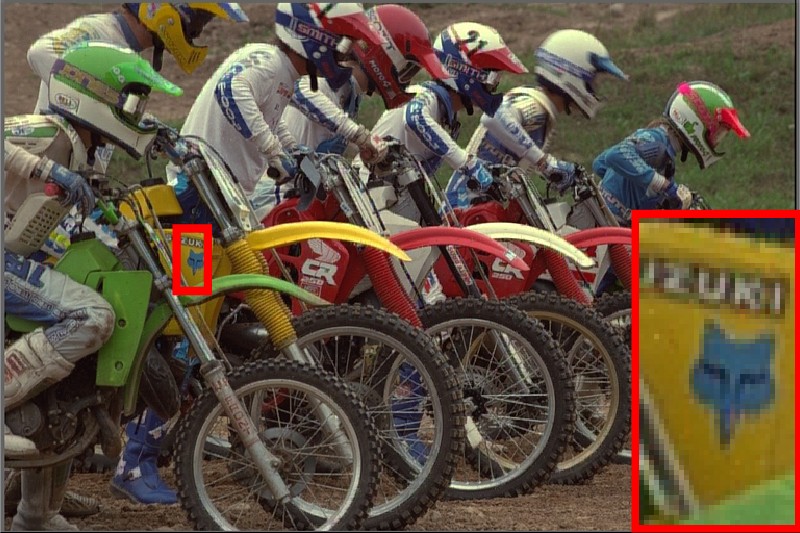}   &
\includegraphics[width=0.115\textwidth]{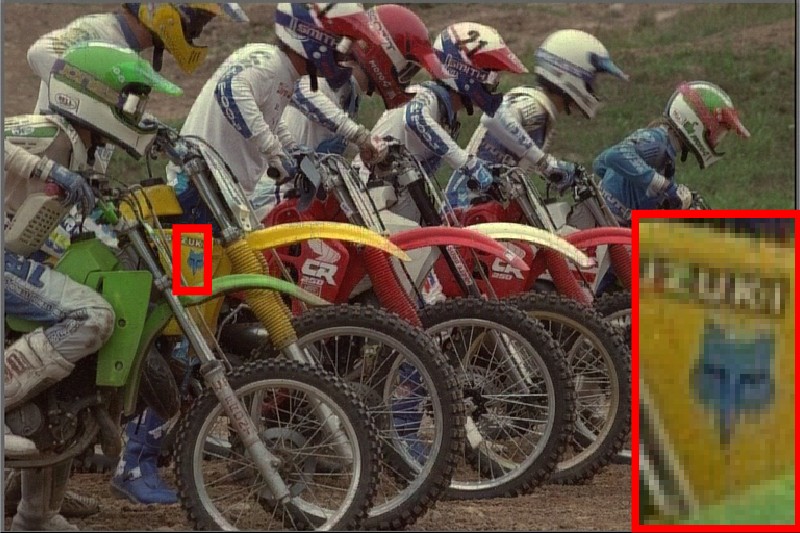}  & 
\includegraphics[width=0.115\textwidth]{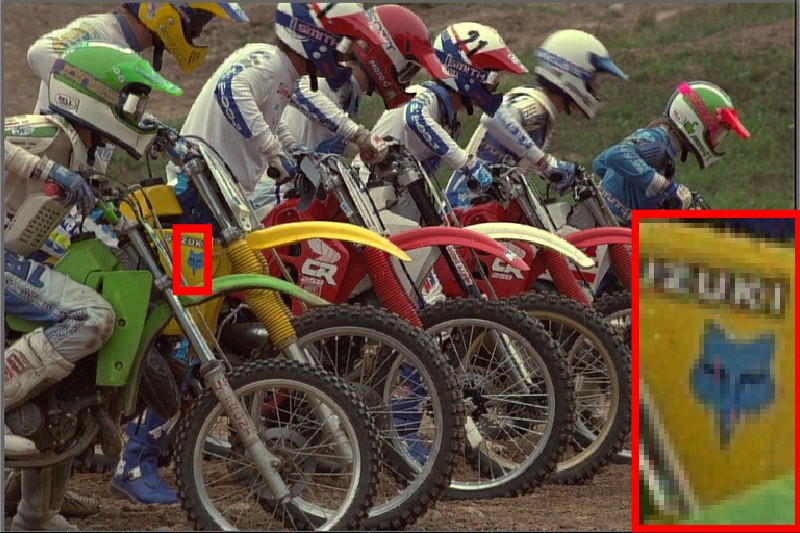}\\

\includegraphics[width=0.115\textwidth]{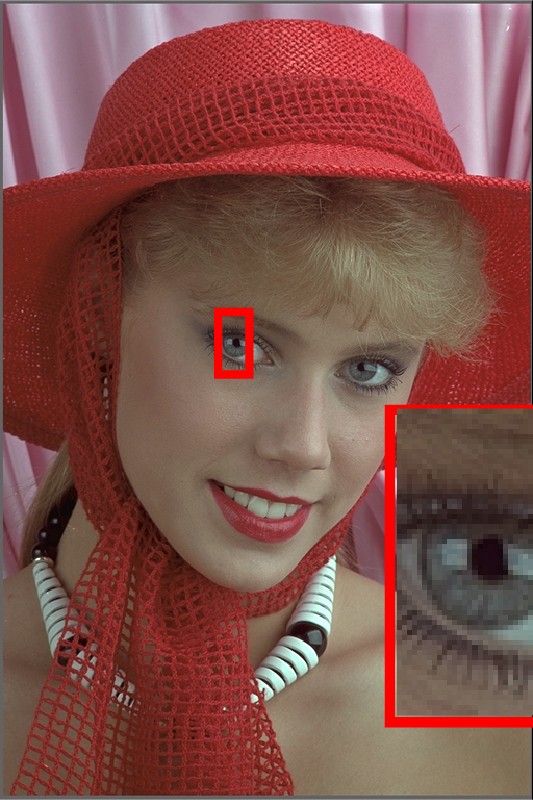}   &
\includegraphics[width=0.115\textwidth]{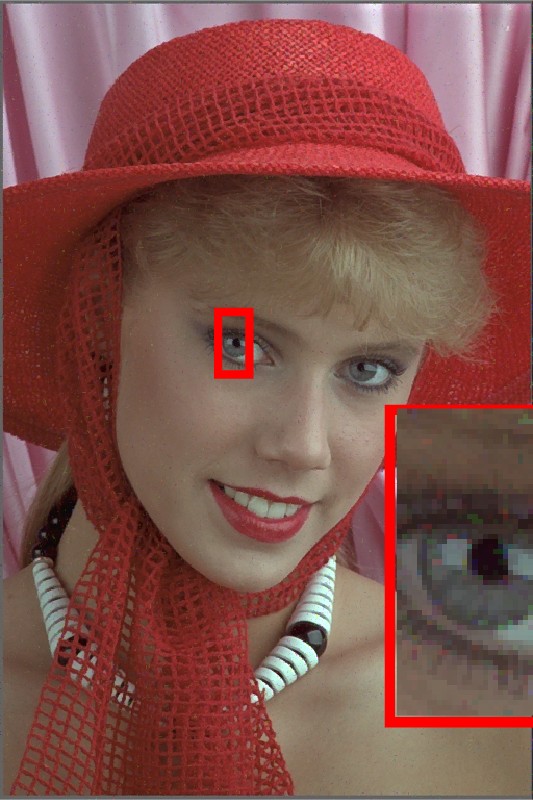}   &
\includegraphics[width=0.115\textwidth]{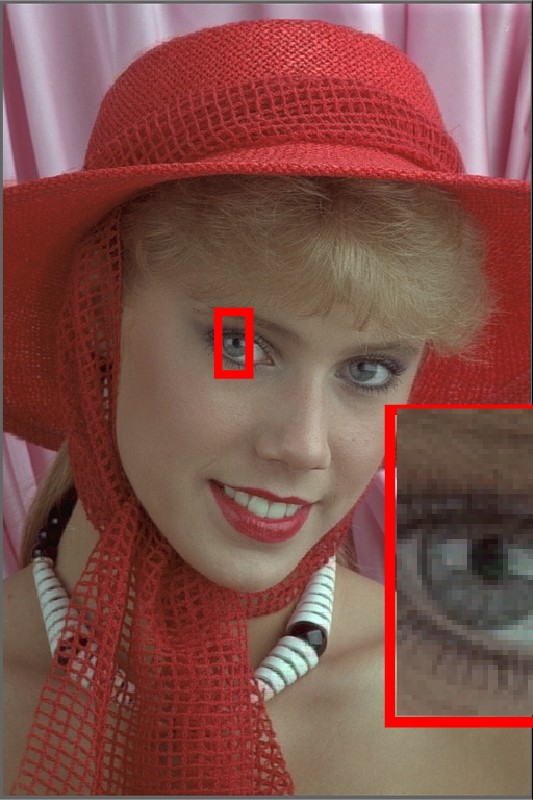}  &
\includegraphics[width=0.115\textwidth]{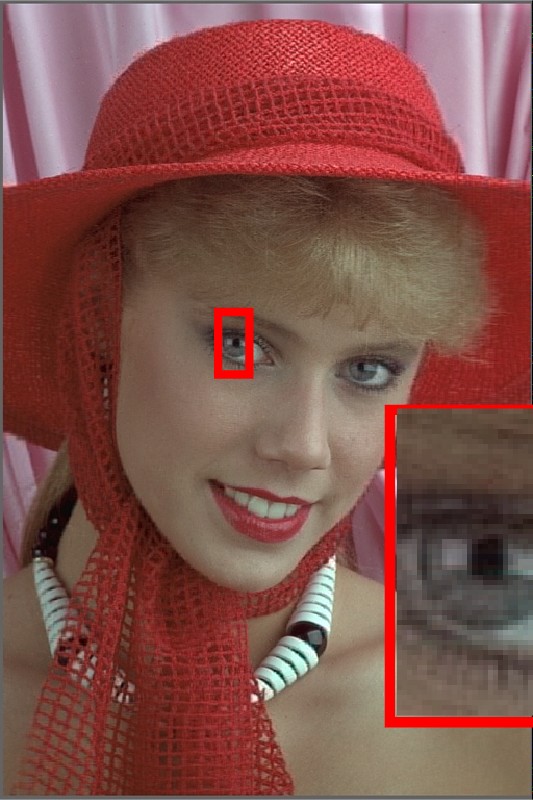}  & 
\includegraphics[width=0.115\textwidth]{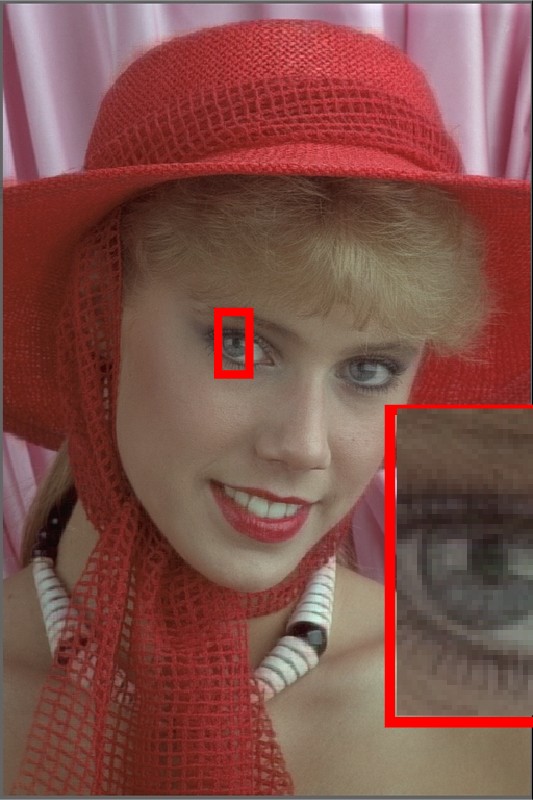}   &
\includegraphics[width=0.115\textwidth]{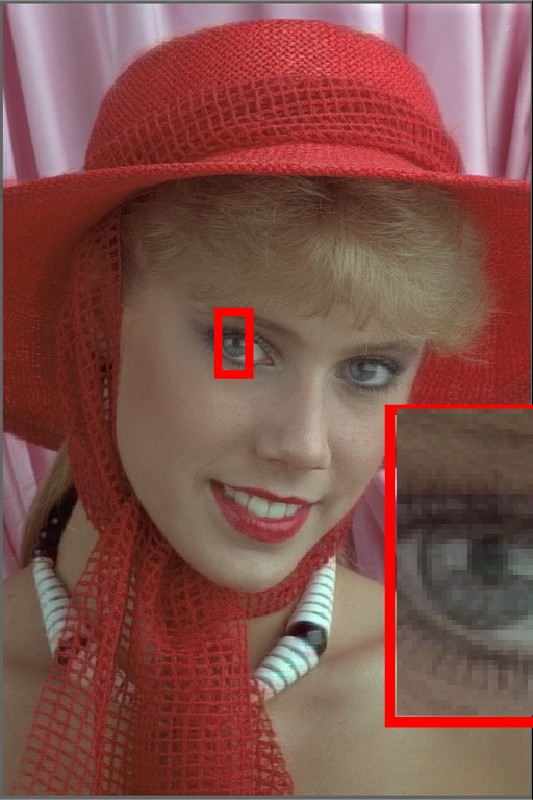}   &
\includegraphics[width=0.115\textwidth]{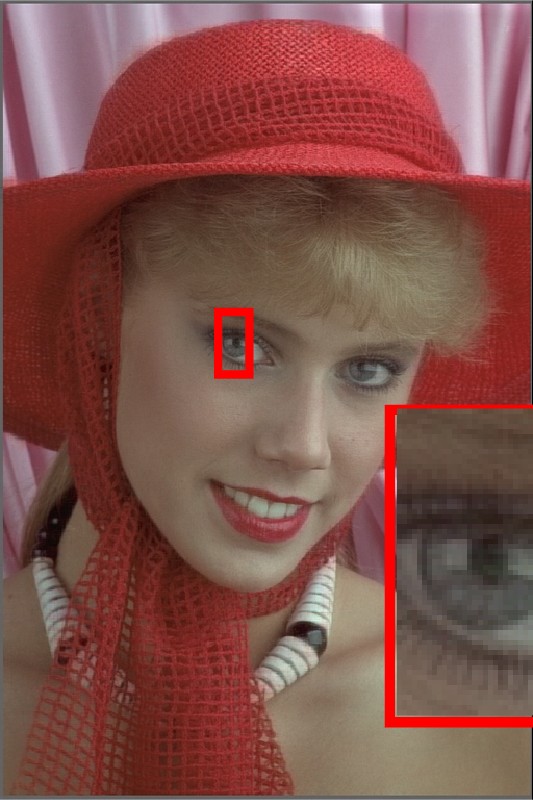}   & 
\includegraphics[width=0.115\textwidth]{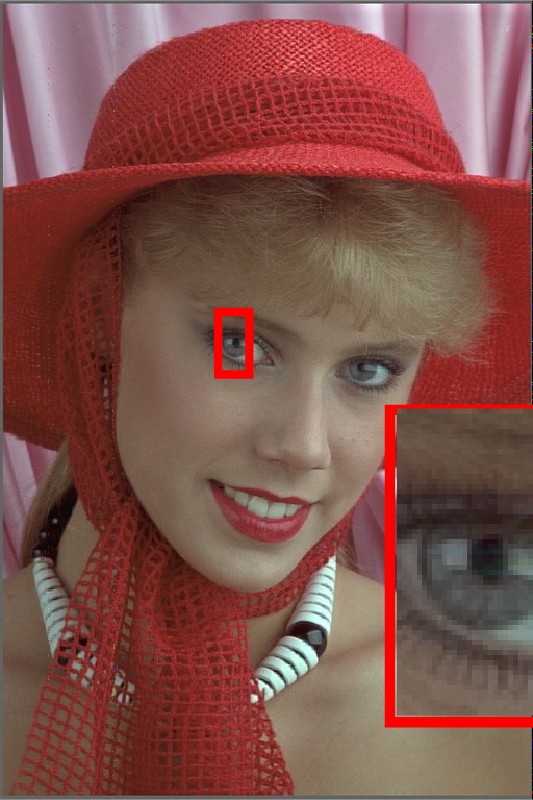}\\

\includegraphics[width=0.115\textwidth]{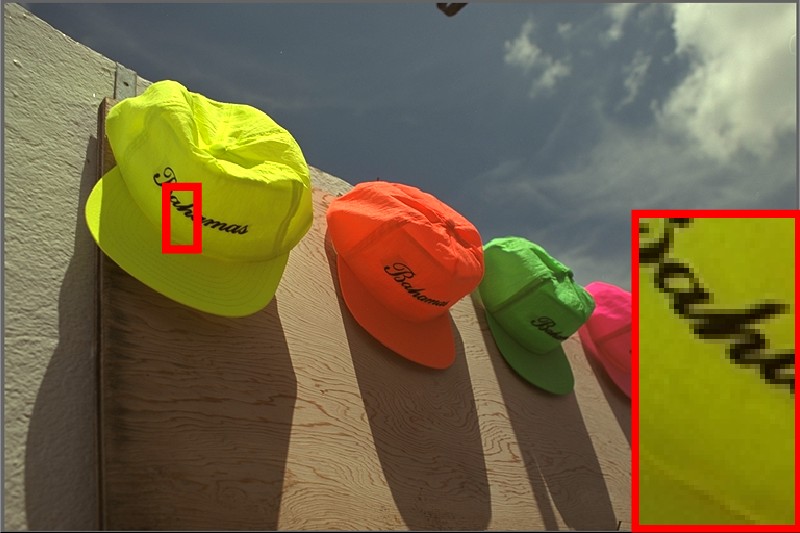}   &
\includegraphics[width=0.115\textwidth]{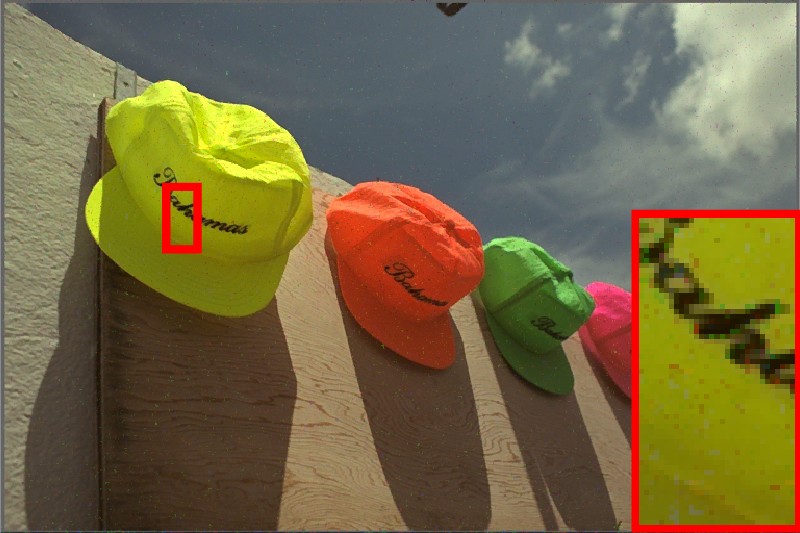}   &
\includegraphics[width=0.115\textwidth]{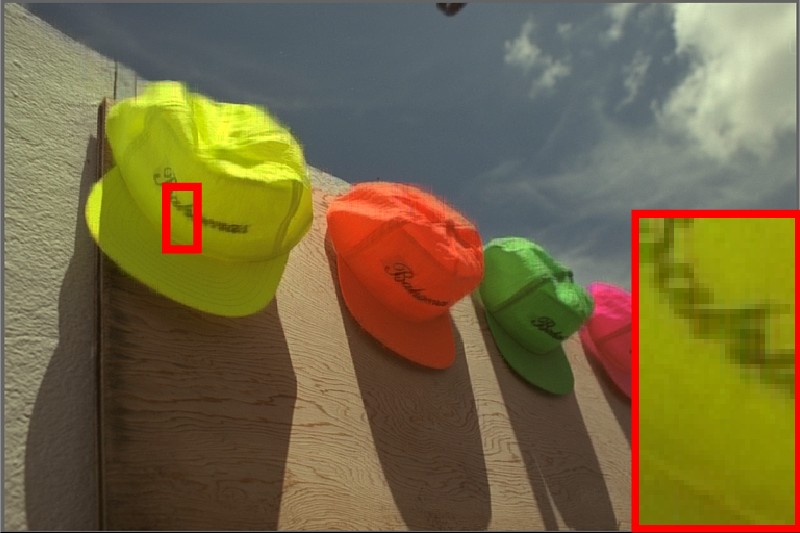}  &
\includegraphics[width=0.115\textwidth]{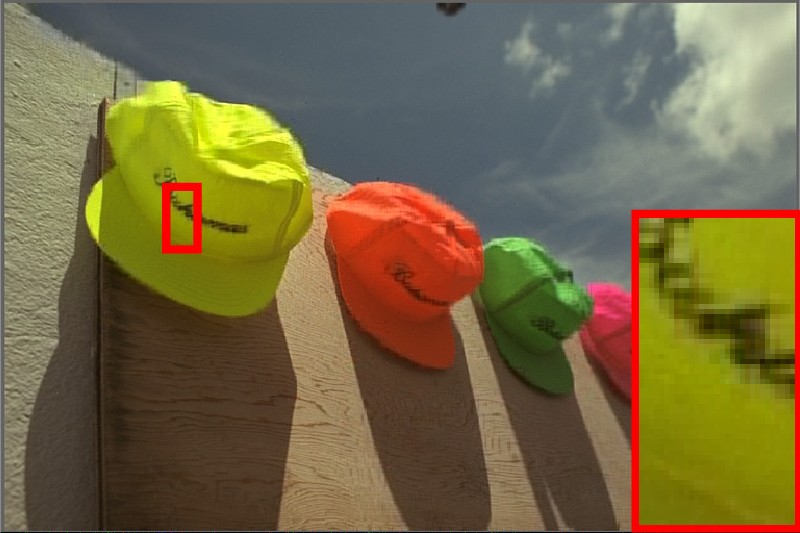}  & 
\includegraphics[width=0.115\textwidth]{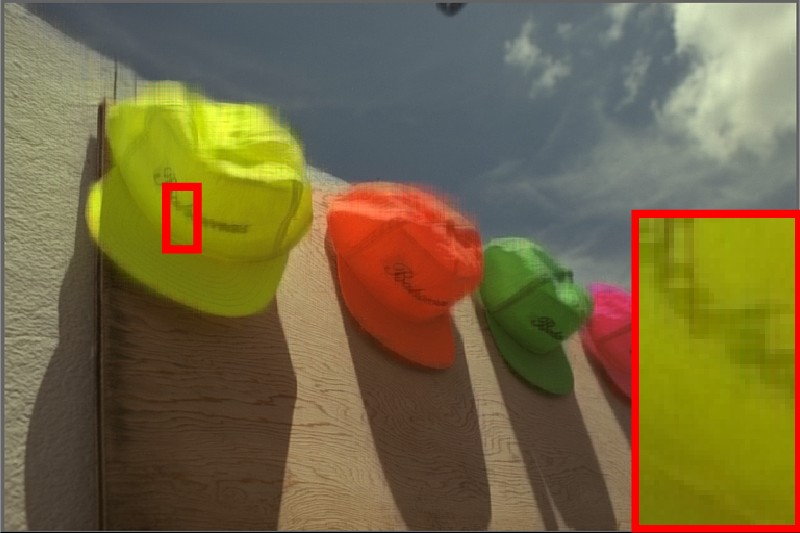}   &
\includegraphics[width=0.115\textwidth]{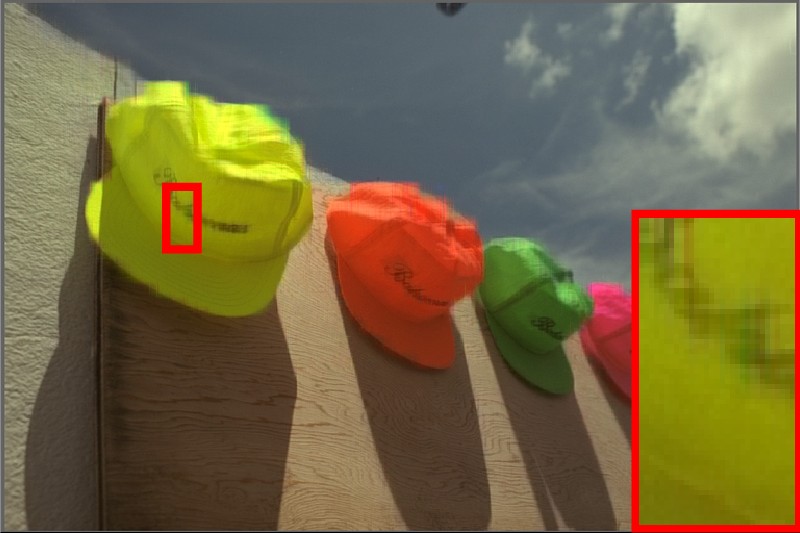}   &
\includegraphics[width=0.115\textwidth]{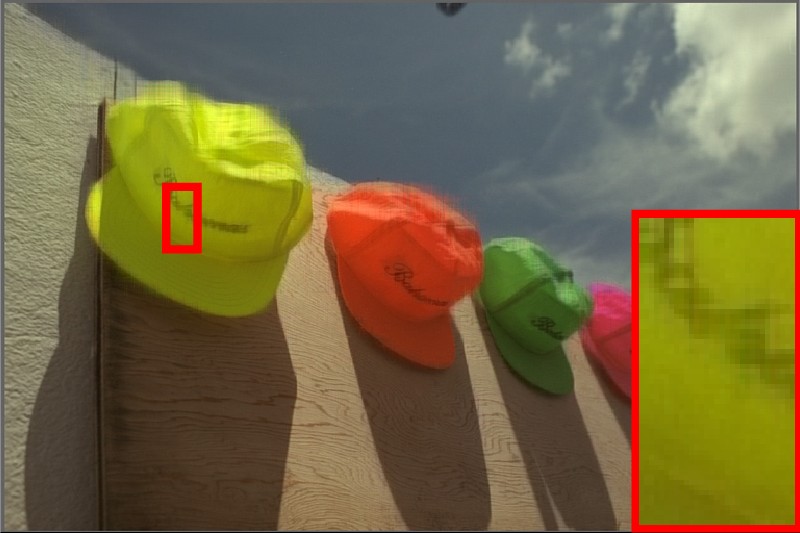}   & 
\includegraphics[width=0.115\textwidth]{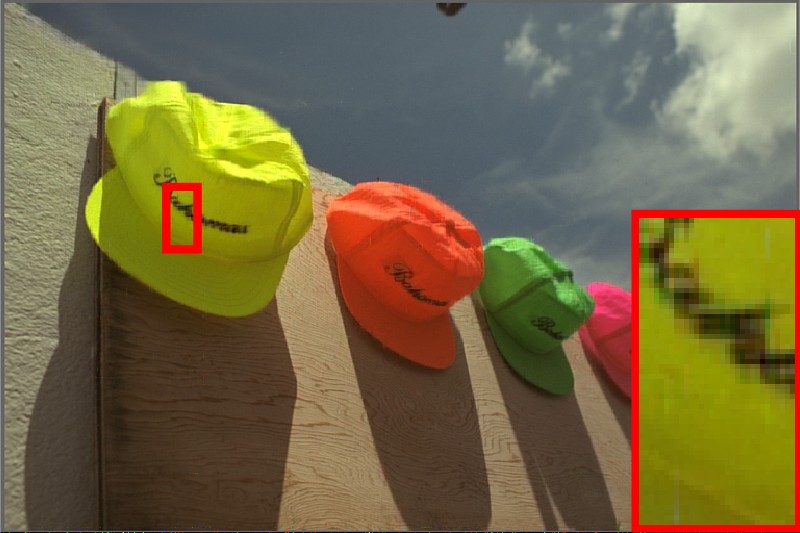}
\end{tabular}

\caption{Comparison of color image Gaussian noise removal performance on four examples.
}\label{fig:gauss}
\end{figure*}

Quantitative evaluations based on PSNR and SSIM are presented in \eqref{tab:gauss}, and the corresponding restored images are displayed in \eqref{fig:gauss}. Our observations indicate that VBI$_{\rm PSTNN}$ consistently outperforms the other methods in terms of PSNR, achieving at least a 0.5 improvement and matching the best-performing methods in SSIM values. Additionally, the restoration of the ``hat" image by VBI$_{\rm PSTNN}$ and BTRTF shows significantly clearer text compared to other methods. However, some artifacts are noted in the ``moto" image restored by BTRTF. In contrast, our method exhibits fewer artifacts across all cases.

\begin{table}[htbp]
\footnotesize\tabcolsep=1.5mm
\caption{Quantitative comparisons of sparse noise removal results obtained by different methods  }\label{tab:gauss}
\begin{center}
  \begin{tabular}{|c|c|c|c|c|c|c|c|c|c|c|c|} \hline
   Data & Index & LRTV  & $S_{wp}(0.9)$ & BTRTF &  TNN & PSTNN  &  VBI$_{\rm TNN}$ & VBI$_{\rm PSTNN}$ \\ \hline
    \multirow{2}{*}{house} & PSNR &  26.167 & 28.028 & 25.930 & 27.030 & 27.522 & 26.878 & {\bf 28.565} \\
     & SSIM & 0.9517 & 0.9717 & 0.9374 & 0.9655 & 0.9691 & 0.9596 & {\bf 0.9741}    \\ \hline 
    \multirow{2}{*}{moto} & PSNR & 27.617 & 28.003 & 24.871 & 26.373 & 27.724 & 25.945 & {\bf 28.781} \\ 
     & SSIM & 0.9590 & 0.9702 & 0.9130 & 0.9554 & 0.9672 & 0.9440 & {\bf 0.9719}\\ \hline 
         \multirow{2}{*}{face} & PSNR & 32.524 & 34.061 & 32.500 & 30.770 & 31.543 & 30.704 & {\bf 34.150} \\ 
     & SSIM & 0.9529 & {\bf 0.9759}  & 0.9405 & 0.9509 & 0.9557 & 0.9475 & 0.9694\\ \hline
     \multirow{2}{*}{hat} & PSNR & 32.626 & 32.787 & 32.558 & 29.453 & 30.895 & 29.755 & {\bf 33.478} \\ 
     & SSIM & 0.9435 & {\bf 0.9750}  & 0.9581 & 0.9473 & 0.9558 & 0.9516 & 0.9735 \\ \hline
     \multirow{2}{*}{mean} & PSNR & 29.733  & 30.720 & 28.965 & 28.407 & 29.421 & 28.321 & {\bf 31.244} \\ 
     & SSIM & 0.9518 & {\bf 0.9732}  & 0.9375 & 0.9548 & 0.9620 & 0.9507 & 0.9722 \\ \hline
  \end{tabular}
\end{center}
\end{table}

\begin{table}[htbp]
\footnotesize\tabcolsep=1.5mm
\caption{Quantitative comparisons of mixed noise removal results obtained by different methods  }\label{tab:mixed}
\begin{center}
  \begin{tabular}{|c|c|c|c|c|c|c|c|c|c|c|c|} \hline
   Data & Index & 3DTNN  & $S_{wp}(0.9)$ & BTRTF &  TNN & PSTNN  &  VBI$_{\rm TNN}$ &  VBI$_{\rm PSTNN}$ \\ \hline
    \multirow{2}{*}{kid} & PSNR& 26.670 & 31.806 & 32.071 & 28.691 & 29.542 & 29.446 & {\bf 32.802}  \\
     & SSIM & 0.9364 & {\bf 0.9752} & 0.9593 & 0.9487 & 0.9558 & 0.9521 & 0.9720   \\ \hline 
       \multirow{2}{*}{house} & PSNR& 27.448 & 32.302 & 30.791 & 29.765 & 30.459 & 29.862 & {\bf 32.496} \\
     & SSIM & 0.9292 & {\bf 0.9708} & 0.9370 & 0.9474 & 0.9532 & 0.9414 & 0.9659   \\ \hline 
            \multirow{2}{*}{river} & PSNR& 24.606 & 26.388 & 23.818 & 25.985 & 26.439 & 25.367 & {\bf 26.968} \\
     & SSIM & 0.9319 & 0.9471 & 0.8606 & 0.9466 & {\bf 0.9515} & 0.9291 & 0.9504   \\ \hline 
     \multirow{2}{*}{hat} & PSNR & 28.017 & 32.771 & 32.553 & 29.449 & 30.891 & 29.753 & {\bf 33.463} \\
     & SSIM  & 0.9359 & 0.9747 & 0.9581 & 0.9471 & 0.9555 & 0.9514 & {\bf 0.9733}   \\ \hline 
         \multirow{2}{*}{mean} & PSNR & 26.685 & 30.817 & 29.808 & 28.473 & 29.333 & 28.607 & {\bf 31.432} \\
     & SSIM  & 0.9334 & {\bf 0.9670} &  0.9288 &  0.9475 & 0.9540 &  0.9436 &  0.9654   \\ \hline 
  \end{tabular}
\end{center}
\end{table}


\subsubsection{Image with mixed  noise}
In this subsection, we perform experiments on four distinct images: ``kid", ``house", ``river", and ``hat".  Initially, each image is corrupted with sparse noise, following the procedure of our previous experiment. Subsequently, we introduce Gaussian noise to each pixel, modeled by the distribution $\mathcal{N}(0, 10^{-3})$. The resultant observation, represented mathematically by $\mathcal{X} = \mathcal{L} + \mathcal{S} + \mathcal{E}$, consists of the real image $\mathcal{L}$, augmented by sparse noise $\mathcal{S}$ and Gaussian noise $\mathcal{E}$. To verify that our method's effectiveness is robust to initial conditions, we set the initial values of $\theta_1$ to 100, and $\theta_2$ and $\theta_3$ to 1, as the same as the ones used in the sparse noise-only scenario.

\begin{figure*}[h]
\centering
\scriptsize
\setlength{\tabcolsep}{0.2em}
\begin{tabular}{cccccccccc}
{ Clean} & { 3DTNN}  & { $S_{wp}$(0.9)} & { BTRTF }  & { TNN}  & { PSTNN} & { VBI$_{\rm TNN} $ }  & { VBI$_{\rm PSTNN}$ }\\
\includegraphics[width=0.115\textwidth]{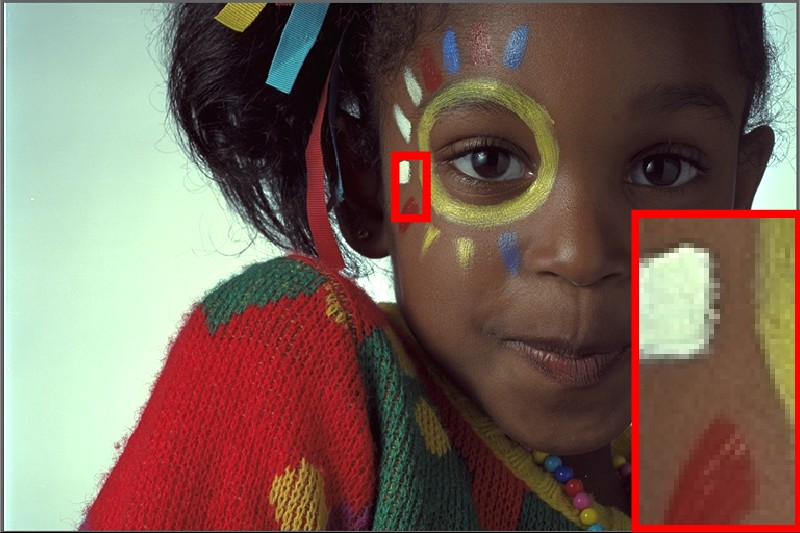}  &
\includegraphics[width=0.115\textwidth]{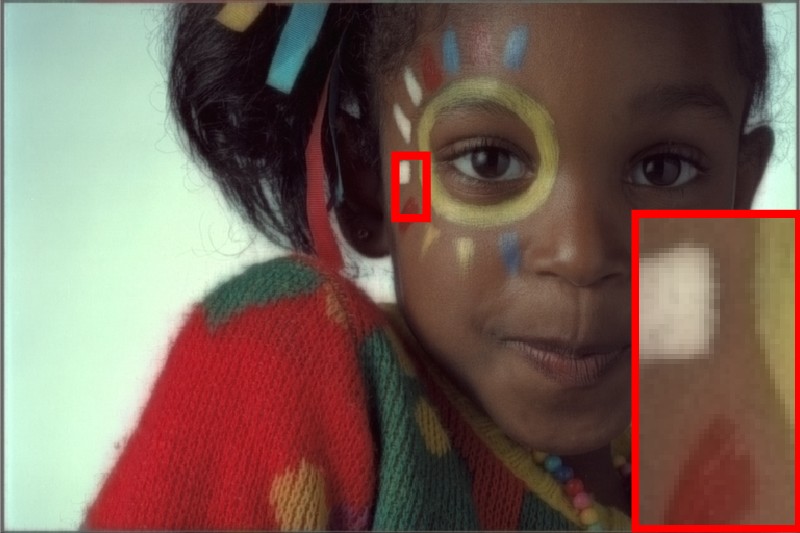}  &
\includegraphics[width=0.115\textwidth]{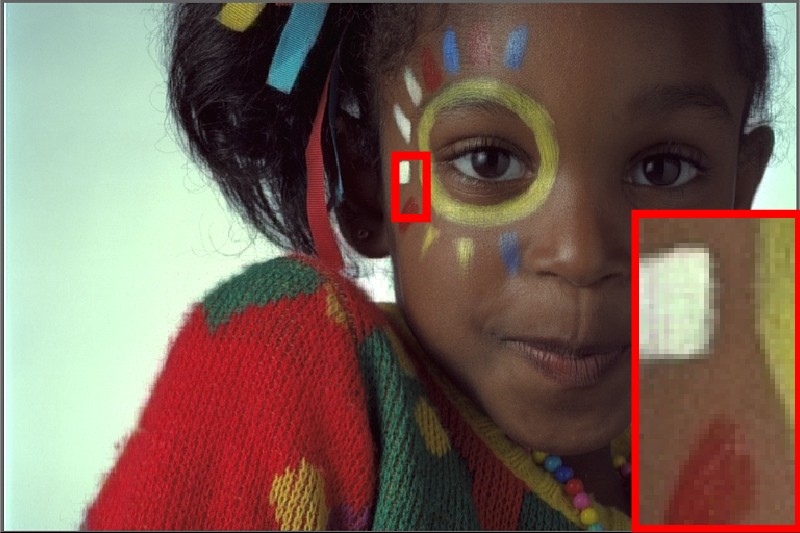} &
\includegraphics[width=0.115\textwidth]{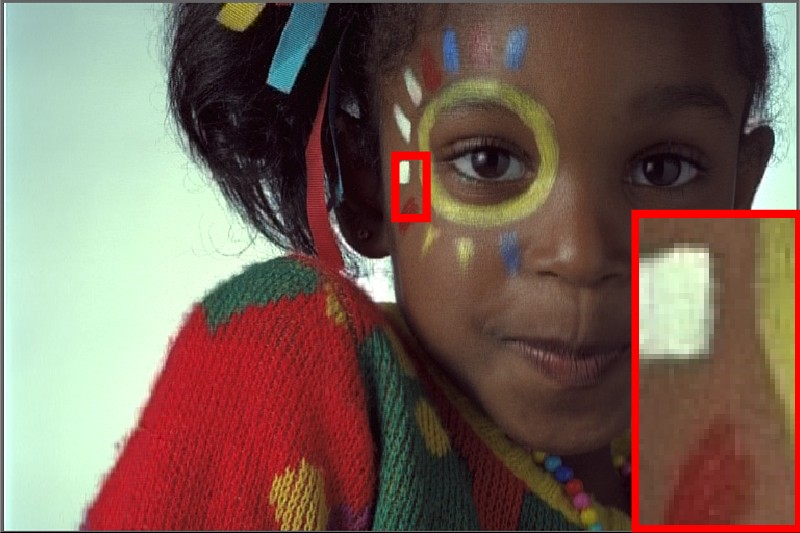} & 
\includegraphics[width=0.115\textwidth]{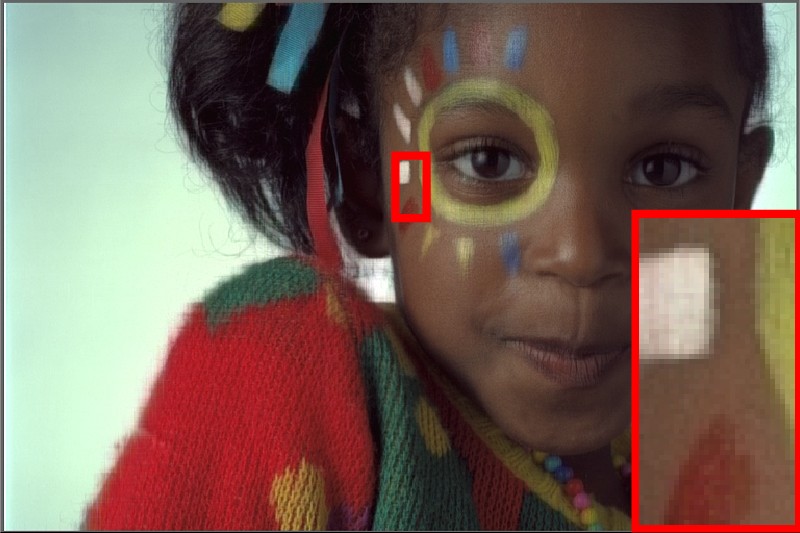}  &
\includegraphics[width=0.115\textwidth]{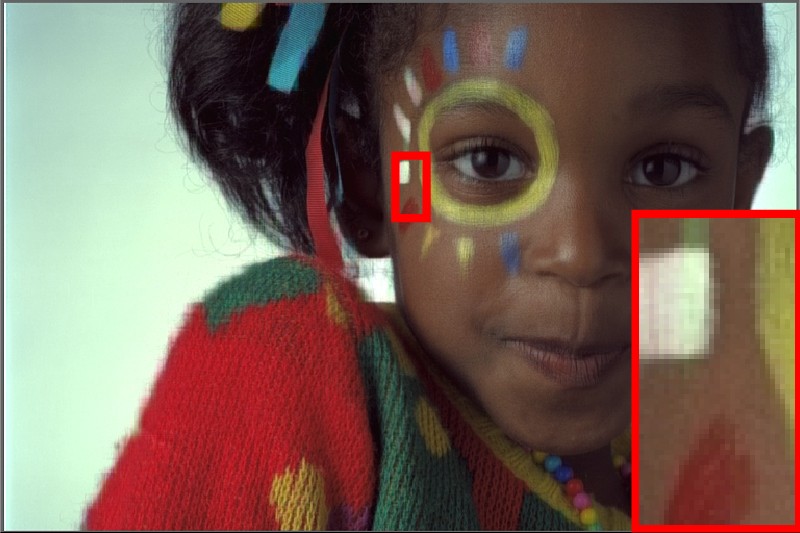}  &
\includegraphics[width=0.115\textwidth]{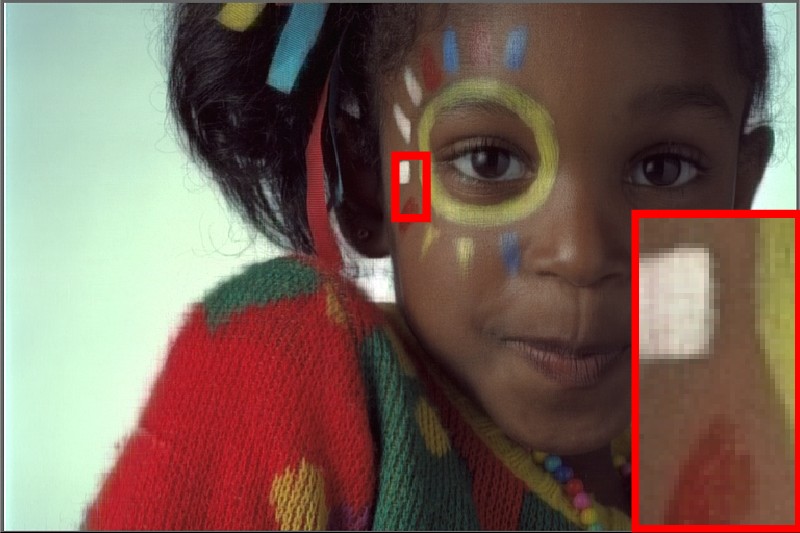} & 
\includegraphics[width=0.115\textwidth]{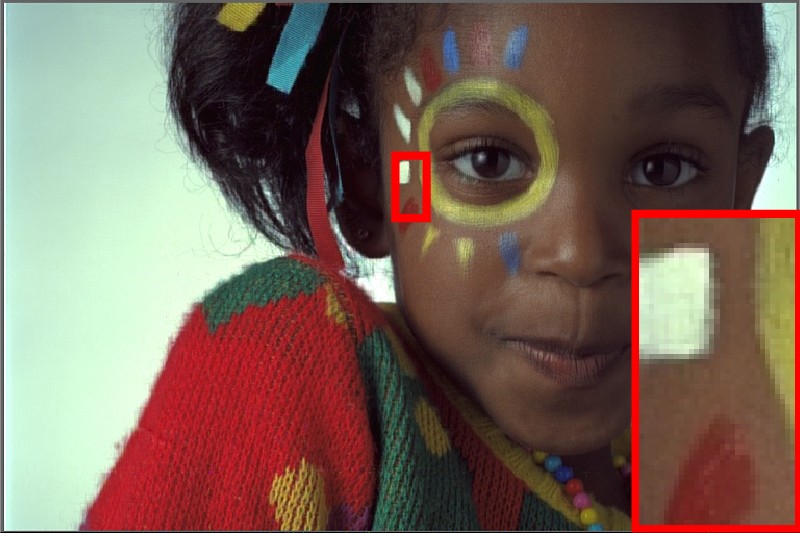}
\\
  \includegraphics[width=0.115\textwidth]{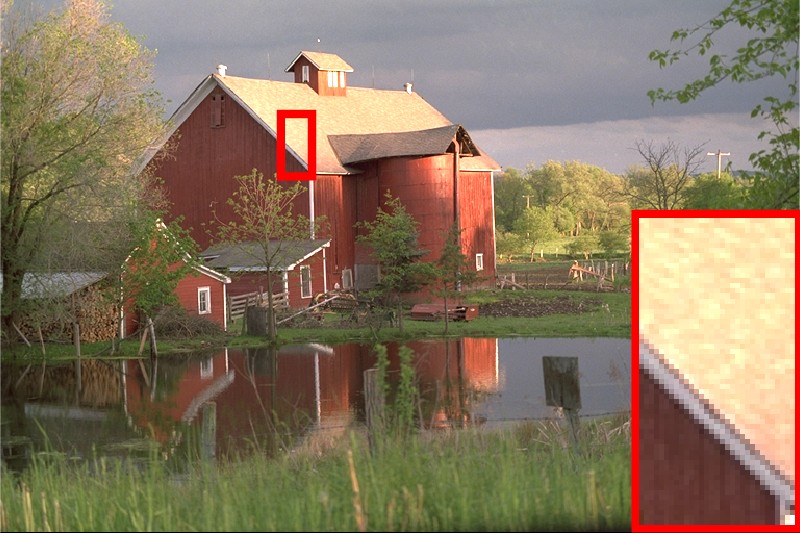}  &
\includegraphics[width=0.115\textwidth]{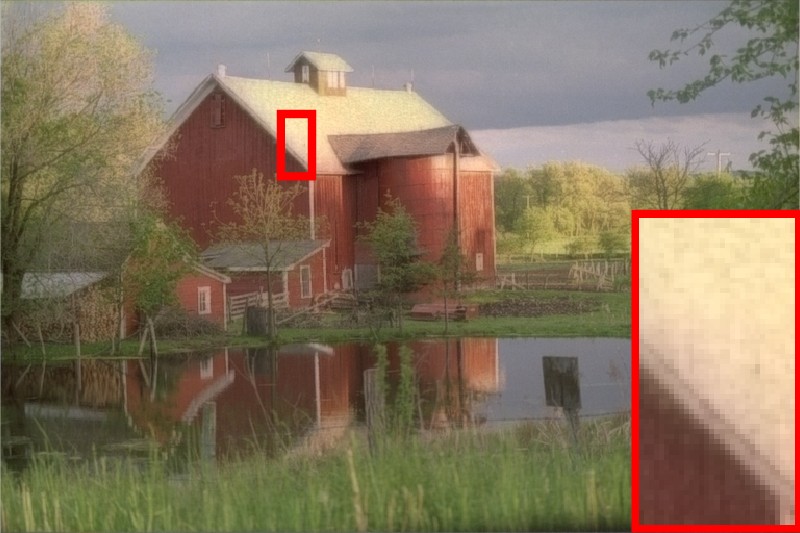}  &
\includegraphics[width=0.115\textwidth]{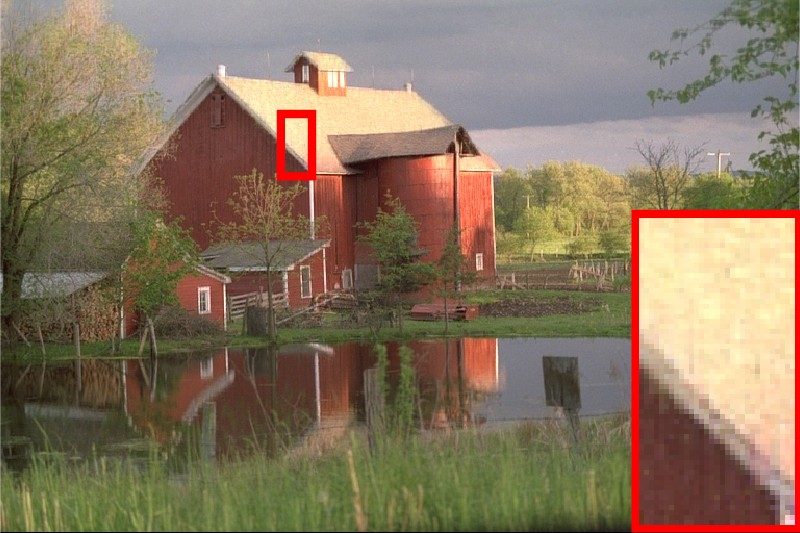} &
\includegraphics[width=0.115\textwidth]{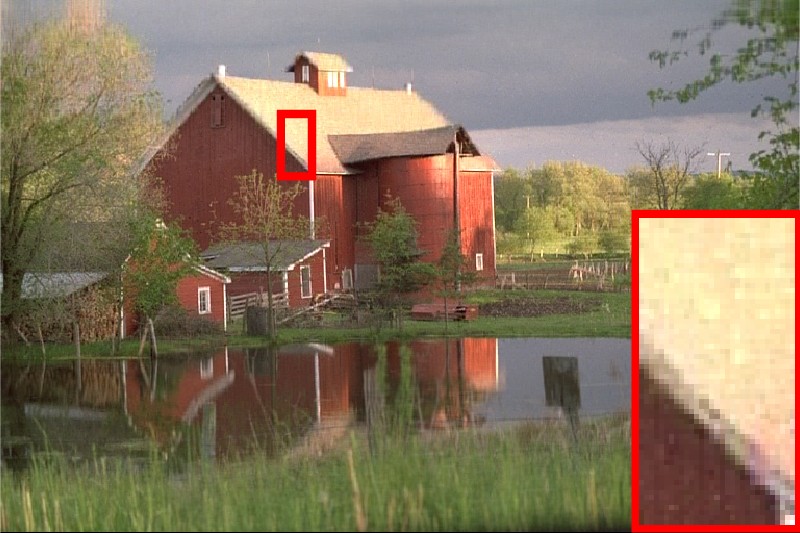} & 
\includegraphics[width=0.115\textwidth]{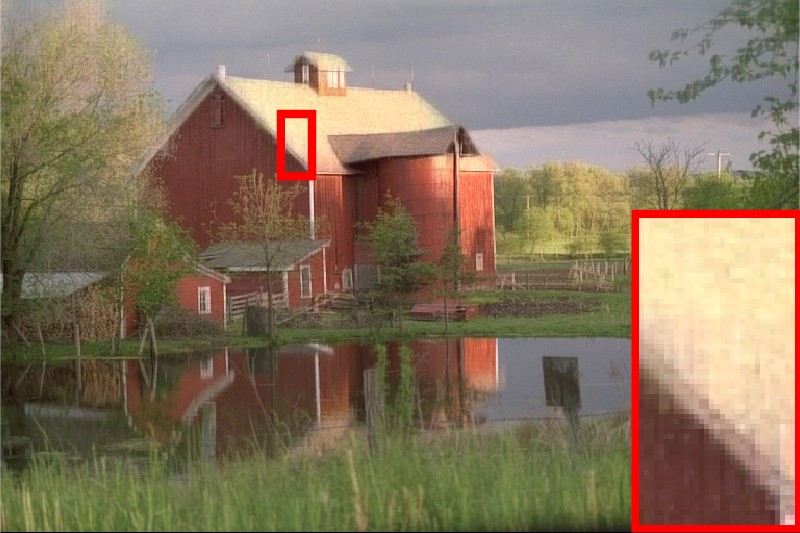}  &
\includegraphics[width=0.115\textwidth]{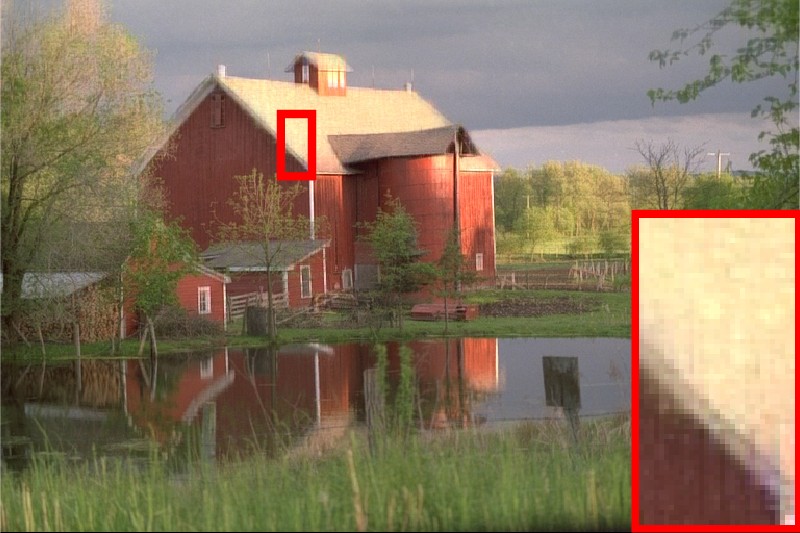}  &
\includegraphics[width=0.115\textwidth]{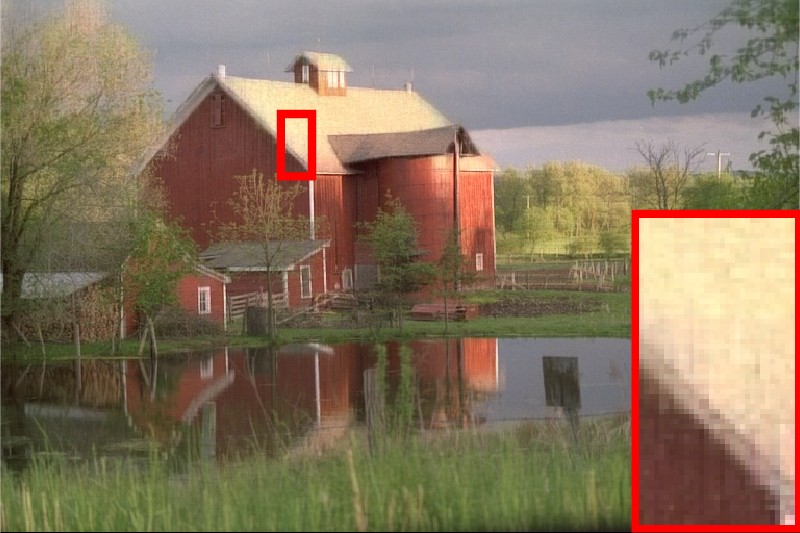} & 
\includegraphics[width=0.115\textwidth]{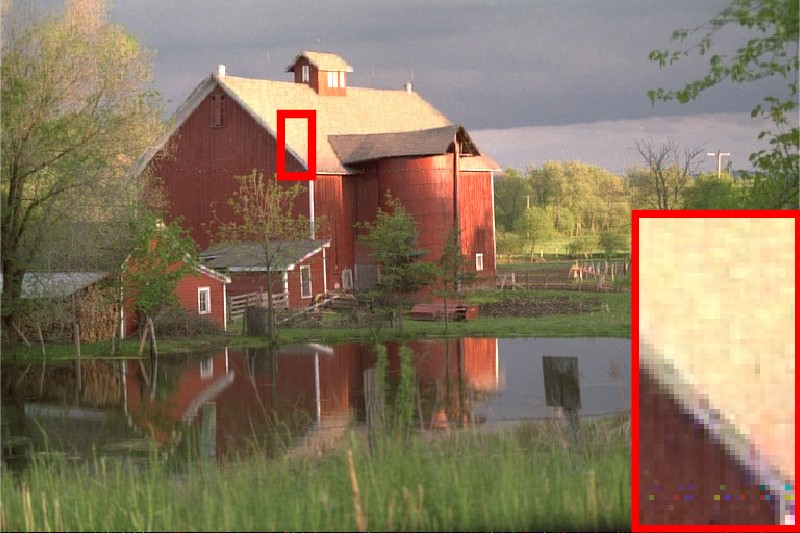}
  \\
 \includegraphics[width=0.115\textwidth]{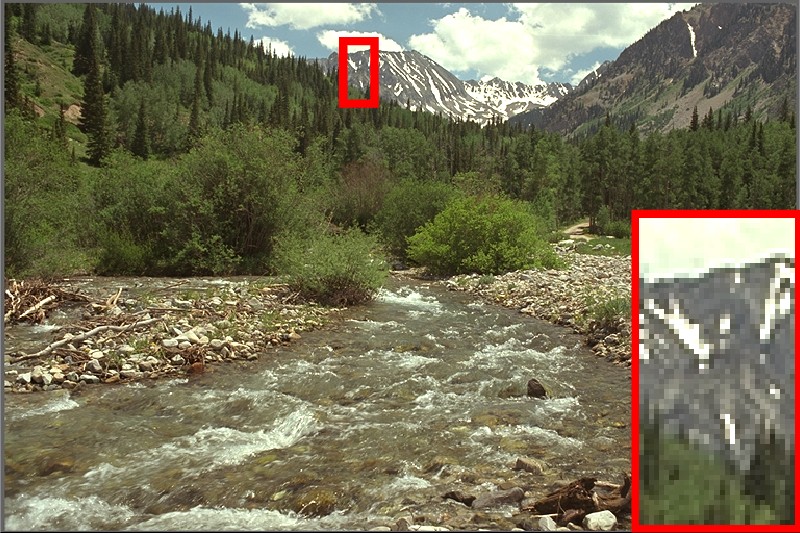}  &
\includegraphics[width=0.115\textwidth]{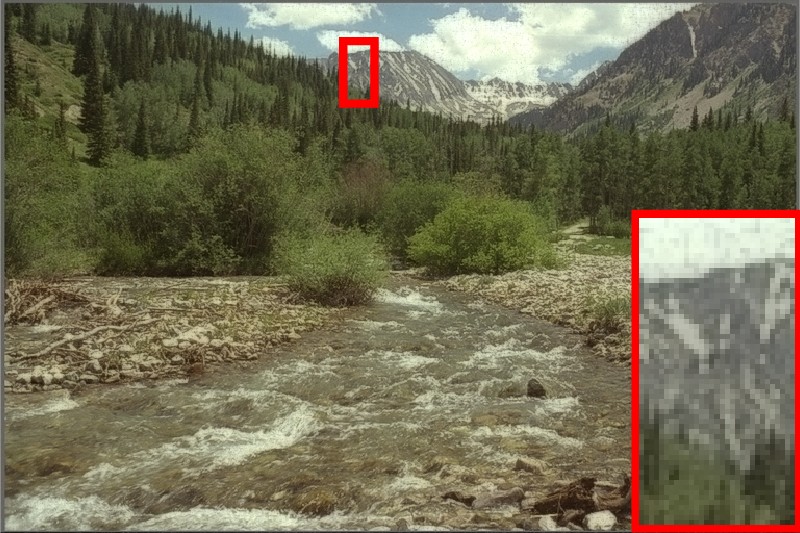}  &
\includegraphics[width=0.115\textwidth]{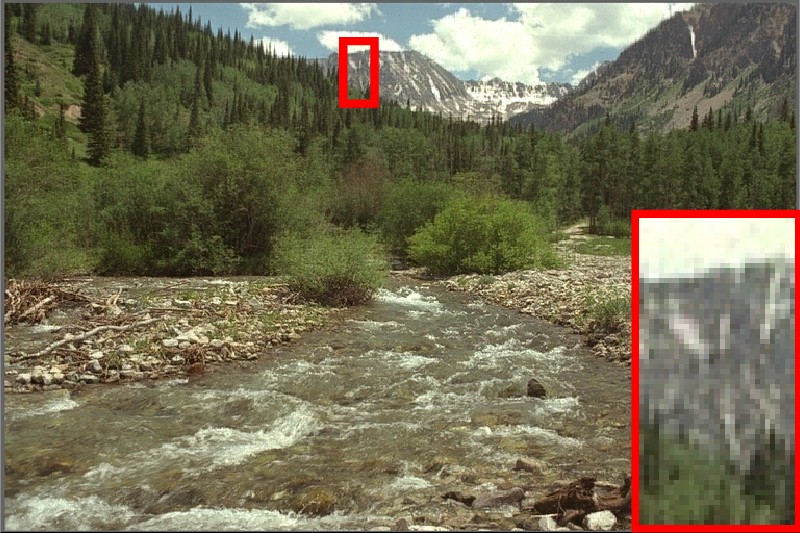}  &
\includegraphics[width=0.115\textwidth]{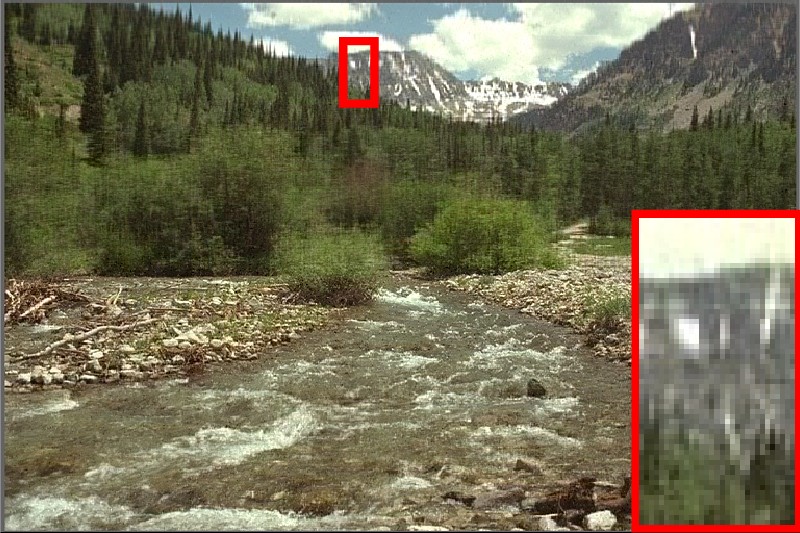}  & 
\includegraphics[width=0.115\textwidth]{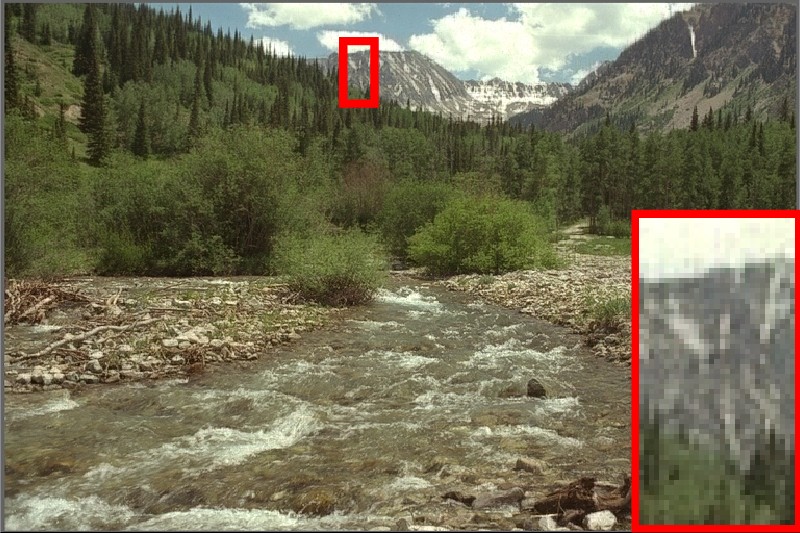}   &
\includegraphics[width=0.115\textwidth]{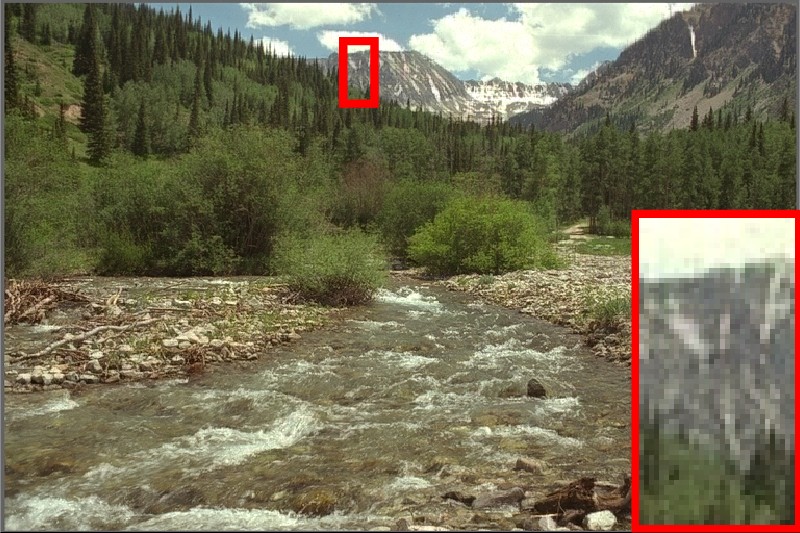}   &
\includegraphics[width=0.115\textwidth]{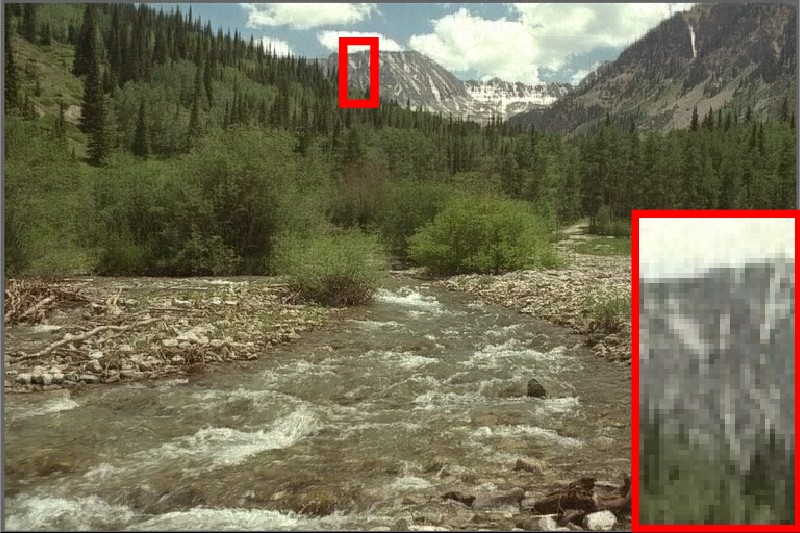}  & 
\includegraphics[width=0.115\textwidth]{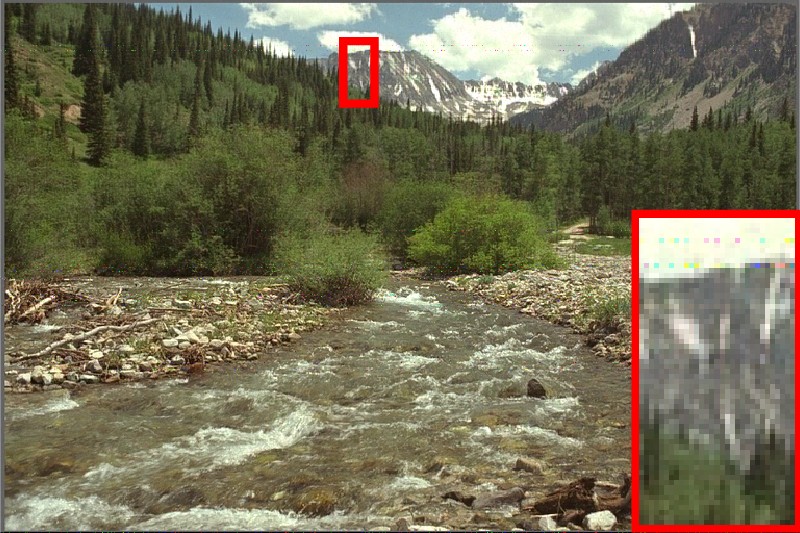} \\
   \includegraphics[width=0.115\textwidth]{Figure/kodim03_ori.jpg}  &
\includegraphics[width=0.115\textwidth]{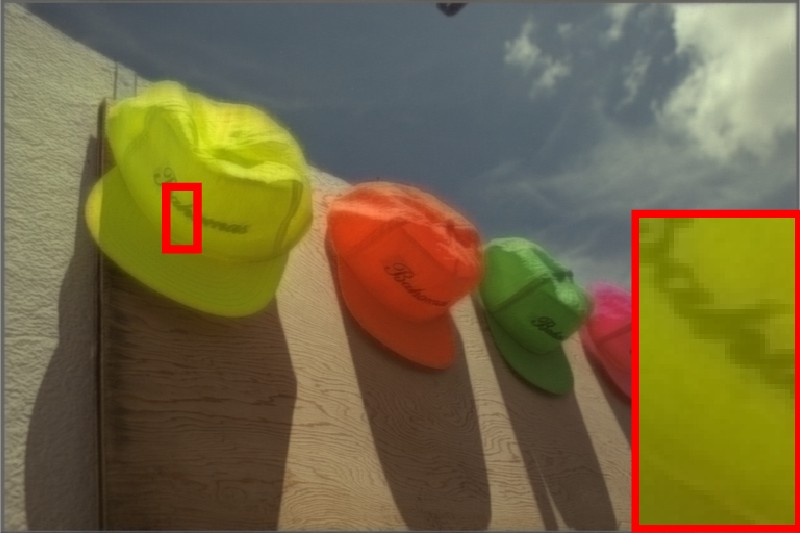}   &
\includegraphics[width=0.115\textwidth]{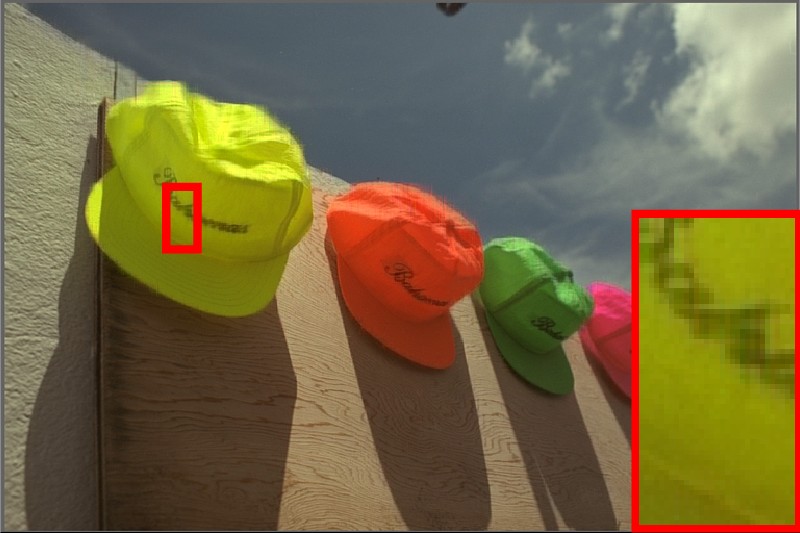}  &
\includegraphics[width=0.115\textwidth]{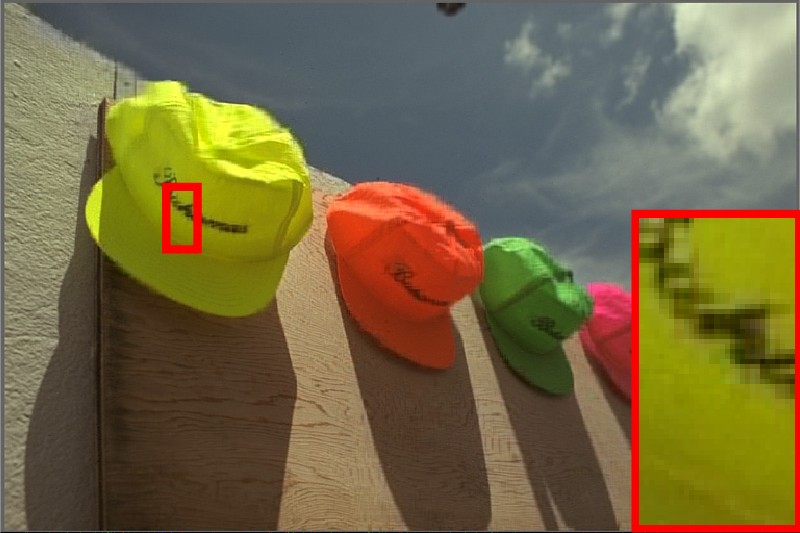}  & 
\includegraphics[width=0.115\textwidth]{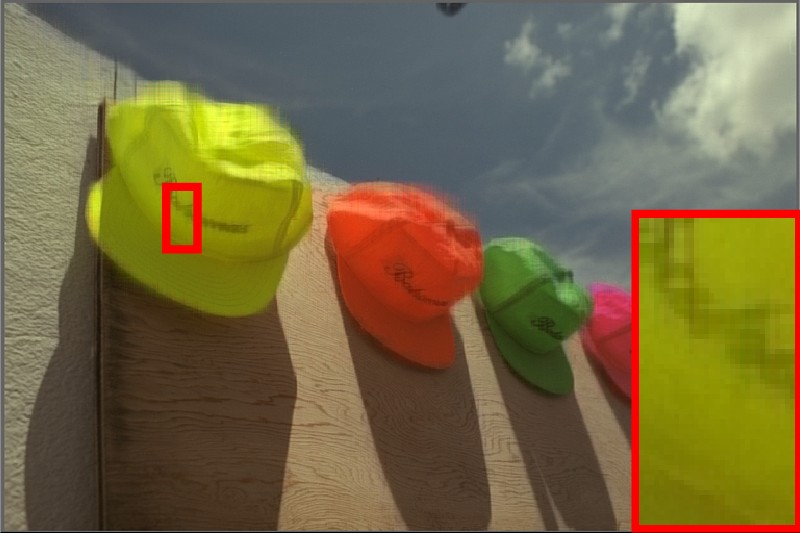}   &
\includegraphics[width=0.115\textwidth]{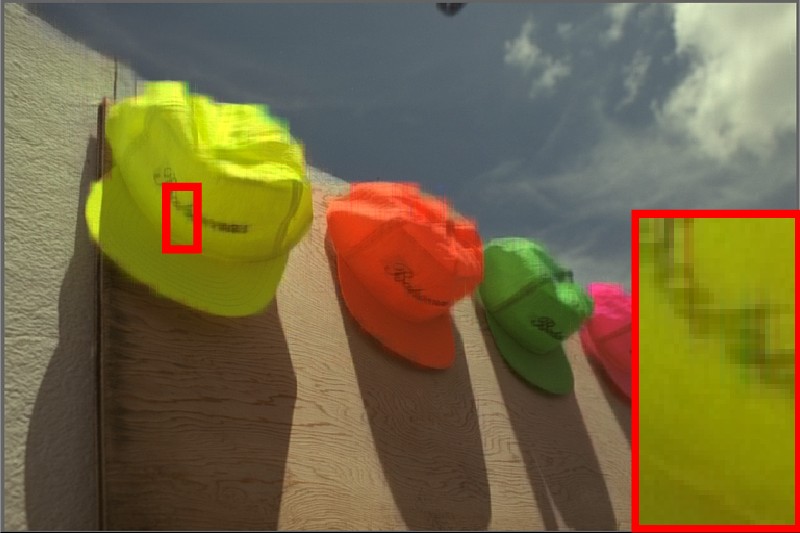}   &
\includegraphics[width=0.115\textwidth]{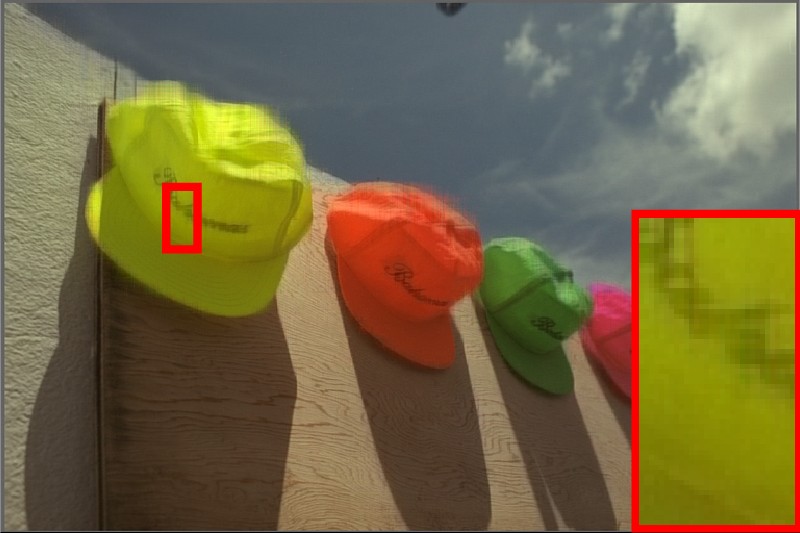} & 
\includegraphics[width=0.115\textwidth]{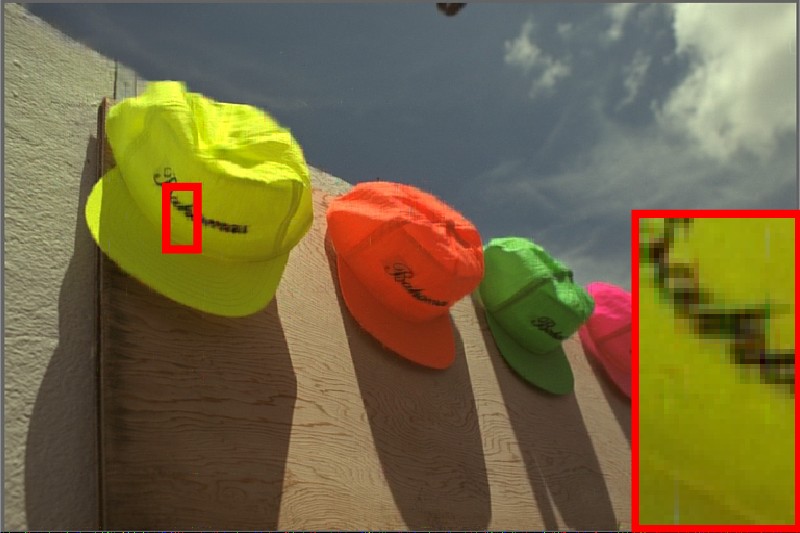}
\end{tabular}

\caption{Comparison of color image mixed noise removal performance on four examples.
}\label{fig:mixed}
\end{figure*}

We benchmark our proposed algorithm against several state-of-the-art methods, including 3DTNN \cite{zheng2019mixed}, $S_{wp}(0.9)$ \cite{YANG2022108311}, BTRTF \cite{zhou2019bayesian}, TNN \cite{lu2019tensor}, and PSTNN \cite{jiang2020multi}. Performance metrics such as PSNR and SSIM are detailed in \eqref{tab:mixed}, with visual results presented in \eqref{fig:mixed}. Notably, our algorithm outperforms both TNN and PSTNN—methods that utilize similar regularization techniques—across all test cases in terms of PSNR, achieving an average improvement of 0.6 dB over the best-reported results. Qualitatively, the images restored by VBI$_{\rm PSTNN}$ exhibit notably sharper boundaries compared to those produced by the other methods, which tend to exhibit some degree of blurring.

\subsection{Background modeling}
The background modeling problem focuses on distinguishing foreground objects from the background in video sequences. This is commonly achieved by modeling the background as a low-rank tensor, which represents the relatively static scenes across different frames, and treating the moving foreground objects as sparse components. In the context of Tensor Robust Principal Component Analysis (TRPCA), these are represented by the low-rank tensor $\mathcal{L}_0$ and the sparse tensor $\mathcal{S}_0$, respectively.

\begin{figure*}[h]
\centering
\scriptsize
\setlength{\tabcolsep}{0.2em}
\begin{tabular}{cccccccccc}
{Orignal} & { 3DTNN}  & { $S_{wp}$(0.9)} & { BTRTF }  & { TNN}  & { PSTNN} & { VBI$_{\rm TNN} $ }  & { VBI$_{\rm PSTNN}$ }\\
\includegraphics[width=0.115\textwidth]{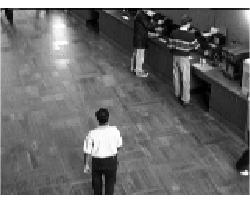}   &
\includegraphics[width=0.115\textwidth]{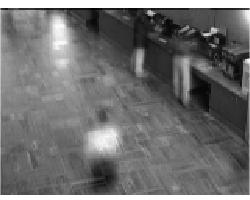}   &
\includegraphics[width=0.115\textwidth]{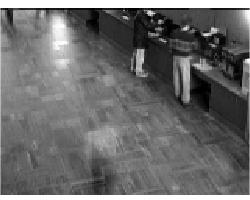}  &
\includegraphics[width=0.115\textwidth]{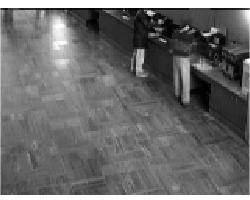}  & 
\includegraphics[width=0.115\textwidth]{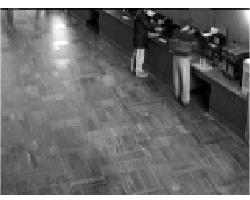}   &
\includegraphics[width=0.115\textwidth]{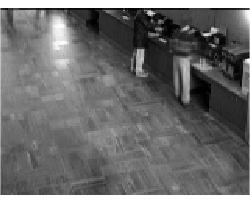}   &
\includegraphics[width=0.115\textwidth]{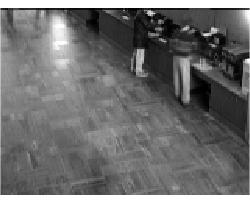}  & 
\includegraphics[width=0.115\textwidth]{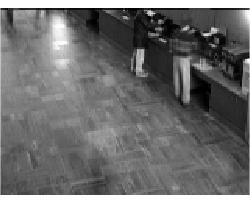}
\\
  
   &
\includegraphics[width=0.115\textwidth]{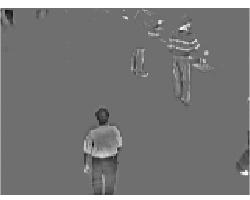}   &
\includegraphics[width=0.115\textwidth]{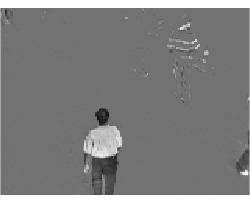}  &
\includegraphics[width=0.115\textwidth]{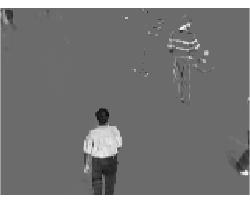}  & 
\includegraphics[width=0.115\textwidth]{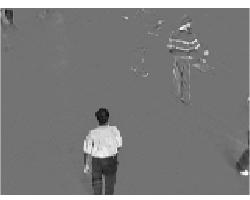}   &
\includegraphics[width=0.115\textwidth]{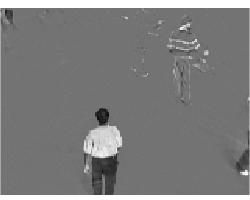}   &
\includegraphics[width=0.115\textwidth]{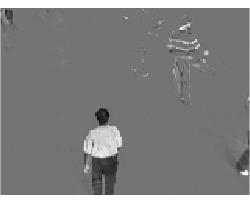}  & 
\includegraphics[width=0.115\textwidth]{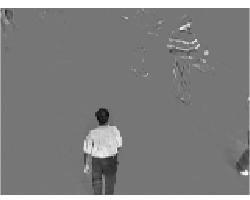} \\

\includegraphics[width=0.115\textwidth]{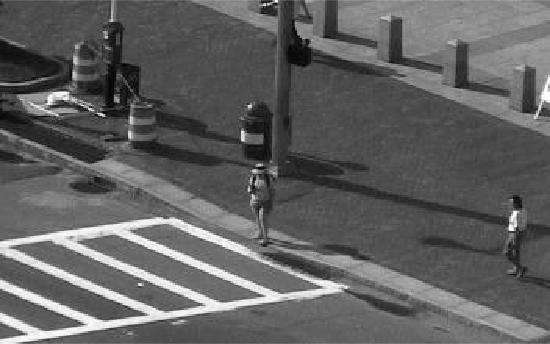}   &
\includegraphics[width=0.115\textwidth]{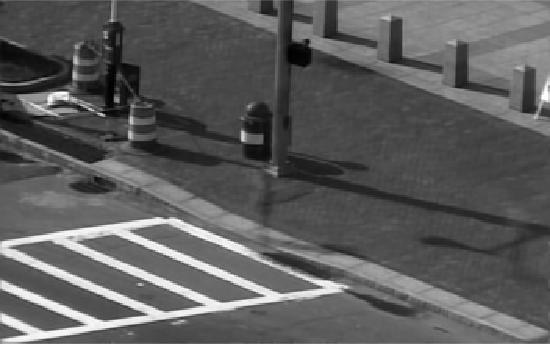}   &
\includegraphics[width=0.115\textwidth]{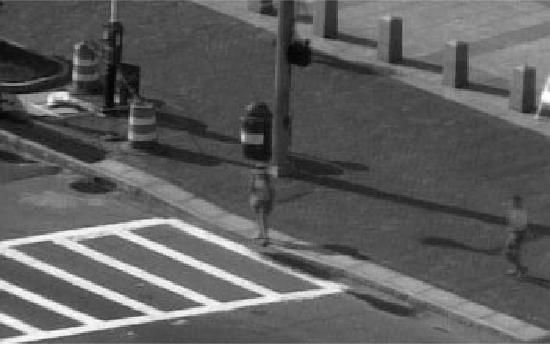}  &
\includegraphics[width=0.115\textwidth]{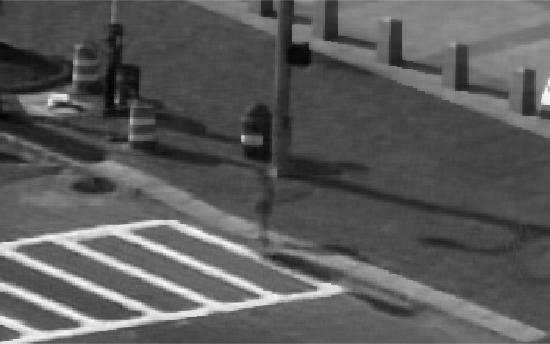}  & 
\includegraphics[width=0.115\textwidth]{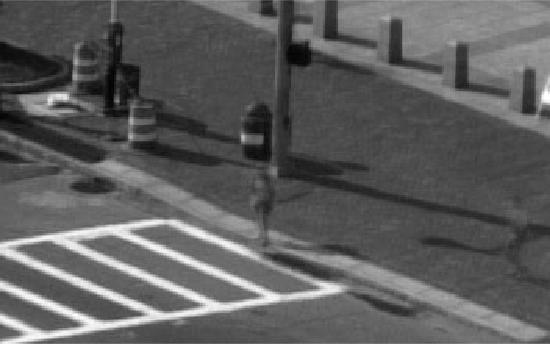}   &
\includegraphics[width=0.115\textwidth]{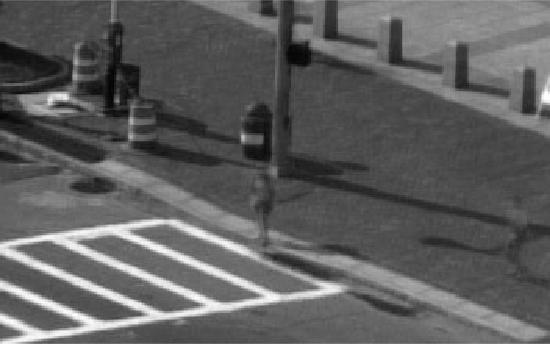}   &
\includegraphics[width=0.115\textwidth]{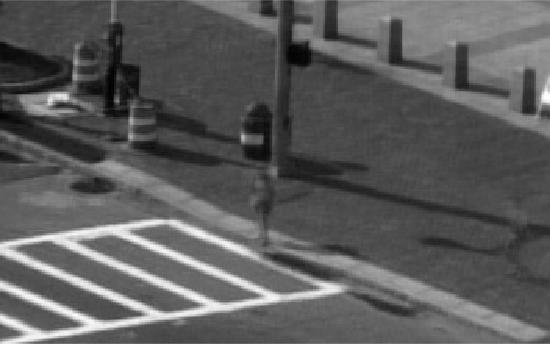}  & 
\includegraphics[width=0.115\textwidth]{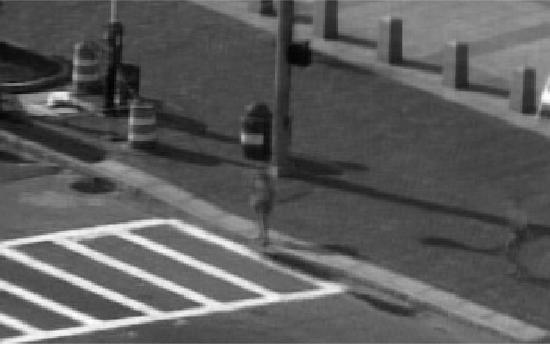}
\\
  
   &
\includegraphics[width=0.115\textwidth]{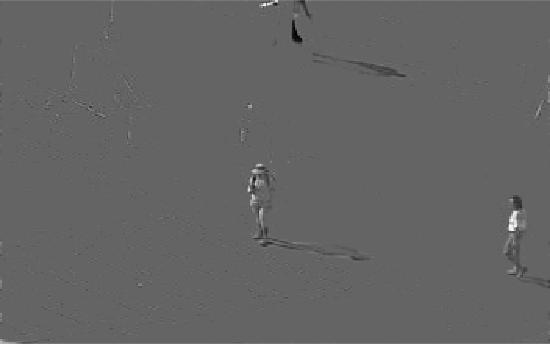}   &
\includegraphics[width=0.115\textwidth]{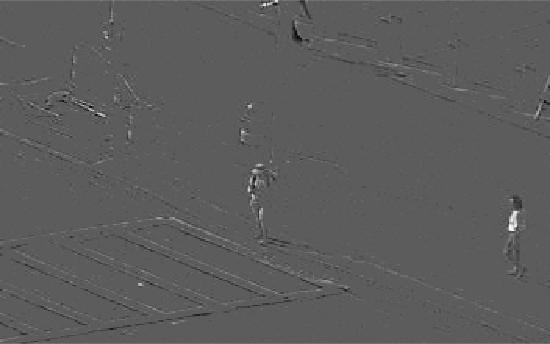}  &
\includegraphics[width=0.115\textwidth]{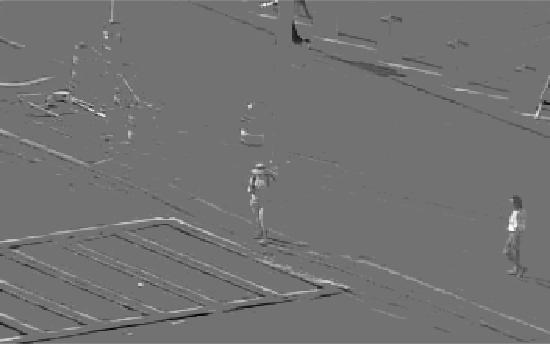}  & 
\includegraphics[width=0.115\textwidth]{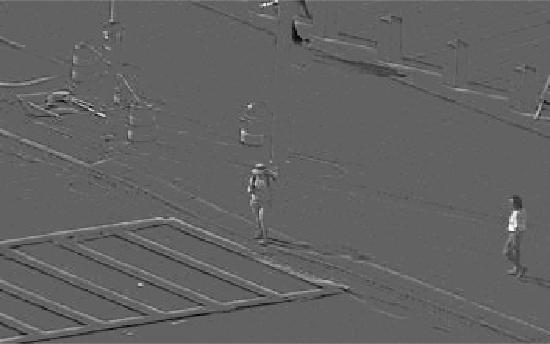}   &
\includegraphics[width=0.115\textwidth]{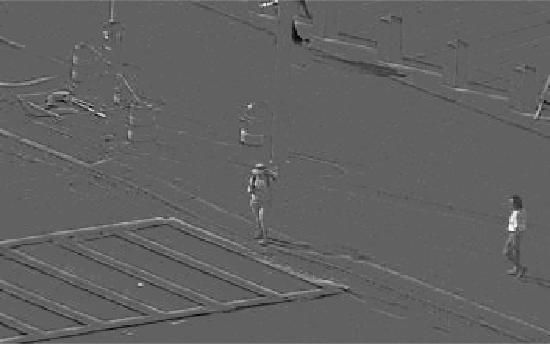}   &
\includegraphics[width=0.115\textwidth]{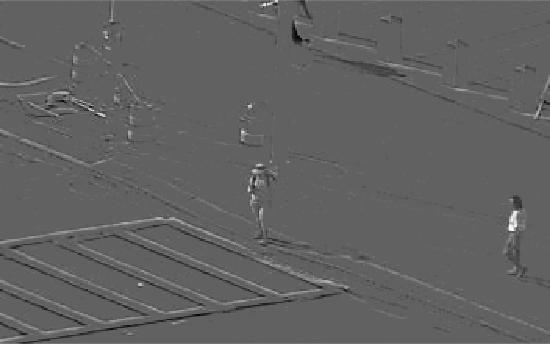}  & 
\includegraphics[width=0.115\textwidth]{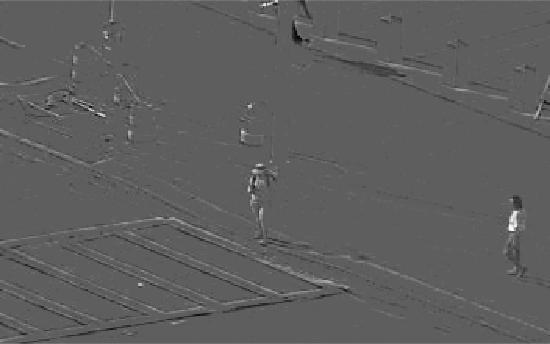}

\end{tabular}

\caption{Background modeling results of two surveillance video sequences.
}\label{fig:bg}
\end{figure*}

We evaluated our models on sequences from the 12R dataset \cite{li2004statistical}, specifically the  ``bootstrap" ($120\times160\times400$), and ``sidewalk'' ($220\times 352 \times 400 $) videos, all characterized by slow-moving objects against varying backgrounds. Our models were compared with several others, including 3DTNN, TNN, BTRTF, PSTNN,  and $\mathrm{t}\mbox{-}S_{w, p}(0.9)$. For VBI$_{\rm PSTNN}$, the truncated parameter $K$ is set as 5, while the initial values of $\theta_1, \theta_2, \theta_3$ are set as 1, 1, 100, respectively.
The results of these comparisons are visually presented in Figure~\ref{fig:bg}. Each video's analysis starts with a frame from the sequence as shown in column (a) of Figure~\ref{fig:bg}, followed by background images generated by the respective methods, from 3DTNN to our approach VBI$_{\rm PSTNN}$. Additionally, the motion in each scene is depicted in the second row for each video.
In the ``bootstrap" video, except for 3DTNN, all the methods achieved superior background separation with fewer ghost silhouettes. In the ``sidewalk" videos, all the approaches perform similarly, while 3DTNN has slightly better results.

\section{Conclusions}\label{sec:concl}
In this paper, we presented a method for recovering low-rank tensors from observations contaminated by sparse outliers and Gaussian noise. Utilizing variational Bayesian inference, we effectively resolved the tensors while simultaneously selecting model parameters. Numerical evaluations highlight the advantages and superior performance of our approach compared to existing methods. Currently limited to linear and convex relaxations, our future work will explore extending this parameter selection technique to nonconvex approximations within tensor recovery models.


\bibliographystyle{siamplain}
\bibliography{references,new_ref,library}

@article{bhattacharya2025,
  title={On the convergence of coordinate ascent variational inference},
  author={Bhattacharya, Anirban and Pati, Debdeep and Yang, Yun},
  journal={The Annals of Statistics},
  volume={53},
  number={3},
  pages={929--962},
  year={2025},
  publisher={Institute of Mathematical Statistics}
}

@article{golub1979,
	author = {Golub, G. and Heath, M. and Wahba, G.},
	journal = {Technometrics},
	number = {2},
	pages = {215--223},
	title = {Generalized cross-validation as a method for choosing a good ridge parameter},
	volume = {21},
	year = {1979}}

@book{Morozov,
	author = {Morozov, V.},
	publisher = {Springer-Verlag, New York},
	title = {Methods for solving incorrectly posed problems},
	year = {1984}}

@article{2013parameter,
	author = {Almeida, M. and Figueiredo, M.},
	journal = {IEEE Trans. Image Process.},
	number = {7},
	pages = {2751-2763},
	title = {Parameter estimation for blind and non-blind deblurring using residual whiteness measures},
	volume = {22},
	year = {2013}}

@article{2023Poisson-whiteness,
	author = {Bevilacqua, F. and Lanza, A. and Pragliola, M. and Sgallari, F.},
	journal = {Applied Mathematical Modelling},
	pages = {197--218},
	title = {Whiteness-based parameter selection for poisson data in variational image processing},
	volume = {117},
	year = {2023}}

@article{Hansen1993,
	author = {P. Hansen and D. O'Leary},
	journal = {SIAM J. Sci. Comput},
	number = {6},
	pages = {1487-1503},
	title = {The use of the {$L$}-curve in the regularization of discrete ill-posed problems},
	volume = {14},
	year = {1993}}

@article{caijf2010,
	author = {Cai, J and Candes, E and Shen, Z},
	journal = {SIAM Journal on Optimization},
	number = {4},
	pages = {1956--1982},
	title = {{A singular value thresholding algorithm for matrix completion}},
	volume = {20},
	year = {2010}}

@article{tibshirani1996regression,
	author = {Tibshirani, R},
	journal = {Journal of the Royal Statistical Society. Series B (Methodological)},
	pages = {267--288},
	publisher = {JSTOR},
	title = {{Regression shrinkage and selection via the lasso}},
	year = {1996}}

@article{cai2024a,
	title = {Robust Tensor CUR Decompositions: Rapid Low-{Tucker}-Rank Tensor Recovery with Sparse Corruptions},
	volume = {17},
	issn = {1936-4954},
	shorttitle = {Robust {Tensor} {CUR} {Decompositions}},
	language = {en},
	number = {1},
	urldate = {2024-06-14},
	journal = {SIAM J. Imag. Sci.},
	author = {Cai, HanQin and Chao, Zehan and Huang, Longxiu and Needell, Deanna},
	month = mar,
	year = {2024},
	pages = {225--247},
}

@article{zhou2024c,
	title = {Proximal {MCMC} for {Bayesian} Inference of Constrained and Regularized Estimation},
	issn = {0003-1305, 1537-2731},
	language = {en},
	urldate = {2024-08-10},
	journal = {Am. Stat.},
	author = {Zhou, Xinkai and Heng, Qiang and Chi, Eric C. and Zhou, Hua},
	month = feb,
	year = {2024},
	pages = {1--12},
}

@article{helin2022,
	title = {Edge-Promoting Adaptive {Bayesian} Experimental Design for {X}-ray {Imaging}},
	volume = {44},
	language = {en},
	number = {3},
	urldate = {2024-08-12},
	journal = {SIAM J. Sci. Comput.},
	author = {Helin, Tapio and Hyvönen, Nuutti and Puska, Juha-Pekka},
	month = jun,
	year = {2022},
	pages = {B506--B530},
}

@article{menzen2022,
	title = {Alternating Linear Scheme in a {Bayesian} Framework for Low-Rank Tensor Approximation},
	volume = {44},
	issn = {1064-8275, 1095-7197},
	language = {en},
	number = {3},
	urldate = {2024-08-12},
	journal = {SIAM J. Sci. Comput.},
	author = {Menzen, Clara and Kok, Manon and Batselier, Kim},
	month = jun,
	year = {2022},
	pages = {A1116--A1144},
}

@article{bazerque2013,
	title = {Rank Regularization and {Bayesian} Inference for Tensor Completion and Extrapolation},
	volume = {61},
	copyright = {https://ieeexplore.ieee.org/Xplorehelp/downloads/license-information/IEEE.html},
	issn = {1053-587X, 1941-0476},
	url = {http://ieeexplore.ieee.org/document/6579771/},
	doi = {10.1109/TSP.2013.2278516},
	language = {en-US},
	number = {22},
	urldate = {2024-09-01},
	journal = {IEEE Trans. Signal Process.},
	author = {Bazerque, Juan Andres and Mateos, Gonzalo and Giannakis, Georgios B.},
	month = nov,
	year = {2013},
	pages = {5689--5703},
}

@article{yang2018a,
	title = {Fast Low-Rank {Bayesian} Matrix Completion With Hierarchical {Gaussian} Prior Models},
	volume = {66},
	issn = {1053-587X, 1941-0476},
	language = {en-US},
	number = {11},
	urldate = {2024-09-01},
	journal = {IEEE Trans. Signal Process.},
	author = {Yang, Linxiao and Fang, Jun and Duan, Huiping and Li, Hongbin and Zeng, Bing},
	month = jun,
	year = {2018},
	pages = {2804--2817},
}

@article{tong2023,
	title = {Bayesian Tensor {Tucker} Completion With a Flexible Core},
	volume = {71},
	issn = {1053-587X, 1941-0476},
	urldate = {2024-09-01},
	journal = {IEEE Trans. Signal Process.},
	author = {Tong, Xueke and Cheng, Lei and Wu, Yik-Chung},
	year = {2023},
	pages = {4077--4091},
}

@article{uribe2023,
	title = {Horseshoe Priors for Edge-Preserving Linear {Bayesian} Inversion},
	volume = {45},
	issn = {1064-8275, 1095-7197},
	language = {en},
	number = {3},
	urldate = {2024-08-12},
	journal = {SIAM J. Sci. Comput.},
	author = {Uribe, Felipe and Dong, Yiqiu and Hansen, Per Christian},
	month = jun,
	year = {2023},
	pages = {B337--B365},
}

@article{jia2021,
	title = {Variational {Bayes}' method for functions with applications to some inverse problems},
	volume = {43},
	issn = {1064-8275, 1095-7197},
	language = {en-US},
	number = {1},
	urldate = {2024-08-12},
	journal = {SIAM J. Sci. Comput.},
	author = {Jia, Junxiong and Zhao, Qian and Xu, Zongben and Meng, Deyu and Leung, Yee},
	month = jan,
	year = {2021},
	pages = {A355--A383},
}

@article{glaubitz2024,
	title = {Leveraging joint sparsity in hierarchical {Bayesian} learning},
	volume = {12},
	issn = {2166-2525},
	number = {2},
	urldate = {2024-08-12},
	journal = {SIAM/ASA J. Uncertainty Quantif.	},
	author = {Glaubitz, Jan and Gelb, Anne},
	month = jun,
	year = {2024},
	pages = {442--472},
}

@article{law2022,
	title = {Sparse online variational {Bayesian} regression},
	volume = {10},
	issn = {2166-2525},
	number = {3},
	urldate = {2024-08-12},
	journal = {SIAM/ASA J. Uncertainty Quantif.},
	author = {Law, Kody J. H. and Zankin, Vitaly},
	month = sep,
	year = {2022},
	pages = {1070--1100},
}

@article{chantas2009variational,
  title={Variational Bayesian image restoration with a product of spatially weighted total variation image priors},
  author={Chantas, Giannis and Galatsanos, Nikolaos P and Molina, Rafael and Katsaggelos, Aggelos K},
  journal={IEEE Trans. Image Process.},
  volume={19},
  number={2},
  pages={351--362},
  year={2009},
  publisher={IEEE}
}

@article{ABMK,
  title={Sparse {Bayesian} blind image deconvolution with parameter estimation},
  author={Amizic, Bruno and Molina, Rafael and Katsaggelos, Aggelos K},
  journal={EUSIPCO},
  volume={2012},
  number={1},
  pages={1--15},
  year={2012},
  publisher={SpringerOpen}
}

@article{Gamma1,
  title={Parameter estimation in {TV} image restoration using variational distribution approximation},
  author={Babacan, S Derin and Molina, Rafael and Katsaggelos, Aggelos K},
  journal={IEEE Trans. Image Process.},
  volume={17},
  number={3},
  pages={326--339},
  year={2008},
  publisher={IEEE}
}

@article{MCMC,
  title={{MCMC}-based image reconstruction with uncertainty quantification},
  author={Bardsley, Johnathan M},
  journal={SIAM J. Sci. Comput.},
  volume={34},
  number={3},
  pages={A1316--A1332},
  year={2012},
  publisher={SIAM}
}

@article{OBF,
  title={Adaptive total variation image deblurring: a majorization--minimization approach},
  author={Oliveira, Joao P and Bioucas-Dias, Jos{\'e} M and Figueiredo, M{\'a}rio AT},
  journal={Signal Process.},
  volume={89},
  number={9},
  pages={1683--1693},
  year={2009},
  publisher={Elsevier}
}

@article{Gamma-non,
  title={Sparse Bayesian learning and the relevance vector machine},
  author={Tipping, Michael E},
  journal={J. Mach. Learn. Res.},
  volume={1},
  number={Jun},
  pages={211--244},
  year={2001}
}

@article{weakly,
  title={A weakly informative default prior distribution for logistic and other regression models},
  author={Gelman, Andrew and Jakulin, Aleks and Pittau, Maria Grazia and Su, Yu-Sung},
  journal={Ann Appl. Stat.},
  volume={2},
  number={4},
  pages={1360--1383},
  year={2008},
  publisher={Institute of Mathematical Statistics}
}

@article{wang2013variational,
  title={Variational inference in nonconjugate models},
  author={Wang, Chong and Blei, David M},
  journal={J. Mach. Learn. Res.},
  year={2013}
}

@article{lalush1992simulation,
  title={Simulation evaluation of {Gibbs} prior distributions for use in maximum a posteriori {SPECT} reconstructions},
  author={Lalush, David S and Tsui, Benjamin MW},
  journal={IEEE Trans. Med. Imaging},
  volume={11},
  number={2},
  pages={267--275},
  year={1992},
  publisher={IEEE}
}

@inproceedings{lefkimmiatis2023,
  title={Learning Sparse and Low-Rank Priors for Image Recovery via Iterative Reweighted Least Squares Minimization},
  author={Stamatios Lefkimmiatis and Iaroslav Koshelev},
  booktitle={Proc. ICLR},
  year={2023}
}

@article{wang2004image,
  title={Image quality assessment: from error visibility to structural similarity},
  author={Wang, Zhou and Bovik, Alan C and Sheikh, Hamid R and Simoncelli, Eero P},
  journal={IEEE Trans. Image Process.},
  volume={13},
  number={4},
  pages={600--612},
  year={2004},
  publisher={IEEE}
}

@article{kolda2009tensor,
  title={Tensor decompositions and applications},
  author={Kolda, Tamara G and Bader, Brett W},
  journal={SIAM Rev.},
  volume={51},
  number={3},
  pages={455--500},
  year={2009},
  publisher={SIAM}
}

@article{kilmer2011factorization,
  title={Factorization strategies for third-order tensors},
  author={Kilmer, Misha E and Martin, Carla D},
  journal={Linear Algebra Appl.},
  volume={435},
  number={3},
  pages={641--658},
  year={2011},
  publisher={Elsevier}
}

@article{he2015total,
  title={Total-variation-regularized low-rank matrix factorization for hyperspectral image restoration},
  author={He, Wei and Zhang, Hongyan and Zhang, Liangpei and Shen, Huanfeng},
  journal={IEEE Trans. Geosci. Remote Sens.},
  volume={54},
  number={1},
  pages={178--188},
  year={2015},
  publisher={IEEE}
}

@article{jolliffe2016principal,
  title={Principal component analysis: a review and recent developments},
  author={Jolliffe, Ian T and Cadima, Jorge},
  journal={Phil. Trans. R. Soc. A.	},
  volume={374},
  number={2065},
  pages={20150202},
  year={2016},
  publisher={The Royal Society Publishing}
}

@article{blei2017variational,
  title={Variational inference: A review for statisticians},
  author={Blei, David M and Kucukelbir, Alp and McAuliffe, Jon D},
  journal={Journal of the American statistical Association},
  volume={112},
  number={518},
  pages={859--877},
  year={2017},
  publisher={Taylor \& Francis}
}

@ARTICLE{lu2019tensor,
  author={Lu, Canyi and Feng, Jiashi and Chen, Yudong and Liu, Wei and Lin, Zhouchen and Yan, Shuicheng},
  journal={IEEE Trans. Pattern Anal. Mach. Intell.}, 
  title={Tensor Robust Principal Component Analysis with a New Tensor Nuclear Norm}, 
  year={2020},
  volume={42},
  number={4},
  pages={925-938}}

@article{zhou2019bayesian,
  title={Bayesian low-tubal-rank robust tensor factorization with multi-rank determination},
  author={Zhou, Yang and Cheung, Yiu-Ming},
  journal={IEEE Trans. Pattern Anal. Mach. Intell.},
  volume={43},
  number={1},
  pages={62--76},
  year={2019},
  publisher={IEEE}
}

@article{zheng2019mixed,
  title={Mixed noise removal in hyperspectral image via low-fibered-rank regularization},
  author={Zheng, Yu-Bang and Huang, Ting-Zhu and Zhao, Xi-Le and Jiang, Tai-Xiang and Ma, Tian-Hui and Ji, Teng-Yu},
  journal={IEEE Trans. Geosci. Remote Sens.},
  volume={58},
  number={1},
  pages={734--749},
  year={2019},
  publisher={IEEE}
}

@article{jiang2020multi,
  title={Multi-dimensional imaging data recovery via minimizing the partial sum of tubal nuclear norm},
  author={Jiang, Tai-Xiang and Huang, Ting-Zhu and Zhao, Xi-Le and Deng, Liang-Jian},
  journal={J. Comput. Appl. Math.},
  volume={372},
  pages={112680},
  year={2020},
  publisher={Elsevier}
}

@inproceedings{phan2020stable,
  title={Stable low-rank tensor decomposition for compression of convolutional neural network},
  author={Phan, Anh-Huy and Sobolev, Konstantin and Sozykin, Konstantin and Ermilov, Dmitry and Gusak, Julia and Tichavsk{\`y}, Petr and Glukhov, Valeriy and Oseledets, Ivan and Cichocki, Andrzej},
  booktitle={Proc. ECCV},
  pages={522--539},
  year={2020},
}

@article{chang2020weighted,
  title={Weighted low-rank tensor recovery for hyperspectral image restoration},
  author={Chang, Yi and Yan, Luxin and Zhao, Xi-Le and Fang, Houzhang and Zhang, Zhijun and Zhong, Sheng},
  journal={IEEE Trans. Cybern.},
  volume={50},
  number={11},
  pages={4558--4572},
  year={2020},
  publisher={IEEE}
}

@article{panagakis2021tensor,
  title={Tensor methods in computer vision and deep learning},
  author={Panagakis, Yannis and Kossaifi, Jean and Chrysos, Grigorios G and Oldfield, James and Nicolaou, Mihalis A and Anandkumar, Anima and Zafeiriou, Stefanos},
  journal={Proc. IEEE},
  volume={109},
  number={5},
  pages={863--890},
  year={2021},
  publisher={IEEE}
}

@article{su2022hyperspectral,
  title={Hyperspectral image denoising via weighted multidirectional low-rank tensor recovery},
  author={Su, Yanchi and Zhu, Haoran and Wong, Ka-Chun and Chang, Yi and Li, Xiangtao},
  journal={IEEE Trans. Cybern.},
  volume={53},
  number={5},
  pages={2753--2766},
  year={2022},
  publisher={IEEE}
}

@article{YANG2022108311,
title = {Nonconvex {3D} array image data recovery and pattern recognition under tensor framework},
author = {Ming Yang and Qilun Luo and Wen Li and Mingqing Xiao},
journal = {Pattern Recognit.},
volume = {122},
pages = {108311},
year = {2022},
}

@article{li2004statistical,
  title={Statistical modeling of complex backgrounds for foreground object detection},
  author={Li, Liyuan and Huang, Weimin and Gu, Irene Yu-Hua and Tian, Qi},
  journal={IEEE Trans. Image Process.},
  volume={13},
  number={11},
  pages={1459--1472},
  year={2004},
  publisher={IEEE}
}

@article{candes2012exact,
  title={Exact matrix completion via convex optimization},
  author={Candes, Emmanuel and Recht, Benjamin},
  journal={Commun. ACM},
  volume={55},
  number={6},
  pages={111--119},
  year={2012},
  publisher={ACM New York, NY, USA}
}

@article{kiers2000towards,
  title={Towards a standardized notation and terminology in multiway analysis},
  author={Kiers, Henk AL},
  journal={J. Chemom.},
  volume={14},
  number={3},
  pages={105--122},
  year={2000},
  publisher={Wiley Online Library}
}

@article{gandy2011tensor,
  title={Tensor completion and low-$n$-rank tensor recovery via convex optimization},
  author={Gandy, Silvia and Recht, Benjamin and Yamada, Isao},
  journal={Inverse Probl.},
  volume={27},
  number={2},
  pages={025010},
  year={2011},
  publisher={IOP Publishing}
}

@article{tucker1966some,
  title={Some mathematical notes on three-mode factor analysis},
  author={Tucker, Ledyard R},
  journal={Psychometrika},
  volume={31},
  number={3},
  pages={279--311},
  year={1966},
  publisher={Springer}
}

@article{kilmer2013third,
  title={Third-order tensors as operators on matrices: A theoretical and computational framework with applications in imaging},
  author={Kilmer, Misha E and Braman, Karen and Hao, Ning and Hoover, Randy C},
  journal={SIAM J. Matrix Anal. Appl.},
  volume={34},
  number={1},
  pages={148--172},
  year={2013},
  publisher={SIAM}
}

@inproceedings{mu2014square,
  title={Square deal: Lower bounds and improved relaxations for tensor recovery},
  author={Mu, Cun and Huang, Bo and Wright, John and Goldfarb, Donald},
  booktitle={Proc. ICML},
  pages={73--81},
  year={2014},
  organization={PMLR}
}

@article{goldfarb2014robust,
  title={Robust low-rank tensor recovery: Models and algorithms},
  author={Goldfarb, Donald and Qin, Zhiwei},
  journal={SIAM J. Matrix Anal. Appl.},
  volume={35},
  number={1},
  pages={225--253},
  year={2014},
  publisher={SIAM}
}

@article{wang2023guaranteed,
  title={Guaranteed tensor recovery fused low-rankness and smoothness},
  author={Wang, Hailin and Peng, Jiangjun and Qin, Wenjin and Wang, Jianjun and Meng, Deyu},
  journal={IEEE Trans. Pattern Anal. Mach. Intell.},
  volume={45},
  number={9},
  pages={10990--11007},
  year={2023},
  publisher={IEEE}
}

@article{zheng2024scale,
  title={A Scale-Invariant Relaxation in Low-Rank Tensor Recovery with an Application to Tensor Completion},
  author={Zheng, Huiwen and Lou, Yifei and Tian, Guoliang and Wang, Chao},
  journal={SIAM J. Imag. Sci.},
  volume={17},
  number={1},
  pages={756--783},
  year={2024},
  publisher={SIAM}
}

@article{mu2020weighted,
  title={Weighted tensor nuclear norm minimization for tensor completion using {tensor-SVD}},
  author={Mu, Yang and Wang, Ping and Lu, Liangfu and Zhang, Xuyun and Qi, Lianyong},
  journal={Pattern Recognit. Lett.},
  volume={130},
  pages={4--11},
  year={2020},
  publisher={Elsevier}
}
\end{document}